\documentclass{amsart}
\usepackage[T1]{fontenc}
\usepackage{textcomp}
\usepackage[utf8x]{inputenc}
\usepackage[swedish,english]{babel}
\usepackage{ucs}
\usepackage{graphicx}
\usepackage{mathrsfs}
\usepackage{rotating} 
\usepackage{tikz}
\usetikzlibrary{calc}
\usepackage{amssymb}
\usepackage[titletoc,toc]{appendix}
\usepackage[hidelinks]{hyperref}
\usepackage{todonotes}
\usepackage{pgfplots}
\usepackage[miktex]{gnuplottex}
\numberwithin{equation}{section} 
\newtheorem{thm}{Theorem}[section]

\newtheorem{cor}[thm]{Corollary}
\newtheorem{lem}[thm]{Lemma}
\newtheorem{prop}[thm]{Proposition}

\theoremstyle{definition}
\newtheorem{definition}[thm]{Definition}
\theoremstyle{remark}
\newtheorem{rem}[thm]{Remark}
\numberwithin{equation}{section}

\newcommand{\id}{\text{id}}

\newcommand{\Pn}[1]{\mathbb{P}^{#1}}

\newcommand{\Z}{\mathbb{Z}}

\newcommand{\derham}[2]{H^{#1} \! \left(#2\right)}

\newcommand{\Ht}[1]{\mathcal{H}_{#1}[2]}
\newcommand{\Frob}{F}
\newcommand{\jac}[1]{\mathrm{Jac}\! \left( #1 \right)}

\newcommand{\pglt}{\mathrm{PGL}(3)}
\newcommand{\gline}[2]{#1 #2}
\newcommand{\field}{K}
\newcommand{\Mgl}[2]{\mathcal{M}_{#1}[#2]}
\newcommand{\pic}[1]{\mathrm{Pic}\left( #1 \right)}
\newcommand{\symp}[1]{\mathrm{Sp}\!\left( #1 \right)}
\newcommand{\Mg}[1]{\mathcal{M}_{#1}}
\newcommand{\Hg}[1]{\mathcal{H}_{#1}}
\newcommand{\Q}{\mathcal{Q}}
\newcommand{\Pts}{\mathcal{P}^2_7}
\newcommand{\Fq}[1]{\mathbb{F}_{#1}}
\newcommand{\Fqbar}[1]{\overline{\mathbb{F}}_{#1}}
\newcommand{\Hetc}[2]{H^{#1}_{\text{\'et},c}(#2,\mathbb{Q}_l)}
\newcommand{\Qbtg}{\mathcal{Q}_{\mathrm{btg}}}
\newcommand{\Qbtgbar}{\mathcal{Q}_{\overline{\mathrm{btg}}}}

\newcommand{\poincareserre}[1]{PS_{#1}}
\makeindex
\begin{document}

\title[Equivariant Cohomology via Point Counts]
{Equivariant Cohomology of the Moduli Space of Genus Three Curves with Symplectic Level Two Structure via Point Counts}%
\author{Olof Bergvall}%
\address{Humboldt-Universität zu Berlin, Institut für Mathematik, 10099 Berlin, Germany}
\email{olof.bergvall@hu-berlin.de}

\begin{abstract}
  We make cohomological computations related to the moduli space of genus three curves with
  symplectic level two structure by means of counting points over finite fields.
  In particular, we determine the cohomology groups of the
  quartic locus as representations of the symmetric group
  on seven elements.
\end{abstract}

\maketitle

\section{Introduction}
\label{intro}
Let $n$ be a positive integer and let $C$ be a curve.
A level $n$ structure on $C$ is a choice of basis
for the $n$-torsion of the Jacobian of $C$.
The purpose of this paper is to study the cohomology
of the moduli space $\Mg{3}[2]$ of genus $3$ curves
with symplectic level $2$ structure. 

A genus $3$ curve which is not hyperelliptic is
embedded as a plane quartic via its canonical linear system.
The corresponding locus in $\Mg{3}[2]$ is called
the quartic locus and it is denoted $\Q[2]$.
A plane quartic with level $2$ structure is specified,
up to isomorphism, by an ordered septuple of points in general position
in $\Pn{2}$, up to the action of $\mathrm{PGL}(3)$.
This identification will be the basis for our investigation
of $\Q[2]$.

Our main focus will be on $\Q[2]$ but we will also
consider its complement in $\Mg{3}[2]$, i.e. the hyperelliptic locus $\Hg{3}[2]$.
In both cases, the computations will be carried out via
point counts over finite fields. By virtue of the Lefschetz
trace formula, such point counts give cohomological information
in the form of Euler characteristics. However, both
$\Q[2]$ and $\Hg{3}[2]$ satisfy certain strong purity conditions
which allow us to deduce Poincaré polynomials from these Euler
characteristics.

The group $\symp{6,\Z/2\Z}$
acts on $\Mg{3}[2]$ as well as on $\Q[2]$ and $\Hg{3}[2]$ by changing level structures. 
The cohomology groups thus become $\symp{6,\Z/2\Z}$-representations
and our computations will therefore be equivariant. However,
the action of \linebreak $\symp{6,\Z/2\Z}$ is rather subtle on $\Q[2]$
when $\Q[2]$ is identified with the space of septuples of points
in general position in $\Pn{2}$. On the other hand, the action
of the symmetric group $S_7$ on seven elements is very clear
and we will therefore restrict our attention to this subgroup.
The full action of $\symp{6,\Z/2\Z}$ is the topic of ongoing research.

The main results are presented in Table~\ref{Qcohtable} and Table~\ref{hypS7tab}
where we give the cohomology groups
of $\Q[2]$ and $\Hg{3}[2]$ as representations of $S_7$.

\subsection*{Acknowledgements}
This paper is based on parts of my thesis \cite{bergvallthesis} written 
at Stockholms Universitet. I would like to thank my advisors Carel
Faber and Jonas Bergström for their help as well as Dan Petersen
for useful comments on an an early version of this paper.

\section{Symplectic level structures}
\label{sympstrsec}
Let $\field$ be an algebraically closed field of characteristic zero
and let $C$ be a smooth and irreducible curve of genus $g$ over $\field$.
The $n$-torsion part $\jac{C}\![n]$ of the Jacobian of $C$ is isomorphic
to $\left( \Z/n\Z \right)^{2g}$ as an abelian group
and the Weil pairing is a nondegenerate and alternating bilinear form on
$\jac{C}\![n]$.

\pagebreak[2]
\begin{definition}
 A \emph{symplectic level $n$ structure} on a curve $C$ is an ordered basis
 \linebreak $(D_1, \ldots, D_{2g})$ of $\jac{C}\![n]$ such that the Weil
 pairing has basis
 \begin{equation*}
  \left( \begin{array}{cc}
  0 & I_g \\
  -I_g & 0
  \end{array}\right),
 \end{equation*}
 with respect to this basis. Here, $I_g$ denotes the $g \times g$ identity matrix.
\end{definition}
\pagebreak[2]

For more information about the Weil pairing and level structures,
see for example \cite{arbarellocgh} or \cite{griffithsharris}. Since we shall
only consider symplectic level structures we shall refer to symplectic level structures
simply as level structures.

A tuple $(C,D_1,\ldots, D_{2g})$ where $C$ is a smooth and irreducible curve and
\linebreak $(D_1, \ldots, D_{2g})$ is a level $n$ structure on $C$ is
called a curve with level $n$ structure. Let $(C',D'_1, \ldots,D'_{2g})$
be another curve with level $n$ structure.
An isomorphism of curves with level $n$ structures is an isomorphism
of curves $\phi:C \to C'$ such that $\phi^*(D'_i)=D_i$ for $i=1, \ldots, n$.
We denote the moduli space of genus $g$ curves with level $n$ structure
by $\Mgl{g}{n}$. We remark that we shall consider these moduli spaces
as coarse spaces and not as stacks. For $n \geq 3$, this remark is somewhat vacuous,
see \cite{harrismorrison}, but for $n=2$ this is not the case.
The group $\symp{2g,\Z/n\Z}$ acts on $\Mgl{g}{n}$
by changing level structures.

In the following we shall only be interested in level $2$ structures.
A concept closely related to level $2$ structures is that of
theta characteristics.

\pagebreak[2]
\begin{definition}
Let $C$ be a smooth and irreducible curve and let $K_C$ be its
canonical class. An element $\theta \in \pic{C}$ such that
$2 \theta = K_C$ is called a \emph{theta characteristic}.
We denote the set of theta characteristics of $C$ by $\Theta(C)$.
\end{definition}
\pagebreak[2]

Let $C$ be a curve of genus $g$. 
Given two theta characteristics $\theta_1$ and $\theta_2$ on 
$C$ we obtain an element $D \in \jac{C}\![2]$ by taking
the difference $\theta_1 -\theta_2$.
Conversely, given a theta characteristic $\theta$ and a $2$-torsion
element $D$ we obtain a new theta characteristic as $\theta'=\theta+D$.
More precisely we have that $\Theta(C)$ is a $\jac{C}\![2]$-torsor
and the set $\tilde{\Theta}(C)=\Theta(C) \cup \jac{C}\![2]$ is
a vector space of dimension $2g+1$ over the field $\Z/2\Z$ of two elements.

\pagebreak[2]
\begin{definition}
 An ordered basis $A=(\theta_1, \ldots, \theta_{2g+1})$ of theta characteristics
 of the vector space $\tilde{\Theta}(C)$ is called an
 \emph{ordered Aronhold basis} if the expression
 \begin{equation*}
 h^0(\theta) \quad \mathrm{mod} \, 2,
 \end{equation*}
 only depends on the number of elements in $A$ that is required to express
 $\theta$ for any theta characteristic $\theta$.
\end{definition}
\pagebreak[2]

\pagebreak[2]
\begin{prop}
\label{aronholdprop}
 Let $C$ be a smooth an irreducible curve.
 There is a bijection between the set of ordered Aronhold bases on
 $C$ and the set of level $2$ structures on $C$.
\end{prop}
\pagebreak[2]

For a proof of Proposition~\ref{aronholdprop} as well as a more thorough
treatment of theta characteristics and Aronhold bases we refer
to \cite{grossharris}.

Proposition~\ref{aronholdprop} provides a more geometric way
to think about level $2$ structures. In the case of a plane
quartic curve, which shall be the case of most importance to us,
we point out that each theta characteristic occurring in an
Aronhold basis is cut out by a bitangent line. Thus, in the case
of plane quartics one can think of ordered Aronhold bases 
as ordered sets of bitangents (although not every ordered set of bitangents
constitute an ordered Aronhold basis). 

\section{Plane quartics}
Let $\field$ be an algebraically closed field of characteristic zero
and let $C$ be a smooth and irreducible curve of genus $g$ over $\field$.
If $C$ is not hyperelliptic it is embedded into $\Pn{g-1}$
via its canonical linear system.
Thus, a non-hyperelliptic curve of genus $3$ is embedded into
$\Pn{2}$ and by the genus-degree formula we see see that
the degree of the image is $4$. We shall therefore refer to
the complement of the hyperelliptic locus in $\Mg{3}$ as the
\emph{quartic locus} and denote it by $\Q=\Mg{3} \setminus \Hg{3}$. 
Similarly, we denote the complement of the hyperelliptic locus in
$\Mgl{3}{2}$ by $\Q[2]$. Clearly, the action of $\symp{6,\Z/2\Z}$
on $\Mgl{3}{2}$ restricts to an action on $\Q[2]$.

The purpose of this section is to give an explicit, combinatorial description
of $\Q[2]$. This description will be in terms of points
in general position. Intuitively, a set of points in the
projective plane is in general position if there is no ``unexpected''
curve passing through all of them. In our case, this is made precise
in by the following definition.

\pagebreak[2]
\begin{definition}
Let $(P_1, \ldots, P_7)$ be a septuple of points in $\Pn{2}$.
We say that the septuple is in \emph{general position} if there is no
line passing through any three of the points and no conic passing
through any six of them. We denote the moduli space
of septuples of points in general position up to projective
equivalence by $\Pts$.
\end{definition}
\pagebreak[2]

Let $T=(P_1, \ldots, P_7)$ be a septuple of points in general position
in the projective plane and let $\mathcal{N}_T$ be the net of cubics
passing through $T$. If we let $F_0$, $F_1$ and $F_2$ be generators
for $\mathcal{N}_T$, then the equation
\begin{equation*}
 \mathrm{det}\left( \frac{\partial F_i}{\partial x_j} \right) = 0, \quad i,j = 0,1,2,
\end{equation*}
describes a plane sextic curve $S_T$ with double points precisely at
$P_1, \ldots, P_7$. By the genus-degree formula we see that $S_T$
has geometric genus $3$ and it turns out that its smooth model
is not hyperelliptic. Moreover, if we let $\rho: C_T \to S_T$
be a resolution of the singularities, then $D_i = \rho^{-1}(P_i)$ is a 
theta characteristic and $(D_1, \ldots, D_7)$ is an ordered Aronhold
basis.

\pagebreak[2]
\begin{thm}[van Geemen \cite{dolgachevortland}]
 Sending a septuple $T=(P_1, \ldots, P_7)$ of points in general position
 in the projective plane to $(C_T, D_1, \ldots, D_7)$ gives
 a $\symp{6,\Z/2\Z}$-equivariant isomorphism
 \begin{equation*}
  \Pts \to \Q[2].
 \end{equation*}
\end{thm}
\pagebreak[2]

It should be pointed out that while the action of $\symp{6,\Z/2\Z}$
is clear on $\Q[2]$ its action on $\Pts$ is much more subtle.
However, we can at least plainly see the symmetric group $S_7 \subset \symp{6,\Z/2\Z}$ act on $\Pts$ by permuting points.

\section{The Lefschetz trace formula}
We are interested in the spaces $\Mg{3}[2]$, $\Q[2]$ and $\Hg{3}[2]$
and in particular we want to know their cohomology.
The Lefschetz trace formula provides a way to
obtain cohomological information about a space via point
counts over finite fields.

Let $p$ be a prime number, let $n\geq 1$ be an integer and let $q=p^n$.
Also, let $\Fq{q}$ denote a finite field with $q$ elements, let
$\Fq{q^m}$ denote a degree $m$ extension of $\Fq{q}$ and
let $\Fqbar{q}$ denote an algebraic closure of $\Fq{q}$.
Let $X$ be a scheme defined over $\Fqbar{q}$ and let $\Frob$
denote its geometric Frobenius endomorphism induced from $\Fq{q}$.
Finally, let $l$ be another prime number, different from $p$, and
let $\Hetc{k}{X}$ denote the $k$'th compactly supported étale cohomology
group of $X$ with coefficients in $\mathbb{Q}_l$.

Let $\Gamma$ be a finite group of rational automorphisms of $X$. Then each cohomology group 
$\Hetc{k}{X}$ is a $\Gamma$-representation. The Lefschetz trace formula allows us to obtain information about these representations by counting the number of fixed points of 
$\Frob \sigma$ for different $\sigma \in \Gamma$.

\pagebreak[2]
\begin{thm}[Lefschetz trace formula]
\label{lefschetz}
 Let $X$ be a separated scheme of finite type over 
 $\Fqbar{q}$ with Frobenius endomorphism $\Frob$ and let $\sigma$
 be a rational automorphism of $X$ of finite order. Then
 \begin{equation*}
  |X^{\Frob \sigma}| = 
  \sum_{k \geq 0} (-1)^k \cdot \mathrm{Tr}\left(\Frob \sigma, \Hetc{k}{X} \right),
 \end{equation*}
 where $X^{\Frob \sigma}$ denotes the fixed point set of $\Frob \sigma$.
\end{thm}
\pagebreak[2]

For a proof, see \cite{deligne}, Rapport - Th\'eor\`eme 3.2.

\begin{rem}
  This theorem is usually only stated in terms of $\Frob$. 
  To get the above version one simply applies
  the ``usual'' theorem to the twist of $X$ by $\sigma$.
\end{rem}

\begin{rem}
 \label{conjrem}
  If $\Gamma$ is a finite group of rational automorphisms of $X$ and $\sigma \in \Gamma$,
  then $|X^{\Frob \sigma}|$ will only depend on the conjugacy class of $\sigma$.
 \end{rem}
 
 Let $R(\Gamma)$ denote the representation ring of $\Gamma$ and
 let the compactly supported $\Gamma$-equivariant Euler characteristic of $X$ be 
 defined as the virtual representation
 \begin{equation*}
  \mathrm{Eul}^{\Gamma}_{X,c} = \sum_{k \geq 0} (-1)^k \cdot \Hetc{k}{X} \in R(\Gamma).
 \end{equation*}
 By Theorem \ref{lefschetz} we may determine $\mathrm{Eul}^{\Gamma}_{X,c}$ by computing
 $|X^{\Frob \sigma}|$ for each $\sigma \in \Gamma$ and by Remark \ref{conjrem} it is enough to
 do so for one representative of each conjugacy class. This motivates the
 following definition.
 
 \pagebreak[2]
 \begin{definition}
  Let $X$ be a separated scheme of finite type over $\Fq{q}$ with 
  Frobenius endomorphism $\Frob$
  and let $\Gamma$ be a finite group of rational automorphisms of $X$.
  The determination of $|X^{\Frob \sigma}|$ for all $\sigma \in \Gamma$ 
  is then called a \emph{$\Gamma$-equivariant point count} of $X$ over $\Fq{q}$.
 \end{definition}
 \pagebreak[2]

 \section{Minimal purity}
  \label{finiteminimalpurity}
 Let $X$ be a scheme over the finite field $\Fq{q}$ and let
 $\Gamma$ be a group of rational automorphisms of $X$.
 We define the compactly supported $\Gamma$-equivariant Poincaré polynomial of $X$
 as
 \begin{equation*}
  P^{\Gamma}_{X,c}(t) = \sum_{k \geq 0} \Hetc{k}{X} \cdot t^k \in R(\Gamma)[t].
 \end{equation*}
 In the previous section we saw that equivariant point counts give 
 equivariant Euler characteristics. Poincaré polynomials contain more
 information and are therefore more desirable to obtain but are typically
 more complicated to compute. However, if $X$ satisfies a certain
 purity condition one can recover the Poincaré polynomial from the Euler characteristic.
 See also \cite{bergstrom2}, \cite{bergstrom}, \cite{vandenbogaartedixhoven}
 and \cite{bergstromtommasi} where similar phenomena for compact spaces have been exploited.
 
 \pagebreak[2]
 \begin{definition}[Dimca and Lehrer \cite{dimcalehrer}]
 \label{poscharmp}
  Let $X$ be an irreducible and separated scheme of finite type over $\Fqbar{q}$ with Frobenius 
  endomorphism $\Frob$ and let $l$ be a prime not dividing $q$. The scheme $X$ is called
  \emph{minimally pure} if $\Frob$ acts on $\Hetc{k}{X}$ with all eigenvalues equal to 
  $q^{k-\mathrm{dim}(X)}$.
  
  A pure dimensional and separated scheme $X$ of finite type over $\Fqbar{q}$ is minimally
  pure if for any collection $\{X_1, \ldots, X_r\}$ of irreducible components of 
  $X$, the irreducible scheme $X_1 \setminus (X_2 \cup \cdots \cup X_r)$ is
  minimally pure.
 \end{definition}
 \pagebreak[2] 
 
 Thus, if $X$ is minimally pure, then a term $q^{k-\mathrm{dim}(X)}$ in $|X^{\Frob}|$
 can only come from $\Hetc{k}{X}$ and we can determine the $\Gamma$-equivariant Poincaré
 polynomial of $X$ via the relation
 \begin{equation*}
  \mathrm{Eul}^{\Gamma}_X(\sigma) = q^{-2\mathrm{dim}(X)} \cdot P^{\Gamma}_X(\sigma)(-q^2).
 \end{equation*}
 We will see that the moduli space $\Q[2]$ is minimally pure.
 
 Let $C \subset \mathbb{P}^2$ be a plane quartic, let $P \in C$ be a point 
 and let $T_PC$ denote the tangent line of $C$ at $P$. 
 We say that $P$ is a \emph{bitangent point} if
 \begin{equation*}
 C \cdot T_PC = 2P + 2Q
 \end{equation*}
 for some point $Q$ that might coincide with $P$. 
 If $P \neq Q$ we say that $P$ is a \emph{genuine bitangent point}.
 We denote the moduli space of plane quartics with level $2$ structure 
 marked with a bitangent point by $\Qbtgbar[2]$ and we 
 denote the moduli space of plane quartics with level $2$ structure 
 marked with a genuine bitangent point by $\Qbtg[2]$.
 The space $\Qbtg[2]$ is an open subvariety of $\Qbtgbar[2]$.
 
 \pagebreak[2]
 \begin{lem}
 \label{btgpuritylemma}
  $\Qbtgbar[2]$ is minimally pure.
 \end{lem} 
 \pagebreak[2] 
 
 \begin{proof}
 Looijenga \cite{looijenga} has shown that $\Qbtg[2]$ is isomorphic
 to a finite disjoint union of varieties, each isomorphic to the complement
 of an arrangement of tori in an ambient torus of dimension $6$. 
 Dimca and Lehrer \cite{dimcalehrer}
 has shown that such complements of arrangements are minimally pure
 and it thus follows that $\Qbtg[2]$ is minimally pure.
 On the other hand, Looijenga \cite{looijenga} has shown that there
 is an injection
 \begin{equation*}
 \Hetc{k}{\Qbtgbar[2]} \hookrightarrow \Hetc{k}{\Qbtg[2]}.
 \end{equation*}
 We thus see that $\Hetc{k}{\Qbtgbar[2]}$ is an $F$-invariant
 subspace of $\Hetc{k}{\Qbtg[2]}$. Since the eigenvalues of
 $F$ on $\Hetc{k}{\Qbtg[2]}$ are all equal to $q^{k-6}$ we conclude that
 the same is true for $\Hetc{k}{\Qbtgbar[2]}$. Hence,
 $\Qbtgbar[2]$ is minimally pure.
 \end{proof}
 
 \pagebreak[2]
 \begin{prop}
 \label{q2purityprop}
 $\Q[2]$ is minimally pure.
 \end{prop}
 \pagebreak[2]
 
 \begin{proof}
  A plane quartic has $28$ bitangents so the morphism
  \begin{equation*}
   \pi: \Qbtgbar[2] \to \Q[2],
  \end{equation*}
  forgetting the marked bitangent point, is finite of degree $2 \cdot 28 = 56$.
  Thus, the map
  \begin{equation*}
   \pi_* \circ \pi^*: \Hetc{k}{\Q[2]} \to \Hetc{k}{\Q[2]}
  \end{equation*}
  is multiplication with $\mathrm{deg}(\pi) = 56$. In particular,
  the map
  \begin{equation*}
   \pi: \Hetc{k}{\Q[2]} \to \Hetc{k}{\Qbtgbar[2]}
  \end{equation*}
  is injective and we may thus conclude that $\Q[2]$ is minimally pure
  as in the proof of Lemma~\ref{btgpuritylemma}.
 \end{proof}
 
 Since $\Q[2]$ is isomorphic to $\Pts$, we may compute the cohomology
 of $\Q[2]$ as a representation of $S_7$ by making $S_7$-equivariant
 point counts of $\Pts$.
 
  \section{Equivariant point counts}
 In this section we shall perform a $S_7$-equivariant point count of $\Pts$. This
 amounts to the computation of $\left|\left( \Pts \right)^{\Frob \sigma} \right|$
 for one representative $\sigma$ of each of the fifteen conjugacy classes of $S_7$.
 The computations will be rather different in the various cases but at least the underlying
 idea will be the same. Throughout this section we shall work over a finite field
 $\Fq{q}$ where $q$ is odd.
 
 Let $U$ be a subset of $\left( \Pn{2} \right)^7$
 and interpret each point of 
 $U$
 as an ordered septuple of points in $\Pn{2}$.
  Define the \emph{discriminant locus} 
 $\Delta \subset U$
 as the subset consisting of septuples which are not in general position. 
 If $U$ contains the subset of $\left( \Pn{2} \right)^7$ consisting
 of all septuples which are in general position, then
 \begin{equation*}
  \Pts = \left( U \setminus \Delta \right)/\pglt.
 \end{equation*}
 An element of $\pglt$ is completely specified by where
 it takes four points in general position. Therefore, 
 the points of $\Pts$ do not have any automorphisms
 and we have the simple relation
 \begin{equation}
 \label{quotientformula}
 \left|\left( \Pts \right)^{\Frob \sigma} \right| =
 \frac{\left| U^{\Frob \sigma} \right| - \left| \Delta^{\Frob \sigma} \right| }{|\pglt|}.
 \end{equation}
 We will choose the set $U$ in such a way that counting fixed points
 of $\Frob \sigma$ in $U$ is easy. We shall therefore focus on 
 the discriminant locus.
 
 The discriminant locus can be decomposed as
 \begin{equation*}
  \Delta = \Delta_l \cup \Delta_c,
 \end{equation*}
 where $\Delta_l$ consists of septuples where at least three points lie on a line
 and $\Delta_c$ consists of septuples where at least six points lie on a conic.
 The computation of $|\Delta^{\Frob \sigma}|$ will consist of the following
 three steps:
 \begin{itemize}
  \item the computation of $|\Delta_l^{\Frob \sigma}|$,
  \item the computation of $|\Delta_c^{\Frob \sigma}|$,
  \item the computation of $|(\Delta_l \cap \Delta_c)^{\Frob \sigma}|$.
 \end{itemize}
 We can then easily determine $|\Delta^{\Frob \sigma}|$ via the principle of inclusion
 and exclusion.
 
\begin{figure}[ht]
\centering
\resizebox{0.8\linewidth}{!}{
 \begin{tikzpicture}
  \begin{axis}[axis lines=none]
    \addplot +[id=intro_con,black,no markers, raw gnuplot, thick, empty line = jump 
      ] 
      gnuplot {
      set contour base;
      set cntrparam levels discrete 0.003;
      unset surface;
      set view map;
      set isosamples 500;
      set xrange [-3:3];
      set yrange [-3:3];
      splot 3*x^2+3*y^2+2*x*y-4;
    };
    \addplot[thick,white,dashed,domain=-2:2]{-2};
    \addplot[thick,white,dashed,domain=-2:2]{2};
    \addplot[thick,black,dashed,domain=-1.7:1.7]{-1.23};
    \addplot +[mark=none,thick,black,dashed] coordinates {(-1.23, -1.7) (-1.23, 1.7)};
    \addplot[thick,black,dashed,domain=-0.28:1.7]{-x+1.42};
    \node [circle, black, fill, minimum size= 2pt, inner sep=0, label={[black]above right:$P$}] at (77, 77) {};
    \node [circle, black, fill, minimum size= 2pt, inner sep=0, label={[black]above right:$Q$}] at (174, 174) {};
    \node [circle, black, fill, minimum size= 2pt, inner sep=0, label={[black]above right:$R$}] at (271, 271) {};
  \end{axis}
\end{tikzpicture}
}
\caption{A conic $C$ with a point $P$ on the outside of $C$, a point $Q$ on the inside of $C$ and a point $R$ on $C$.}
\label{insideoutsidefig}
\end{figure}
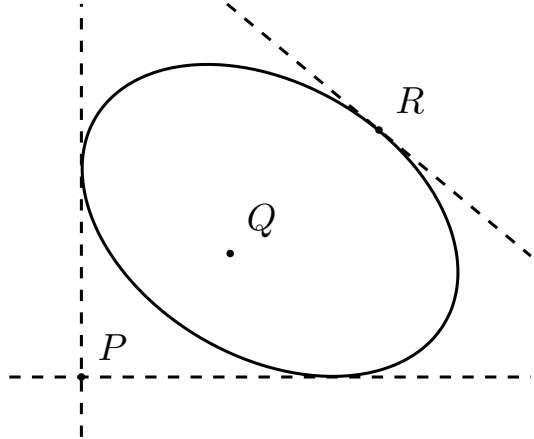

 In the analysis of $\Delta_l \cap \Delta_c$ the following definition will sometimes
 be useful, see Figure~\ref{insideoutsidefig}.
 
 \pagebreak[2]
 \begin{definition}
  Let $C$ be a smooth conic over $\Fq{q}$ and let $P \in \Pn{2}$ be a $\Fq{q}$-point. We
  then say
  \begin{itemize}
   \item that $P$ is on the \emph{$\Fq{q}$-inside} of $C$ if there is no $\Fq{q}$-tangent
   to $C$ passing through $P$,
   \item that $P$ is on $C$ if there is precisely one $\Fq{q}$-tangent to $C$ passing
   through $P$,
   \item that $P$ is on the \emph{$\Fq{q}$-outside} of $C$ if there are two $\Fq{q}$-tangents
   to $C$ passing through $P$.
  \end{itemize}
 \end{definition}
 \pagebreak[2] 
 
 For a motivation of the terminology, see Figure~\ref{insideoutsidefig}.
 
 It is not entirely clear that the above definition makes sense. To see that it
 in fact does, we need the following lemma.
 
 \pagebreak[2]
 \begin{lem}
\label{threetan}
 Let $C \subset \Pn{2}$ be a smooth conic over a field $k$. If there is a point
 $P$ such that three tangents of $C$ pass through $P$, then
 the characteristic of $k$ is $2$.
\end{lem}
 \pagebreak[2] 
 
\begin{proof}
 Since $C$ is smooth, 
 the three points of tangency will be in general position, so by a projective
 change of coordinates they can be transformed to $[1:0:0]$, $[0:1:0]$
 and $[0:0:1]$ and $C$ will then be given by a polynomial $F=XY+\alpha XZ+ \beta YZ$,
 where $\alpha, \beta \in k^*$. The tangent lines thus become
 $Y+\alpha Z$, $X+\beta Z$ and $\alpha X+\beta Y$. 
 
 Let the coordinates of $P$ be $[a:b:c]$. Since these lines all pass through
 $P$, the first tangent equation gives that $b=-\alpha c$ and the second gives
 $a = - \beta c$. Inserting these expressions into the third tangent equation
 gives $-2 \alpha \beta c =0$. If $c=0$, then also $a=b=0$ which is impossible.
 Since also $\alpha$ and $\beta$ are nonzero we see that the only possibility is
 that the characteristic of $k$ is $2$. 
\end{proof}
 
 Let $\sigma^{-1}=(i_1\ldots i_r)$ be a cycle in $S_7$. An ordered septuple $(P_1, \ldots, P_7)$ of points in $\Pn{2}$
 will be fixed by $\Frob\sigma$ if and only if $\Frob P_{i_j}=P_{i_{j+1}}$ for $i=1, \ldots,r-1$ and $\Frob P_{i_r}=P_{i_1}$.
 This is the motivation for the following definition.
 
 \pagebreak[2]
 \begin{definition}
  Let $X$ be a $\Fq{q}$-scheme with Frobenius endomorphism $\Frob$ and let
  $Z \subset X_{\Fqbar{q}}$ be a subscheme. If 
  \begin{equation*}
  \left|\left\{ \Frob^iZ \right\}_{i \geq 0}\right| = m,
  \end{equation*}  
  we say that $Z$ is a \emph{strict $\Fq{q^m}$}-subscheme.
 
  If $Z$ is a strict $\Fq{q^m}$-subscheme, the $m$-tuple $(Z, \ldots, \Frob^{m-1}Z)$
  is called a \emph{conjugate $m$-tuple}.
  Let $r$ be a positive integer and let $\lambda = [1^{\lambda_1}, \ldots, r^{\lambda_r}]$
  be a partition of $r$. An $r$-tuple $(Z_1, \ldots, Z_r)$
  of closed subschemes of $X$ is called a \emph{conjugate $\lambda$-tuple} if
  it consists of $\lambda_1$ conjugate $1$-tuples, $\lambda_2$-conjugate $2$-tuples
  and so on.
  We denote the set of conjugate $\lambda$-tuples of closed points of $X$ by $X(\lambda)$.
 \end{definition}
 \pagebreak[2]
  
 Since the conjugacy class of an element in $S_7$ is given by its cycle type, 
 we want to count the number of conjugate $\lambda$-tuples in both $U$ and $\Delta$
 for each partition of seven. In this pursuit, the following formula is helpful.
 Its proof is a simple
 application of the principle of inclusion and exclusion.
 
 \pagebreak[2]
 \begin{lem}
 \label{conjlem}
  Let $X$ be a $\Fq{q}$-scheme and let $\lambda=[1^{\lambda_1}, \ldots, n^{\lambda_n}]$ be a partition.
  Then
  \begin{equation*}
 |X(\lambda)| = \prod_{i=1}^{\nu} \prod_{j=0}^{\lambda_i-1} \left[ \left( \sum_{d \mid i} \mu\left( \frac{i}{d} \right) \cdot |X(d)| \right)  - i \cdot j \right],
\end{equation*}
where $\mu$ is the Möbius function. 
 \end{lem}
 \pagebreak[2] 
 
 We now recall a number of basic results regarding point counts.
 First, note that if we apply Lemma \ref{conjlem} to $X= \left( \Pn{n} \right)^{\vee}$, the dual projective space, we see that the number of conjugate $\lambda$-tuples of hyperplanes
 is equal to the number of conjugate $\lambda$-tuples of points in $\Pn{n}$.
 We also recall that
\begin{equation*}
\left| \Pn{n}_{\Fq{q}} \right| = \sum_{i=0}^n q^i,
\end{equation*} 
 and that 
\begin{equation*}
|\pglt|=q^3 \cdot (q^3-1) \cdot (q^2-1).
\end{equation*}
A slightly less elementary result is that the number of smooth 
conics defined over $\Fq{q}$ is 
\begin{equation*}
q^5-q^2.
\end{equation*}
 To see this,
 note that there is a $\Pn{5}$ of conics. Of these there are $q^2+q+1$ double $\Fq{q}$-lines, $\frac{1}{2} \cdot (q^2+q+1) \cdot (q^2+q)$
 intersecting pairs of $\Fq{q}$-lines and $\frac{1}{2} \cdot (q^4-q)$ conjugate pairs of $\Fq{q^2}$-lines.

 We are now ready for the task of counting the number of conjugate $\lambda$-tuples for each element of $S_7$.
 
 \begin{rem}
 Since $\Pts$ is minimally pure, Equation \ref{quotientformula} gives that $\left|\left( \Pts \right)^{\Frob \sigma} \right|$
 is a monic polynomial in $q$ of degree six so it is in fact enough to make counts for
 six different finite fields and interpolate. This is however hard to carry out in practice, even with a computer,
 as soon as $\lambda$ contains parts of large enough size (where ``large enough'' means $3$ or $4$).
 However, one can always obtain partial information which provides important checks for our computations.
\end{rem}

 \subsection{The case \texorpdfstring{$\lambda = [7]$}{[7]}}
 Let $\lambda=[7]$. Since we only need to make the computation for one permutation
 $\sigma$ of cycle type $\lambda$, we may as well assume that $\sigma^{-1}=(1234567)$
 so that $\Frob$ acts as $\Frob P_i=P_{i+1}$ for $i=1,\ldots,6$ and $\Frob P_7=P_1$.
 In this case, we simply take $\left( \Pn{2} \right)^7$ as our set $U$.

 The main observation is the following.
 
 \pagebreak[2]
 \begin{lem}
 \label{7lines}
 If $(P_1, \ldots, P_7)$ is a $\lambda$-tuple with three of its points on a
 line, then all seven points lie on a line defined over $\Fq{q}$.
 \end{lem}
 \pagebreak[2] 
 
 \begin{proof}
  Suppose that the set $S=\{P_i,P_j,P_k\}$ is contained in the line $L$.
  Then $L$ is either defined over $\Fq{q}$ or $\Fq{q^7}$.
  One easily checks that for each of the $\binom{7}{3}=35$ possible choices
  of $S$ there is an integer $1 \leq i \leq 6$ such that $|\Frob^iS \cap S|=2$.
  Since a line is defined by any two points on it we have that $L=F^iL$.
  Hence, we have that $L$ is defined over $\Fq{q}$ and that
  $\{P_i,\Frob P_i, \ldots, \Frob^6P_i\}=\{P_1, \ldots, P_7\} \subset L$.
 \end{proof}
  
  \pagebreak[2]
  \begin{lem}
  \label{7conics}
  If $(P_1, \ldots, P_7)$ is a $\lambda$-tuple with six of its points on a
  smooth conic, then all seven points lie on a smooth conic defined over $\Fq{q}$.
  \end{lem}
 \pagebreak[2]  
  
  \begin{proof}
   Suppose that the set $S=\{P_{i_1}, \ldots, P_{i_6}\}$ lies on a smooth conic $C$.
   We have $|\Frob S \cap S|=5$ and since a conic is defined by any five points
   on it we have $\Frob C=C$. Hence, we have that $C$ is defined over $\Fq{q}$ and that
   all seven points lie on $C$.
  \end{proof}
  
  We conclude that $\Delta_l$ and $\Delta_c$ are disjoint.
  We obtain $|\Delta_l|$ by first choosing a $\Fq{q}$-line $L$ and then
  picking a $\lambda$-tuple on $L$. We thus have
  \begin{equation*}
   |\Delta_l| = (q^2+q+1) \cdot (q^7-q).
  \end{equation*}
  To obtain $|\Delta_c|$ we first choose a smooth conic $C$ and then
  a conjugate $\lambda$-tuple on $C$. We thus have
  \begin{equation*}
   |\Delta_c| = (q^5-q^2) \cdot (q^7-q).
  \end{equation*}
  Equation \ref{quotientformula} now gives
  \begin{equation*}
   \left|\left( \Pts \right)^{\Frob \sigma} \right| = q^6+q^3.
  \end{equation*}

  \subsection{The case \texorpdfstring{$\lambda = [1,6]$}{[1,6]}}
 Let $\lambda=[1,6]$. Since we only need to make the computation for one permutation
 $\sigma$ of cycle type $\lambda$, we may as well assume that $\sigma^{-1}=(123456)(7)$
 so that $\Frob$ acts as $\Frob P_i=P_{i+1}$ for $i=1,\ldots,5$, $\Frob P_6=P_1$ and 
 $\Frob P_7=P_7$. Also in this case we take $\left( \Pn{2} \right)^7$ as our set $U$.

 The main observation is the following.
 
 \pagebreak[2]
 \begin{lem}
  \label{61lines}
  If a $\lambda$-tuple has three points on a line, then either
 \begin{itemize}
  \item[(1)] the first six points of the $\lambda$-tuple lie on a $\Fq{q}$-line or,
  \item[(2)] the first six points lie on two conjugate $\Fq{q^2}$-lines, 
  the $\Fq{q^2}$-lines contain three $\Fq{q^6}$-points each and these triples
  are interchanged by $\Frob$, or,
  \item[(3)] the first six points lie pairwise on three conjugate $\Fq{q^3}$-lines
  which intersect in $P_7$.
 \end{itemize}
 \end{lem}
 \pagebreak[2]
 
 \begin{proof}
  Suppose that $S=\{P_i, P_j, P_k\}$ lie on a line $L$.
  Then $L$ is either defined over $\Fq{q}$, $\Fq{q^2}$, $\Fq{q^3}$ 
  or $\Fq{q^6}$.
  One easily checks that for each of the $\binom{7}{3}=35$ possible choices
  of $S$ there is an integer $1 \leq i \leq 3$ such that $|F^iS \cap S|=2$
  so $L$ is defined over $\Fq{q}$, $\Fq{q^2}$ or $\Fq{q^3}$,
  i.e. we are in one of the three cases above.
 \end{proof}
 
 Let $\Delta_{l,i}$ be the subset of $\Delta_l$ corresponding
 to case (i) in Lemma \ref{61lines}. The set $\Delta_{l,1}$ is clearly disjoint
 from the other two.
 
 \pagebreak[2]
 \begin{lem}
 \label{61cons}
 If six of the points of a $\lambda$-tuple lie on a smooth conic,
 then all of the first six points of the tuple lie on the conic and the conic is defined
 over $\Fq{q}$.
 \end{lem}
 \pagebreak[2] 
 
 \begin{proof}
 Suppose that $S=\{P_{i_1}, \ldots, P_{i_6}\}$ lie on a smooth conic $C$.
 Then $|\Frob S \cap S| \geq 5$ so $\Frob C=C$. Let $P \in S$ be a $\Fq{q^6}$-point.
 Then $\{P,FP, \ldots, F^5\}=\{P_1, \ldots, P_6\} \subset C$.
 \end{proof} 
 
 Since a smooth conic does not contain a line, we have that $\Delta_c$
 only intersects $\Delta_{l,3}$. 
 
 We compute $|\Delta_{l,1}|$ by first choosing a $\Fq{q}$-line $L$
 and then a $\Fq{q^6}$ point on $L$. Finally we choose
 a $\Fq{q}$-point $P_7$ anywhere. We thus have
 \begin{equation*}
 |\Delta_{l,1}| = (q^2+q+1) \cdot (q^6-q^3-q^2+q) \cdot (q^2+q+1).
 \end{equation*}
 To obtain $|\Delta_{l,2}|$ we first choose a $\Fq{q^2}$-line, $L$. By Lemma \ref{conjlem}
 there are $q^4-q$ such lines. The other $\Fq{q^2}$-line must then be $\Frob L$. 
 We then choose a $\Fq{q^6}$-point $P_1$ on
$L$. The points $P_2=\Frob P_1, \ldots, P_6=\Frob^5 P_1$ will then be the rest of our
conjugate sextuple. By Lemma \ref{conjlem} (with $\Fq{q^2}$ as the ground field) 
there are $q^6-q^2$ choices. We now have two $\Fq{q^2}$-lines with three of our
six $\Fq{q^6}$-points on each so all that remains is to choose a $\Fq{q}$-point
anywhere we want in $q^2+q+1$ ways. Hence,
\begin{equation*}
 |\Delta_{l,2}|=(q^4-q) (q^6-q^2) (q^2+q+1).
\end{equation*}
To count $|\Delta_{l,3}|$ we first choose a $\Fq{q}$-point $P_7$ in $q^2+q+1$ ways.
There is a $\Pn{1}$ of lines through $P_7$ and we want to choose
a $\Fq{q^3}$-line $L$ through $P$. By Lemma \ref{conjlem} there are
$q^3-q$ choices. Finally, we choose a $\Fq{q^6}$-point $P_1$ on $L$.
By Lemma \ref{conjlem} there are $q^6-q^3$ possible choices. We
thus have
\begin{equation*}
 |\Delta_{l,3}|=(q^2+q+1) (q^3-q) (q^6-q^3).
\end{equation*}
In order to finish the computation of $\Delta_l$, we need to compute
$|\Delta_{l,2} \cap \Delta_{l,3}|$.  
We first choose a pair of conjugate $\Fq{q^2}$-lines in $\frac{1}{2}(q^4-q)$
ways. These intersect in a $\Fq{q}$-point and we choose $P_7$ away from this point in $q^2+q$ ways.
We then choose a $\Fq{q^3}$-line through $P_7$ in $q^3-q$ ways. This line intersects 
the two $\Fq{q^2}$-lines
in $2$ distinct points which clearly must be defined over $\Fq{q^6}$. 
We choose one of them to become $P_1$ in $2$ ways. Thus, in total we have
\begin{equation*}
  |\Delta_{l,2} \cap \Delta_{l,3}| = (q^4-q) \cdot (q^2+q) \cdot (q^3-q).
\end{equation*}

To compute $|\Delta_c|$ we first choose a smooth
conic $C$ in $q^5-q^2$ ways and then use Lemma \ref{conjlem}
to see that we have $q^6-q^3-q^2+q$ ways of choosing a conjugate sextuple
on $C$. Finally, we choose $P_7$ anywhere we want in $q^2+q+1$ ways.
We thus see that
\begin{equation*}
 |\Delta_c| = (q^5-q^2) (q^6-q^3-q^2+q) (q^2+q+1).
\end{equation*}

It remains to compute the size of the intersection between
$\Delta_l$ and $\Delta_c$.
To compute do this,
we begin by choosing a smooth conic $C$ in $q^5-q^2$ ways and then
a $\Fq{q}$-point $P_7$ not on $C$ in $q^2+q+1-(q+1)=q^2$
ways. By Lemma \ref{conjlem} there are $q^3-q$ strict $\Fq{q^3}$-lines passing
through $P$. All of these intersect $C$ in two $\Fq{q^3}$-points
since, by Lemma \ref{threetan}, these lines cannot be tangent to $C$ since the 
characteristic of $\Fq{q}$ is odd. More precisely, choosing any of the
$q^3-q$ strict $\Fq{q^3}$-points of $C$ gives a strict $\Fq{q^3}$-line, and since every such line
cuts $C$ in exactly two points we conclude that there are precisely
$\frac{1}{2}(q^3-q)$ strict $\Fq{q^3}$-lines through $P$ intersecting $C$
in two $\Fq{q^3}$-points. Thus, the remaining 
\begin{equation*}
 q^3-q-\frac{1}{2}(q^3-q) = \frac{1}{2}(q^3-q)
\end{equation*}
$\Fq{q^3}$-lines through $P_7$ will intersect $C$ in two $\Fq{q^6}$-points. 
If we pick one of them and label it $P_1$ we obtain an element in $\Delta_l \cap \Delta_c$.
Hence,
\begin{equation*}
 |\Delta_l \cap \Delta_c| = (q^5-q^2) q^2 (q^3-q).
\end{equation*}
We now conclude that
\begin{equation*}
 \left|\left( \Pts \right)^{\Frob \sigma} \right| = q^6-2q^3+1.
\end{equation*}

  \subsection{The case \texorpdfstring{$\lambda=[2,5]$}{[2,5]}}
 Throughout this section, $\lambda$ will denote the partition
 $[2,5]$. We take $U=\left( \Pn{2} \right)^7$.
 
 \pagebreak[2]
\begin{lem}
\label{52linesandcons}
 If $(P_1, \ldots, P_7)$ is a $\lambda$-tuple with three of its points on a
 line, then all five $\Fq{q^5}$-points lie on a line defined over $\Fq{q}$.
 If six of the points lie on a smooth conic $C$, then all seven points lie on 
 $C$ and $C$ is defined over $\Fq{q}$. 
\end{lem}
\pagebreak[2]

\begin{proof}
The proof is very similar to the proofs of Lemmas \ref{7lines} and \ref{7conics}
and is therefore omitted.
\end{proof}

By Lemma \ref{conjlem}, there are $q^{10}+q^5-q^2-q$ 
conjugate quintuples whereof $(q^2+q+1)(q^5-q)$ lie on
a line. We may thus choose a conjugate quintuple whose points
do not lie on a line in $q^{10}-q^7-q^6+q^3$ ways.
This quintuple defines a smooth conic $C$.
By Lemma \ref{52linesandcons},
it is enough to choose  a conjugate pair outside $C$ in order
to obtain an element of $\left( \Pn{2} \right)^7 \setminus \Delta$
of the desired type.
Since there
are $q^4-q$ conjugate pairs of which $q^2-q$ lie on $C$ there are $q^4-q^2$ 
remaining choices. We thus obtain
\begin{equation*}
 \left|\left( \Pts \right)^{\Frob \sigma} \right| = q^6-q^2.
\end{equation*}

\subsection{The case \texorpdfstring{$\lambda=[1^2,5]$}{[1,1,5]}}
The computation in this case is very similar to that of
the case $\lambda=[2,5]$ and we therefore simply state the
result:
\begin{equation*}
 \left|\left( \Pts \right)^{\Frob \sigma} \right| = q^6-q^2.
\end{equation*}

\subsection{The case \texorpdfstring{$\lambda=[3^1,4^1]$}{[3,4]}}
Throughout this section, $\lambda$ shall mean the partition
$[3^1,4^1]$. Since we only need to make the computation for one permutation,
we shall assume that the Frobenius permutes points $P_1, P_2, P_3, _4$ according
to $(1234)$ and the three points $P_5,P_6,P_7$ according to  $(567)$.
We take $U=\left( \Pn{2} \right)^7$.

\pagebreak[2]
\begin{lem}
\label{43lines}
 If a conjugate $\lambda$-tuple has three points on a line, then either
 \begin{itemize}
  \item[(1)] the four $\Fq{q^4}$-points lie on a $\Fq{q}$-line, or
  \item[(2)] the three $\Fq{q^3}$-points lie on a $\Fq{q}$-line.
 \end{itemize}
\end{lem}
\pagebreak[2]

\begin{proof}
 It is easy to see that if three $\Fq{q^4}$-points lie on a line, then
 all four $\Fq{q^4}$-points lie on that line and even easier to see the
 corresponding result for three $\Fq{q^3}$-points.
 
 Suppose that two $\Fq{q^4}$-points $P_i$ and $P_j$ and a $\Fq{q^3}$-point
 $P$ lie on a line $L$. Since $\Frob^4P_i=P_i$ and $\Frob^4P_j=P_j$
 we see that $\Frob F^4L=L$. However, $\Frob^4P=\Frob P \neq P$. Repeating
 this argument again, with $\Frob P$ in the place of $P$, shows
 that also $\Frob ^2P$ lies on $L$. We are thus in case (1).
 
 If we assume that two $\Fq{q^3}$-points and a $\Fq{q^4}$-point lie on a line,
 then a completely analogous argument shows that all four $\Fq{q^4}$-points 
 lie on that line.
\end{proof}

We decompose $\Delta_l$ as
\begin{equation*}
 \Delta_l = \Delta_{l,1} \cup \Delta_{l,2},
\end{equation*}
where $\Delta_{l,1}$ consists of tuples with the four $\Fq{q^4}$-points on
a line and $\Delta_{l,2}$ consists of tuples with the three $\Fq{q^3}$-points on
a line. The computations of $|\Delta_{l,1}|$, $\Delta_{l,2}$ and $|\Delta_{l,1} \cap \Delta_{l,2}|$
are completely straightforward and we get 
\begin{equation*}
 |\Delta_l| = q^{13} + 2q^{12} - 3q^{10} - 2q^9 + q^8 + q^7 - q^6 - q^5 + q^4 + q^3.
\end{equation*}

To compute $|\Delta_c|$ we start by noting that if
six of the points of a $\lambda$-tuple lie on a smooth conic $C$,
then all seven points lie on $C$ and $C$ is defined over $\Fq{q}$.
Thus, the problem consists of choosing a smooth conic $C$
over $\Fq{q}$ and then picking a $\lambda$-tuple on $C$.
We thus have
\begin{equation*}
 |\Delta_c| = (q^5-q^2)(q^4-q^2)(q^3-q).
\end{equation*}
Since no three points on a smooth conic lie on a line we conclude
that the intersection $\Delta_{l} \cap \Delta_{c}$
is empty.
We now obtain
\begin{equation*}
 \left|\left( \Pts \right)^{\Frob \sigma} \right| = q^6-q^5-2q^4+q^3+q^2:
\end{equation*}

\subsection{The case \texorpdfstring{$\lambda=[1,2,4]$}{[1,2,4]}}
\label{421sec}
Throughout this section, $\lambda$ shall mean the partition
$[1^1,2^1,4^1]$. Since we only need to make the computation for one permutation,
we shall assume that the Frobenius permutes points $P_1, P_2,P_3, P_4$ according
to $(1234)$, switches the two points $P_5,P_6$ and fixes the point $P_7$.
The computation will turn out to be quite a bit more complicated in
this case than in the previous cases,
mainly because both $1$ and $2$ divide
$4$. We take $U=\left( \Pn{2} \right)^7$.

We have the following trivial decomposition of $\Delta_{l}$
\begin{equation*}
 \Delta_l = \bigcup_{i=1}^6 \Delta_{l,i},
\end{equation*}
where
 \begin{itemize}
  \item $\Delta_{l,1}$ consists of $\lambda$-tuples with three $\Fq{q^4}$-points lying on a line, 
  \item $\Delta_{l,2}$ consists of $\lambda$-tuples with two $\Fq{q^4}$-points and a $\Fq{q^2}$-point lying on a line,
  \item $\Delta_{l,3}$ consists of $\lambda$-tuples with two $\Fq{q^4}$-points and the $\Fq{q}$-point lying on a line, 
  \item $\Delta_{l,4}$ consists of $\lambda$-tuples with a $\Fq{q^4}$-point and two $\Fq{q^2}$-points lying on a line, 
  \item $\Delta_{l,5}$ consists of $\lambda$-tuples with a $\Fq{q^4}$-point, a $\Fq{q^2}$-point and a $\Fq{q}$-point lying on a line, and, 
  \item $\Delta_{l,6}$ consists of $\lambda$-tuples with two $\Fq{q^2}$-points and the $\Fq{q}$-point lying on a line.
 \end{itemize}
This decomposition is of course naive and is not
very nice to work with since none of the possible intersections
are empty. The reader can surely think of many other decompositions
which a priori look more promising. However, the more ``clever''
approaches we have tried have turned out to be quite hard to work with in practice.
The positive thing about the above decomposition is that
most intersections are rather easily handled
and that quadruple intersections (and higher) all consist of tuples where
all seven points lie on a $\Fq{q}$-line.

The two slightly more complicated sets in the above list are $\Delta_{l,2}$ 
and $\Delta_{l,3}$. We shall therefore comment a bit about the computations
involving them.

The set $\Delta_{l,2}$ splits into three disjoint subsets
\begin{equation*}
 \Delta_{l,2} = \Delta_{l,2}^1 \cup \Delta_{l,2}^2 \cup \Delta_{l,2}^3.
\end{equation*}
where
\begin{itemize}
 \item $\Delta_{l,2}^1$ consists of $\lambda$-tuples such that the four $\Fq{q^4}$-points and the two $\Fq{q^2}$-points lie on a $\Fq{q}$-line, or,
 \item $\Delta_{l,2}^2$ consists of $\lambda$-tuples such that the two $\Fq{q^4}$-points and the $\Fq{q^2}$-point lie on a $\Fq{q^2}$-line $L$ (and the other two
 $\Fq{q^4}$-points and the second $\Fq{q^2}$-point lie on $\Frob L$), or,
 \item $\Delta_{l,2}^3$ consists of $\lambda$-tuples such that the four $\Fq{q^4}$-points and the two $\Fq{q^2}$-points are intersection points of
 four conjugate $\Fq{q^4}$-lines.
\end{itemize}
The sets $\Delta_{l,2}^2$ and $\Delta_{l,2}^3$ are illustrated in Figure \ref{421Bfig} below. The cardinality
of $\Delta_{l,2}^1$ is easily computed to be $(q^2+q+1)^2(q^4-q^2)(q^2-q)$.
To get the cardinality of $\Delta_{l,2}^2$, we first choose a $\Fq{q^2}$-line $L$
in $q^4-q$ ways and then a $\Fq{q^4}$-point $P_1$ on $L$ in $q^4-q^2$ ways. This determines
all the four $\Fq{q^4}$ points since they must be $P_2=FP_1$, $P_3=F^2P_1$ and $P_4=F^3P_1$.
We must now decide if $P_5$ should lie on $L$ or $\Frob L$. We then choose a $\Fq{q^2}$-point
on the chosen line. The lines $L$ and $\Frob L$ both contain $q^2+1$ points defined over 
$\Fq{q^2}$ of which
precisely one is defined over $\Fq{q}$ (namely the point $L \cap FL$). Hence, there are $q^2$
choices for $P_5$. It now only remains to choose $P_7$ in one of $q^2+q+1$ ways. We thus have
\begin{equation*}
 |\Delta_{l,2}^2| = 2(q^4-q)(q^4-q^2)q^2(q^2+q+1).
\end{equation*}
It remains to compute $|\Delta_{l,2}^3|$.
We first choose a $\Fq{q^2}$-point $P_5$ not defined over $\Fq{q}$ in one of $q^4-q$ ways. There
are $q^4-q^2$ lines $L$ strictly defined over $\Fq{q^4}$ through $P_5$ and we choose one.
We thus get four $\Fq{q^4}$-lines which intersect in the two $\Fq{q^2}$-points $P_5$ and
$P_6$ as well as in four $\Fq{q^4}$-points. We choose one of these to become $P_1$ and
the labels of the other three points are then given. However, we could as well have chosen
the line $\Frob^2L$ and ended up with the same four $\Fq{q^4}$-points. We therefore must divide
by $2$. Finally, we choose any of the $q^2+q+1$ $\Fq{q}$-points to become $P_7$. We thus
have
\begin{equation*}
|\Delta_{l,2}^3| = 2(q^4-q)(q^4-q^2)(q^2+q+1).
\end{equation*}


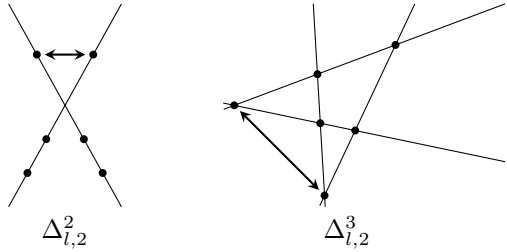
\begin{figure}[!htbp]
\begin{center}
\begin{tikzpicture}[scale=0.3,trans/.style={thick,<->,shorten >=3pt,shorten <=3pt,>=stealth}]
 \def\ptsize{5pt}
 \def\tptsize{10pt}
 \def\const{0.5}
 \draw (0,0) coordinate (b_1)--(5,9) coordinate (b_2);
 \draw (0,9) coordinate (c_1)--(5,0) coordinate (c_2);
 \coordinate (h0B) at (0,6.75);
 \coordinate (g0B) at (5,6.75);
 \coordinate (p_0B) at (intersection of b_1--b_2 and h0B--g0B);
 \coordinate (p_1B) at (intersection of c_1--c_2 and h0B--g0B);
 \fill[black] (p_0B) circle (\ptsize);
 \fill[black] (p_1B) circle (\ptsize);
 \coordinate (h1B) at (0,3);
 \coordinate (g1B) at (5,3);
 \coordinate (p_2B) at (intersection of b_1--b_2 and h1B--g1B);
 \coordinate (p_3B) at (intersection of c_1--c_2 and h1B--g1B);
 \fill[black] (p_2B) circle (\ptsize);
 \fill[black] (p_3B) circle (\ptsize);
 \coordinate (h2B) at (0,1.5);
 \coordinate (g2B) at (5,1.5);
 \coordinate (p_4B) at (intersection of b_1--b_2 and h2B--g2B);
 \coordinate (p_5B) at (intersection of c_1--c_2 and h2B--g2B);
 \fill[black] (p_4B) circle (\ptsize);
 \fill[black] (p_5B) circle (\ptsize);
 \draw[trans] (p_0B)--(p_1B);
 \coordinate (Lab_B) at (2.5,0);
 \node [below] at (Lab_B) {$\Delta_{l,2}^2$};
 
 \coordinate (p21) at (10,4.5);
 \coordinate (p22) at (14,0.5);
 \fill[black] (p21) circle (\ptsize);
 \fill[black] (p22) circle (\ptsize);
 \coordinate (h1) at (13.5,9);
 \coordinate (h2) at (22,2);
 \coordinate (h3) at (22,9);
 \coordinate (h4) at (18,9);
 \draw [shorten <=-0.15cm, shorten >=-0cm] (p22)--(h1);
 \draw [shorten <=-0.15cm, shorten >=-0cm] (p21)--(h2);
 \draw [shorten <=-0.15cm, shorten >=-0cm] (p21)--(h3);
 \draw [shorten <=-0.15cm, shorten >=-0cm] (p22)--(h4);
 \coordinate (p41) at (intersection of p22--h1 and p21--h3);
 \fill[black] (p41) circle (\ptsize);
 \coordinate (p42) at (intersection of p22--h4 and p21--h3);
 \fill[black] (p42) circle (\ptsize);
 \coordinate (p43) at (intersection of p22--h4 and p21--h2);
 \fill[black] (p43) circle (\ptsize);
 \coordinate (p44) at (intersection of p22--h1 and p21--h2);
 \fill[black] (p44) circle (\ptsize);
 \draw[trans] (p21)--(p22);
 \coordinate (Lab_C) at (15,0);
 \node [below] at (Lab_C) {$\Delta_{l,2}^3$};
\end{tikzpicture}
\end{center}
\caption{Illustration of elements of the sets $\Delta_{l,2}^2$ and $\Delta_{l,2}^3$.}
\label{421Bfig}
\end{figure}

The set $\Delta_{l,3}$ splits into two disjoint
subsets
\begin{equation*}
 \Delta_{l,3} = \Delta_{l,3}^1 \cup \Delta_{l,3}^2,
\end{equation*}
where
\begin{itemize}
 \item $\Delta_{l,3}^1$ consists of $\lambda$-tuples such that the four $\Fq{q^4}$-points and the $\Fq{q}$-point lie on a $\Fq{q}$-line, or,
 \item $\Delta_{l,3}^2$ consists of $\lambda$-tuples such that there are two conjugate $\Fq{q^2}$-lines intersecting in the $\Fq{q}$-point, each
 $\Fq{q^2}$-line containing two of the $\Fq{q^4}$-points.
\end{itemize}
To compute $|\Delta_{l,3}^1|$ we first choose a
$\Fq{q}$-line $L$, then a conjugate quadruple and a $\Fq{q}$-point on $L$ and
finally a conjugate pair of $\Fq{q^2}$-points anywhere. Hence
\begin{equation*}
|\Delta_{l,3}^1| = (q^2+q+1)(q^4-q^2)(q+1)(q^4-q).
\end{equation*}
To compute $|\Delta_{l,3}^2|$ we first choose a $\Fq{q^2}$-line $L$ not defined over
$\Fq{q}$, then a $\Fq{q^4}$-point $P_4$ not defined over $\Fq{q^2}$ on $L$ and finally
a pair of conjugate $\Fq{q^2}$-points anywhere. We thus have
\begin{equation*}
 |\Delta_{l,3}^2| = (q^4-q)^2(q^4-q^2).
\end{equation*}

We now consider the intersection $\Delta_{l,2} \cap \Delta_{l,3}$.
The decompositions above yield a decomposition
\begin{equation*}
 \Delta_{l,2} \cap \Delta_{l,3} = \bigcup_{i,j} \Delta_{l,2}^i \cap \Delta_{l,3}^j.
\end{equation*}
The intersection $\Delta_{l,2}^1 \cap \Delta_{l,3}^1$ consists of configurations where all seven points
 lie on a $\Fq{q}$-line. There are
 \begin{equation*}
  (q^2+q+1)(q^4-q^2)(q^2-q)(q+1)
 \end{equation*}
 such $\lambda$-tuples. Both the intersections $\Delta_{l,2}^1 \cap \Delta_{l,3}^2$  and  $\Delta_{l,2}^2 \cap \Delta_{l,3}^1$ are empty.
 
 To compute the cardinality of $\Delta_{l,2}^2 \cap \Delta_{l,3}^2$ we first choose a $\Fq{q^2}$-line $L$
 in $q^4-q$ ways and then a strict $\Fq{q^4}$-point $P_1$ on $L$ in $q^4-q^2$ ways. 
 We must now decide if $P_5$ should lie on $L$ or $\Frob L$. We then choose a strict $\Fq{q^2}$-point
 on the chosen line in $q^2$ ways. We are now sure to have a tuple in $\Delta_{l,2}^2$. To make sure that
 the tuple also lies in $\Delta_{l,3}^2$
 we have no choice but to put $P_7$ at the intersection of $L$ and $\Frob L$. There
 are thus
 \begin{equation*}
  2(q^4-q)(q^4-q^2)q^2
 \end{equation*}
 elements in the intersection $\Delta_{l,2}^2 \cap \Delta_{l,3}^2$.
 
 The intersection $\Delta_{l,2}^3 \cap \Delta_{l,3}^1$ is empty so there is only one intersection
 remaining.
 As explained in the computation of $|\Delta_{l,2}^3|$, there
 are $2(q^4-q)(q^4-q^2)$ ways to obtain four strict $\Fq{q^4}$-points and two strict $\Fq{q^2}$-points which are 
 the intersection points of four conjugate $\Fq{q^4}$-lines. We now note 
 that there are precisely six lines through pairs of points among the four
 $\Fq{q^4}$-points. Four of these lines are of course the four $\Fq{q^4}$-lines. The
 remaining two lines are defined over $\Fq{q^2}$ and therefore intersect
 in a $\Fq{q}$-point. To obtain a tuple in $\Delta_{l,3}^2$ we have no choice but to
 choose $P_7$ as this intersection point. Therefore, there are
 \begin{equation*}
  2(q^4-q)(q^4-q^2),
 \end{equation*}
 elements in the intersection $\Delta_{l,2}^3 \cap \Delta_{l,3}^2$.

The remaining computations are rather straightforward. One obtains
\begin{equation*}
 |\Delta_l| = q^{13}+5q^{12}-4q^{10}-5q^9-3q^8+2q^7+q^6+3q^5+q^4-q^3.
\end{equation*}

We now turn to $\Delta_c$. We have that 
 if six points of a conjugate $\lambda$-tuple lie on a smooth
 conic $C$, then the four $\Fq{q^4}$-points and the two $\Fq{q^2}$-points lie 
 on $C$ and $C$ is defined over $\Fq{q}$.
Thus, the computation of $|\Delta_c|$ consists of choosing
a smooth conic $C$ defined over $\Fq{q}$, choose a $\Fq{q^4}$-point
and a $\Fq{q^2}$-point on $C$ and finally choose a $\Fq{q}$-point anywhere.
We thus have
\begin{equation*}
 |\Delta_c| = (q^5-q^2)(q^4-q^2)(q^2-q)(q^2+q+1).
\end{equation*}

It remains to investigate the intersection $\Delta_l \cap \Delta_c$.
Since a smooth conic cannot contain three collinear points, we only have nonempty
intersection between $\Delta$ and the sets $\delta_{l,3}^2$ and
$\Delta_{l,6}$.

To compute $|\Delta_{l,3}^2 \cap \Delta_c|$ we first choose a smooth conic $C$, then a conjugate
quadruple on $C$ and finally a pair of conjugate $\Fq{q^2}$-points on
$C$. The $\Fq{q}$-point is then uniquely defined as the intersection
point of the two $\Fq{q^2}$-lines through pairs of the four $\Fq{q^4}$-points.
We thus have
\begin{equation*}
 |\Delta_{l,3}^2 \cap \Delta_c| = (q^5-q^2)(q^4-q^2)(q^2-q).
\end{equation*}
The same construction as above works for the intersection 
$\Delta_{l,6} \cap \Delta_c$ if we remember that we now have
some choice for the $\Fq{q}$-point since it can lie anywhere on the
line through the two $\Fq{q^2}$-points. We thus have
\begin{equation*}
 |\Delta_{l,6} \cap \Delta_c| = (q^5-q^2)(q^4-q^2)(q^2-q)(q+1).
\end{equation*}

The only thing that remains to compute is the cardinality of the
triple intersection $\Delta_{l,3}^2 \cap \Delta_{l,6} \cap \Delta_c$.
We first assume that the $\Fq{q}$-point is on the outside of
a smooth conic $C$ containing the other six points. 
We first choose $C$ in $q^5-q^2$ ways. There are
$\frac{1}{2}(q+1)q$ ways to choose two $\Fq{q}$-points $P$ and $Q$ on $C$.
Intersecting the tangents $T_PC$ and $T_QC$ gives a $\Fq{q}$-point
$P_7$ which will clearly lie on the outside of $C$.
Hence, there are precisely $\frac{1}{2}(q+1)q$ ways to choose a 
$\Fq{q}$-point on the outside of $C$. 

We now want to choose a $\Fq{q}$-line through $P_7$
intersecting $C$ in two $\Fq{q^2}$-points. There are $q+1$ $\Fq{q}$-lines through
$P_7$ of which two are tangents to $C$. These tangent lines
contain a $\Fq{q}$-point of $C$ each so there are $q-1$ remaining $\Fq{q}$-points
on $C$. Picking such a point gives a line through this point, $P_7$
and one further point on $C$. We thus see that exactly $\frac{1}{2}(q-1)$ of
the $\Fq{q}$-lines through $P_7$ intersect $C$ in two $\Fq{q}$-points. Hence, there
are precisely
\begin{equation*}
 q+1-2-\frac{1}{2}(q-1)=\frac{1}{2}(q-1)
\end{equation*}
$\Fq{q}$-lines through $P_7$ which intersect $C$ in two $\Fq{q^2}$-points. These points
are clearly conjugate under $\Frob$. We label one of them as $P_5$.

We shall now choose a conjugate pair of $\Fq{q^2}$-lines through $P_7$ intersecting $C$ in four $\Fq{q^4}$-points.
There are $q^2-q$ conjugate pairs of $\Fq{q^2}$-lines through $P_7$. No $\Fq{q^2}$-line
through $P_7$ is tangent to $C$ so each $\Fq{q^2}$-line through $P_7$ will intersect
$C$ in two points. The conic $C$ contains $q^2-q$ points which are defined over $\Fq{q^2}$
but not $\Fq{q}$. Picking such a point gives a line through this point and $P_7$ as
well as one further $\Fq{q^2}$-point not defined over $\Fq{q}$. Thus, there are $\frac{1}{2}(q^2-q)$ lines
obtained in this way. Typically, such a line will be defined over $\Fq{q^2}$ but not $\Fq{q}$.
We saw above that the number of such lines which are defined over $\Fq{q}$ is precisely
$\frac{1}{2}(q-1)$. Thus, there are precisely
\begin{equation*}
 \frac{1}{2}(q^2-q)-\frac{1}{2}(q-1) = \frac{1}{2}(q^2-2q+1)
\end{equation*}
$\Fq{q^2}$-lines, not defined over $\Fq{q}$, which intersect $C$ in two $\Fq{q^4}$-points.
Thus, the remaining
\begin{equation}
\label{421outeq}
 q^2-q-\frac{1}{2}(q^2-2q+1)=\frac{1}{2}(q^2-1)
\end{equation}
$\Fq{q^2}$-lines must intersect $C$ in two $\Fq{q^4}$-points. Picking such a line
and labelling one of the points $P_1$ gives a configuration belonging to
$\Delta_{l,3}^2 \cap \Delta_{l,6} \cap \Delta_c$ and we thus see that
there are
\begin{equation*}
 \frac{1}{2}(q^5-q^2)q(q+1)(q-1)(q^2-1)
\end{equation*}
such configurations with $P_7$ on the outside of $C$.

We now assume that $P_7$ is on the inside of $C$. We first choose $C$ in $q^5-q^2$ ways.
Since the number of $\Fq{q}$-points is $q^2+q+1$ and $q+1$ of these lie on $C$ the number
of $\Fq{q}$-points not on $C$ is precisely $q^2$. We just saw that $\frac{1}{2}(q+1)q$ of these
lie on the outside of $C$ so there must be
\begin{equation*}
 q^2-\frac{1}{2}(q+1)q =\frac{1}{2}(q^2-q)
\end{equation*}
$\Fq{q}$-points which lie on the inside of $C$.

Since $P_7$ now lies on the inside of $C$, every $\Fq{q}$-line through $P_7$
will intersect $C$ in two points. Exactly $\frac{1}{2}(q+1)$ will
intersect $C$ in two $\Fq{q}$-points so the remaining $\frac{1}{2}(q+1)$
will intersect $C$ in two conjugate $\Fq{q^2}$-points. We pick such a 
pair of points and label one of them $P_5$.

We now choose a conjugate pair of $\Fq{q^2}$-lines through $P_7$ intersecting $C$ in 
a conjugate quadruple of $\Fq{q^4}$-points.
The number of $\Fq{q^2}$-lines, not defined over $\Fq{q}$, through $P_7$ is $q^2-q$. Two of these are tangent to
$C$ so, using ideas similar to those above, we see that
\begin{equation*}
 \frac{1}{2}(q^2-q-2) -\frac{1}{2}(q+1) +2= \frac{1}{2}(q^2-2q+1), 
\end{equation*}
of these lines intersect $C$ in points defined over $\Fq{q^2}$. Hence,
the remaining
\begin{equation}
\label{421ineq}
 q^2-q-\frac{1}{2}(q^2-2q+1) = \frac{1}{2}(q^2-1)
\end{equation}
lines intersect $C$ in two $\Fq{q^4}$-points which are not defined over $\Fq{q^2}$.
If we pick one of these points to become $P_1$ we end up with a configuration
in $\Delta_{l,3}^2 \cap \Delta_{l,6} \cap \Delta_c$. We thus have
\begin{equation*}
 \frac{1}{2}(q^5-q^2)(q^2-q)(q+1)(q^2-1)
\end{equation*}
such configurations with $P_7$ on the inside of $C$. One may note
that the expression above actually is the same as the expression
when $P_7$ was on the outside, but this will not always be the case.

We thus have
\begin{equation*}
 |\Delta_l \cap \Delta_c| = q^{12}+q^{11}-4q^{10}-2q^9+3q^8+4q^7-4q^5+q^3,
\end{equation*}
and finally also
\begin{equation*}
 \left|\left( \Pts \right)^{\Frob \sigma} \right| = q^6-q^5-2q^4+q^3-2q^2+3.
\end{equation*}

 \subsection{The case \texorpdfstring{$\lambda=[1^3,4]$}{[1,1,1,4]}}
Throughout this section, $\lambda$ shall mean the partition
$[1^3,4]$. We shall assume that the Frobenius automorphism permutes the points $P_1, P_2,P_3, P_4$ according
to $(1234)$ and fixes the $\Fq{q}$-points $P_5$, $P_6$ and $P_7$.
Let $U$ be the open subset of $(\Pn{2})^7$ consisting of septuples such that
the last three points of the tuple are not collinear. 
In other words, we choose any conjugate quadruple and three
$\Fq{q}$-points which do not lie on a line. 

We can decompose $\Delta_l$ into a disjoint as
\begin{equation*}
 \Delta_l = \Delta_{l,1} \cup \Delta_{l,2},
\end{equation*}
where $\Delta_{l,1}$ consists of septuples such that all four $\Fq{q^4}$-points lie on a $\Fq{q}$-line and
$\Delta_{l,2}$ consists of septuples such that the $\Fq{q^2}$-line through $P_1$ and $P_3$ intersects the $\Fq{q^2}$-line through $P_2$ and $P_4$
  in $P_5$, $P_6$ or $P_7$.

To compute $|\Delta_{l,1}|$
 we simply choose a $\Fq{q}$-line $L$, a conjugate quadruple
on $L$ and finally place the three $\Fq{q}$-points in such a way that they
do not lie on a line. This can be done in
\begin{equation*}
 (q^2+q+1)^2(q^4-q^2)(q^2+q)q^2
\end{equation*}
ways. To compute $|\Delta_{l,2}|$  we first choose 
$P_5$, $P_6$ or $P_7$ then a $\Fq{q^2}$-line $L$
not defined over $\Fq{q}$ through this point. Finally, we choose a
$\Fq{q^4}$-point $P_1$ which is not defined over $\Fq{q^2}$ on $L$ and
make sure that the final two $\Fq{q}$-points are not collinear with the
first. This can be done in
\begin{equation*}
 3(q^2+q+1)(q^2-q)(q^4-q^2)(q^2+q)q^2
\end{equation*}
ways. This gives that
\begin{equation*}
 |\Delta_{l}|=4q^{12}+6q^{11}+q^{10}-4q^9-5q^8-q^7-q^5.
\end{equation*}

We now turn to investigate $\Delta_{c}$. It is not hard to see that
 if six of the points lie on a smooth conic $C$, then the four $\Fq{q^4}$-points 
 must lie on that conic and $C$ must be defined over $\Fq{q}$.
We thus choose a smooth conic $C$ over $\Fq{q}$ and a conjugate
quadruple on $C$. Then we choose $P_5$, $P_6$ or $P_7$ to possibly
not lie on $C$ and place the other two on $C$. Finally, we place
the final $\Fq{q}$-point anywhere except on the line through the
other two $\Fq{q}$-points. This gives the number
\begin{equation*}
 3(q^5-q^2)(q^4-q^2)(q+1)q \cdot q^2.
\end{equation*}
However, we have now counted the configurations where all seven
points lie on a smooth conic three times. There are
\begin{equation*}
 (q^5-q^2)(q^4-q^2)(q+1)q(q-1)
\end{equation*}
such configurations and it thus follows that
\begin{equation*}
 |\Delta_{c}| = 3q^{13}+q^{12}-3q^{11}-2q^{10}-q^9+q^8-q^7+2q^5.
\end{equation*}

We now turn to the intersection $\Delta_{l}\cap \Delta_{c}=\Delta_{l,1} \cap \Delta_{c}$.
We begin by choosing one of the
$\Fq{q}$-points $P_5$, $P_6$ and $P_7$ to not lie on the conic $C$. We call the chosen
point $P$ and the remaining two points $P_i$ and $P_j$ where $i<j$.
We now have three disjoint possibilities:
\begin{itemize}
 \item[(\emph{i})]   the point $P$ may lie on the outside of $C$ with one of the tangents through $P$ also passing through $P_i$,
 \item[(\emph{ii})]  the point $P$ may lie on the outside of $C$ with none of the tangents through $P$ passing through $P_i$,
 \item[(\emph{iii})] the point $P$ may lie on the inside of $C$.
\end{itemize}
We consider the three cases (\emph{i})-(\emph{iii}) separately.

\vspace{7pt}
\noindent \textbf{(\emph{i}).}
We begin by choosing a smooth conic $C$ in $q^5-q^2$ ways and a $\Fq{q}$-point $P$ 
on the outside of $C$ in $\frac{1}{2}(q+1)q$ ways. By Equation \ref{421outeq},
there are now $q^2-1$ ways to choose $P_1$.
Since we require $P_i$ to lie on a tangent to $C$ which passes through $P$, we
only have $2$ choices for $P_i$. Finally, we may choose $P_j$ as any of the $q$ remaining
points on $C$. We thus have
\begin{equation*}
 3(q^5-q^2)(q+1)q(q^2-1)q
\end{equation*}
possibilities in this case.

\vspace{7pt}
\noindent \textbf{(\emph{ii}).}
Again, we begin by choosing a smooth conic $C$ in $q^5-q^2$ ways and a $\Fq{q}$-point $P$ 
on the outside of $C$ in $\frac{1}{2}(q+1)q$ ways and choose $P_1$ in one of $q^2-1$ ways.
The point $P_i$ should now not lie on a tangent to $C$ which passes through $P$ so we
have $q-1$ choices. Since the line between $P$ and $P_i$ is not a tangent to $C$, there
is one further intersection point between this line and $C$. We must choose $P_j$ away
from this point and $P_i$ and thus have $q-1$ possible choices. Hence, we have
\begin{equation*}
 \frac{3}{2}(q^5-q^2)(q+1)q(q^2-1)(q-1)^2
\end{equation*}
possibilities in this case.

\vspace{7pt}
\noindent \textbf{(\emph{iii}).}
We start by choosing a smooth conic $C$ in $q^5-q^2$ ways and then a point $P$ on the inside of $C$ in $\frac{1}{2}(q^2-q)$ ways.
We continue by using Equation \ref{421ineq} to see that we can choose $P_1$ in $q^2-1$ ways.
We now choose $P_i$ as any of the $\Fq{q}$-points on $C$ and thus have $q+1$ choices.
Finally, we may choose $P_j$ as any $\Fq{q}$-point on $C$, except in the intersection of $C$ with the line through $P_i$ and $P$.
Hence, we have $q-1$ choices.
We thus have
\begin{equation*}
 \frac{3}{2} (q^5-q^2)(q^2-q)(q^2-1)(q+1)(q-1)
\end{equation*}
possibilities in this case.

We now conclude that
\begin{equation*}
 |\Delta_l \cap \Delta_c| = 3q^{11}-3q^9-3q^5+3q^3,
\end{equation*}
and, finally,
\begin{equation*}
 \left|\left( \Pts \right)^{\Frob \sigma} \right| = q^6-q^5-2q^4+q^3-2q^2+3.
\end{equation*}

\subsection{The case \texorpdfstring{$\lambda=[1,3^2]$}{[1,3,3]}}

Throughout this section, $\lambda$ shall mean the partition
$[1^1,3^2]$. We shall use
the notation $P_1$, $P_2$, $P_3$ for the first conjugate
triple and $Q_1$, $Q_2$, $Q_3$ for the second. The $\Fq{q}$-point
will be denoted by $R$. We take $U=\left( \Pn{2} \right)^7$.

We decompose $\Delta_l$ as
\begin{equation*}
 \Delta_l = \bigcup_{i=1}^5 \Delta_{l,i},
\end{equation*}
where
 \begin{itemize}
  \item $\Delta_{l,1}$ consists of $\lambda$-tuples such that the points $P_1$, $P_2$ and $P_3$ lie on a $\Fq{q}$-line,
  \item $\Delta_{l,2}$ consists of $\lambda$-tuples such that the points $P_1$, $P_2$ and $P_3$ are the intersection points of a conjugate triple of $\Fq{q^3}$-lines with each of the lines
  containing one of the points $Q_1$, $Q_2$ and $Q_3$,
  \item $\Delta_{l,3}$ consists of $\lambda$-tuples such that the points $Q_1$, $Q_2$ and $Q_3$ are the intersection points of a conjugate triple of $\Fq{q^3}$-lines with each of the lines
  containing one of the points $P_1$, $P_2$ and $P_3$,
  \item $\Delta_{l,4}$ consists of $\lambda$-tuples such that the points $Q_1$, $Q_2$ and $Q_3$ lie on a $\Fq{q}$-line, and
  \item $\Delta_{l,5}$ consists of $\lambda$-tuples such that the point $R$ is the intersection of three conjugate $\Fq{q^3}$-lines with each of the $\Fq{q^3}$-lines containing one of the
  points $P_1$, $P_2$ and $P_3$ and one of the points $Q_1$, $Q_2$ and $Q_3$.
 \end{itemize}

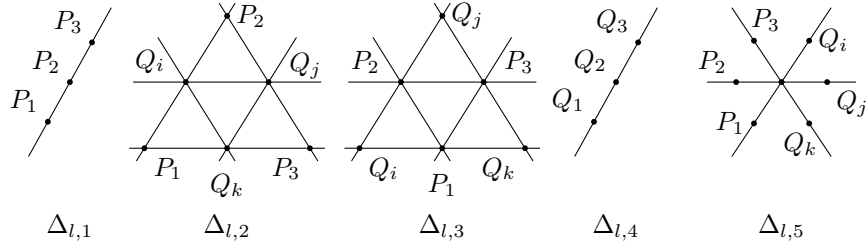
\begin{figure}[!htbp]
\begin{center}
\begin{tikzpicture}[scale=0.22]
 \def\ptsize{5pt}
 \def\tptsize{10pt}
 \def\labht{-3}
 \draw (0,0) coordinate (a_1)--(5,9) coordinate (a_2);
 \coordinate (h0A) at (0,2.1);
 \coordinate (g0A) at (5,2.1);
 \coordinate (p_0A) at (intersection of a_1--a_2 and h0A--g0A);
 \fill[black] (p_0A) circle (\ptsize);
 \node [above left] at (p_0A) {$P_1$};
 \coordinate (h2A) at (0,4.5);
 \coordinate (g2A) at (5,4.5);
 \coordinate (p_2A) at (intersection of a_1--a_2 and h2A--g2A);
 \fill[black] (p_2A) circle (\ptsize);
 \node [above left] at (p_2A) {$P_2$};
 \coordinate (h5A) at (0,6.9);
 \coordinate (g5A) at (5,6.9);
 \coordinate (p_5A) at (intersection of a_1--a_2 and h5A--g5A);
 \fill[black] (p_5A) circle (\ptsize);
 \node [above left] at (p_5A) {$P_3$};
 \coordinate (Lab_A) at (2.5,\labht);
 \node [below] at (Lab_A) {$\Delta_{l,1}$};
 
 \coordinate (b1) at (7,0.5);
 \coordinate (b2) at (12,8.5);
 \coordinate (b3) at (17,0.5);
 \fill[black] (b1) circle (\ptsize);
 \node [below right] at (b1) {$P_1$};
 \fill[black] (b2) circle (\ptsize);
 \node [right] at (b2) {$P_2$};
 \fill[black] (b3) circle (\ptsize);
 \node [below left] at (b3) {$P_3$};
 \draw [shorten <=-0.2cm, shorten >=-0.2cm] (b1)--(b2);
 \draw [shorten <=-0.2cm, shorten >=-0.2cm] (b2)--(b3);
 \draw [shorten <=-0.2cm, shorten >=-0.2cm] (b1)--(b3);
 \coordinate (c1) at (12,0.5);
 \fill[black] (c1) circle (\ptsize);
 \coordinate (c1h) at (12,-0.5);
 \node [below] at (c1h) {$Q_k$};
 \coordinate (h1) at (7,4.5);
 \coordinate (h2) at (17,4.5);
 \coordinate (c2) at (intersection of b1--b2 and h1--h2);
 \path let \p2 = (c2) in node  at (\x2-20,\y2) [above left]{$Q_i$};
 \fill[black] (c2) circle (\ptsize);
 \coordinate (c3) at (intersection of b2--b3 and h1--h2);
 \path let \p3 = (c3) in node  at (\x3+20,\y3-10) [above right]{$Q_j$};
 \fill[black] (c3) circle (\ptsize);
 \coordinate (c3h) at (17.2,4.5);
 \draw [shorten <=-0.2cm, shorten >=-0.7cm] (c1)--(c2);
 \draw [shorten <=-0.7cm, shorten >=-0.7cm] (c2)--(c3);
 \draw [shorten <=-0.7cm, shorten >=-0.2cm] (c3)--(c1);
 \coordinate (Lab_B) at (12,\labht);
 \node [below] at (Lab_B) {$\Delta_{l,2}$};
 
 \coordinate (nb1) at (20,0.5);
 \coordinate (nb2) at (25,8.5);
 \coordinate (nb3) at (30,0.5);
 \fill[black] (nb1) circle (\ptsize);
 \node [below right] at (nb1) {$Q_i$};
 \fill[black] (nb2) circle (\ptsize);
 \node [right] at (nb2) {$Q_j$};
 \fill[black] (nb3) circle (\ptsize);
 \node [below left] at (nb3) {$Q_k$};
 \draw [shorten <=-0.2cm, shorten >=-0.2cm] (nb1)--(nb2);
 \draw [shorten <=-0.2cm, shorten >=-0.2cm] (nb2)--(nb3);
 \draw [shorten <=-0.2cm, shorten >=-0.2cm] (nb1)--(nb3);
 \coordinate (nc1) at (25,0.5);
 \fill[black] (nc1) circle (\ptsize);
 \coordinate (nc1h) at (25,-0.5);
 \node [below] at (nc1h) {$P_1$};
 \coordinate (nh1) at (20,4.5);
 \coordinate (nh2) at (30,4.5);
 \coordinate (nc2) at (intersection of nb1--nb2 and nh1--nh2);
 \path let \p4 = (nc2) in node  at (\x4-20,\y4) [above left]{$P_2$};
 \fill[black] (nc2) circle (\ptsize);
 \coordinate (nc3) at (intersection of nb2--nb3 and nh1--nh2);
 \path let \p5 = (nc3) in node  at (\x5+20,\y5) [above right]{$P_3$};
 \fill[black] (nc3) circle (\ptsize);
 \draw [shorten <=-0.2cm, shorten >=-0.7cm] (nc1)--(nc2);
 \draw [shorten <=-0.7cm, shorten >=-0.7cm] (nc2)--(nc3);
 \draw [shorten <=-0.7cm, shorten >=-0.2cm] (nc3)--(nc1);
 \coordinate (Lab_C) at (25,\labht);
 \node [below] at (Lab_C) {$\Delta_{l,3}$};
 
 \draw (33,0) coordinate (a_11)--(38,9) coordinate (a_21);
 \coordinate (h0D) at (33,2.1);
 \coordinate (g0D) at (38,2.1);
 \coordinate (p_1D) at (intersection of a_11--a_21 and h0D--g0D);
 \fill[black] (p_1D) circle (\ptsize);
 \node [above left] at (p_1D) {$Q_1$};
 \coordinate (h2D) at (33,4.5);
 \coordinate (g2D) at (38,4.5);
 \coordinate (p_2D) at (intersection of a_11--a_21 and h2D--g2D);
 \fill[black] (p_2D) circle (\ptsize);
 \node [above left] at (p_2D) {$Q_2$};
 \coordinate (h5D) at (33,6.9);
 \coordinate (g5D) at (38,6.9);
 \coordinate (p_5D) at (intersection of a_11--a_21 and h5D--g5D);
 \fill[black] (p_5D) circle (\ptsize);
 \node [above left] at (p_5D) {$Q_3$};
 \coordinate (Lab_D) at (35.5,\labht);
 \node [below] at (Lab_D) {$\Delta_{l,4}$};
 
 \draw (41,4.5) coordinate (d_1)--(50,4.5) coordinate (d_2);
 \draw (42.5,0) coordinate (e_1)--(48.5,9) coordinate (e_2);
 \draw (42.5,9) coordinate (f_1)--(48.5,0) coordinate (f_2);
 \coordinate (p_kC) at (intersection of d_1--d_2 and e_1--e_2);
 \fill[black] (p_kC) circle (\ptsize);
 \coordinate (h0C) at (41,7);
 \coordinate (g0C) at (50,7);
 \coordinate (p_0C) at (intersection of e_1--e_2 and h0C--g0C);
 \coordinate (p_1C) at (intersection of f_1--f_2 and h0C--g0C);
 \fill[black] (p_0C) circle (\ptsize);
 \node [right] at (p_0C) {$Q_i$};
 \fill[black] (p_1C) circle (\ptsize);
 \path let \p6 = (p_1C) in node  at (\x6+20,\y6) [above]{$P_3$};
 \coordinate (h1C) at (41,2);
 \coordinate (g1C) at (50,2);
 \coordinate (p_2C) at (intersection of e_1--e_2 and h1C--g1C);
 \coordinate (p_3C) at (intersection of f_1--f_2 and h1C--g1C);
 \fill[black] (p_2C) circle (\ptsize);
 \node [left] at (p_2C) {$P_1$};
 \fill[black] (p_3C) circle (\ptsize);
 \path let \p7 = (p_3C) in node  at (\x7+30,\y7) [below left]{$Q_k$};
 \coordinate (h2C) at (42.75,0);
 \coordinate (g2C) at (42.75,9);
 \coordinate (p_4C) at (intersection of d_1--d_2 and h2C--g2C);
 \fill[black] (p_4C) circle (\ptsize);
 \node [above left] at (p_4C) {$P_2$};
 \coordinate (h3C) at (48.25,0);
 \coordinate (g3C) at (48.25,9);
 \coordinate (p_5C) at (intersection of d_1--d_2 and h3C--g3C);
 \fill[black] (p_5C) circle (\ptsize);
 \node [below right] at (p_5C) {$Q_j$};
 \coordinate (Lab_E) at (45.5,\labht);
 \node [below] at (Lab_E) {$\Delta_{l,5}$};
\end{tikzpicture}
\end{center}
\caption{Illustration of the decomposition of $\Delta_{l}$.}
\end{figure}

\subsubsection{\texorpdfstring{$\Delta_{l,1}$ and $\Delta_{l,4}$}{Case 1 and 4}}
The cardinalities of $\Delta_{l,1}$ and
$\Delta_{l,4}$ are obviously the same so we only make the computation for
$\Delta_{l,1}$. We thus choose a $\Fq{q}$-line $L$, a conjugate $\Fq{q^3}$-tuple
$P_1$, $P_2$, $P_3$ on $L$, a conjugate $\Fq{q^3}$-tuple $Q_1$, $Q_2$, $Q_3$ 
anywhere except equal to the other $\Fq{q^3}$-tuple and,
finally, a $\Fq{q}$-point anywhere. We thus get
\begin{equation*}
 |\Delta_{l,1}|=|\Delta_{l,4}|= (q^2+q+1)^2(q^3-q)(q^6+q^3-q^2-q-3).
\end{equation*}

\subsubsection{\texorpdfstring{$\Delta_{l,2}$ and $\Delta_{l,3}$}{Case 2 and 3}}
The cardinalities of $\Delta_{l,2}$ and $\Delta_{l,3}$ are
of course also the same. To compute $|\Delta_{l,2}|$ we first choose a conjugate
triple of lines, $L$, $\Frob L$, $\Frob^2 L$, which do not intersect in a point. There
are $q^6+q^3-q^2-q$ conjugate triples of lines of which $(q^2+q+1)(q^3-q)$
intersect in a point. There are thus $q^6-q^5-q^4+q^3$ possible triples.
We give the label $P_1$ to the point $L \cap FL$ which determines the labels
of the other two intersection points. We must now choose if $Q_1$ should lie on
$L$, $\Frob L$ or $\Frob^2 L$ and then place $Q_1$ on the chosen line. There are $3(q^3-1)$
ways to do this. Finally, we choose any $\Fq{q}$-point. We thus have
\begin{equation*}
 |\Delta_{l,2}|=|\Delta_{l,3}|=3(q^6-q^5-q^4+q^3)(q^3-1)(q^2+q+1).
\end{equation*}

\subsubsection{\texorpdfstring{$\Delta_{l,5}$}{Case 5}}
We first choose a $\Fq{q}$-point $R$ anywhere and
then a conjugate triple of lines, $L$, $\Frob L$, $\Frob^2L$ through $R$.
We then choose a point $P_1$ somewhere on $L$ in $q^3$ ways. We must
now decide if $Q_1$ should lie on $L$, $\Frob L$ or $\Frob^2L$ and then pick
a point on the chosen line in one of $q^3-1$ ways. We thus see that
\begin{equation*}
 |\Delta_{l,5}|=3(q^2+q+1)(q^3-q)q^3(q^3-1).
\end{equation*}

We must now compute the cardinalities of the different intersections.
Firstly, note that the intersection between $\Delta_{l,1}$ and $\Delta_{l,2}$ is empty.
Secondly, the size of the intersection of $\Delta_{l,1}$ and $\Delta_{l,3}$ is equal
to that of the intersection of $\Delta_{l,2}$ and $\Delta_{l,4}$. To obtain $|\Delta_{l,1} \cap \Delta_{l,3}|$
we first choose a conjugate triple of lines, $L$, $\Frob L$, $\Frob^2L$, 
which do not intersect in a point and label the intersection $L \cap FL$
by $Q_1$. We then choose a $\Fq{q}$-line $L'$ and thus get three
$\Fq{q^3}$-points $L' \cap L$, $L' \cap FL$ and $L' \cap F^2L$. We label
one of these by $P_1$. We may now choose any $\Fq{q}$-point to become
$R$. We thus see that
that
\begin{equation*}
 |\Delta_{l,1} \cap \Delta_{l,3}| = |\Delta_{l,2} \cap \Delta_{l,4}| = 3(q^6-q^5-q^4+q^3)(q^2+q+1)^2.
\end{equation*}

When we consider the intersection between $\Delta_{l,1}$ and $\Delta_{l,4}$ we must
distinguish between the cases where the two triples
lie on the same line and when they do not. A simple
computation then gives
\begin{equation*}
 |\Delta_{l,1} \cap \Delta_{l,4}| = (q^2+q+1)^2(q^3-q)(q^3-q-3) + (q^2+q+1)^2(q^2+q)(q^3-q)^2.
\end{equation*}

We continue by observing that $|\Delta_{l,1} \cap \Delta_{l,5}| = |\Delta_{l,4} \cap \Delta_{l,5}|$. To compute
$|\Delta_{l,1} \cap \Delta_{l,5}|$ we first choose a $\Fq{q}$-point $R$ and then a conjugate $\Fq{q^3}$-tuple
of lines $L$, $\Frob L$, $\Frob^2L$ through $R$. We continue by choosing
a $\Fq{q}$-line $L'$ not through $R$ in one of $q^2$ ways and then
label $L' \cap L$, $L' \cap FL$ or $L' \cap F^2L$ by $Q_1$. Finally,
we choose one of the remaining $q^3-1$ points of $L$ to become $P_1$.
Hence,
\begin{equation*}
 |\Delta_{l,1} \cap \Delta_{l,5}| = |\Delta_{l,4} \cap \Delta_{l,5}| = 3(q^2+q+1)(q^3-q)q^2(q^3-1).
\end{equation*}

The sets $\Delta_{l,2}$ and $\Delta_{l,3}$ do not intersect and neither do the sets $\Delta_{l,3}$ and 
$\Delta_{l,4}$. Hence, there are only two intersections left, namely
the one between $\Delta_{l,2}$ and $\Delta_{l,5}$ and the one between $\Delta_{l,3}$ and $\Delta_{l,5}$.
These have equal cardinalities. To compute $|\Delta_{l,2} \cap \Delta_{l,5}|$
we first choose a conjugate triple of lines, $L$, $\Frob L$, $\Frob^2L$, 
which do not intersect in a point and label the intersection $L \cap FL$
by $Q_1$. We then choose a $\Fq{q}$-point $R$. The lines between
$R$ and the points $Q_1$, $Q_2$ and $Q_3$ intersect the
lines $L$, $\Frob L$ and $\Frob^2L$ in three points and we label
one of them by $P_1$. We thus have
\begin{equation*}
 |\Delta_{l,2} \cap \Delta_{l,5}| = |\Delta_{l,3} \cap \Delta_{l,5}| = 3 (q^6-q^5-q^4+q^3)(q^2+q+1).
\end{equation*}

There is only one triple intersection, namely between $\Delta_{l,1}$, $\Delta_{l,4}$
and $\Delta_{l,5}$. To compute the size of this intersection we first choose
a $\Fq{q}$-point $R$ and then a conjugate triple of lines, $L$, $\Frob L$, $\Frob^2L$ 
through $R$. We then choose a $\Fq{q}$-line $L'$ not through $R$ and
label the intersection $L \cap L'$ by $P_1$. Finally, we choose
another $\Fq{q}$-line $L''$ and label one of the intersections $L'' \cap L$,
$L'' \cap FL$ and $L'' \cap F^2L$ by $Q_1$. This shows that
\begin{equation*}
 |\Delta_{l,1} \cap \Delta_{l,4} \cap \Delta_{l,5}| = 3(q^2+q+1)(q^3-q)q^2(q^2-1).
\end{equation*}

We now turn to the computation of $|\Delta_c|$. If six points of 
a $\lambda$-tuple lie on a smooth conic $C$, then
both of the conjugate triples must lie on $C$ and $C$ must be
defined over $\Fq{q}$. Hence, to obtain $|\Delta_c|$ we only
have to choose a smooth $\Fq{q}$-conic $C$, two conjugate triples
on $C$ and a $\Fq{q}$-point anywhere. We thus have that,
\begin{equation*}
 |\Delta_c| = (q^5-q^2)(q^3-q)(q^3-q-3)(q^2+q+1).
\end{equation*}

Since the sets $\Delta_{l,1}$, $\Delta_{l,2}$, $\Delta_{l,3}$ and $\Delta_{l,4}$ all require three of the $\Fq{q^3}$-points
to lie on a line, they will have empty intersection with $\Delta_c$.
This is however not true for the set $\Delta_{l,5}$. To obtain such a configuration
we first choose a smooth conic $C$ and a $\Fq{q}$-point $R$. Now
choose a $\Fq{q^3}$-point $P_1$ on $C$ which is not defined over $\Fq{q}$
in $q^3-q$ ways. Since both $C$ and $R$ are defined over $\Fq{q}$ we
know that any tangent to $C$ which passes through $R$ must
either be defined over $\Fq{q^2}$ (or $\Fq{q}$). Hence, the line through $R$
and $P_1$ will intersect $C$ in $P_1$ and one point more. We label
this point with $Q_1$, $Q_2$ or $Q_3$. We thus have
\begin{equation*}
 |\Delta_{l} \cap \Delta_c| = 3(q^5-q^2)q^2(q^3-q).
\end{equation*}
We now obtain
\begin{equation*}
 \left|\left( \Pts \right)^{\Frob \sigma} \right| = q^6-2q^5-2q^4-8q^3+16q^2+10q+21.
\end{equation*}

\subsection{The case \texorpdfstring{$\lambda=[2^2,3]$}{[2,2,3]}}
Throughout this section, $\lambda$ shall mean the partition
$[2^2,3^1]$. We shall use
the notation $P_1$, $P_2$, $P_3$ for the conjugate
triple and $Q_1$, $Q_2$ and $R_1$, $R_2$ for the two conjugate pairs.
We let $U=\left( \Pn{2} \right)^7$.  

We can decompose $\Delta_{l}$ as
\begin{equation*}
 \Delta_{l} = \Delta_{l,1} \cup \Delta_{l,2},
\end{equation*}
where $\Delta_{l,1}$ consists of septuples such that the three $\Fq{q^3}$-points
lie on a line and $\Delta_{l,2}$ consists of septuples such that the four $\Fq{q^2}$-points
lie on a line.

We have
\begin{equation*}
 |\Delta_{l,1}| = (q^2+q+1)(q^3-q)(q^4-q)(q^4-q-2),
\end{equation*}
and
\begin{equation*}
 |\Delta_{l,2}| = (q^2+q+1)(q^2-q)(q^2-q-2)(q^6+q^3-q^2-q).
\end{equation*}
The cardinality of the intersection is easily computed to
be
\begin{equation*}
 |\Delta_{l,1} \cap \Delta_{l,2}| =(q^2+q+1)^2(q^3-q)(q^2-q)(q^2-q-2).
\end{equation*}
This allows us to compute $|\Delta_{l}|$. 

We have that if six of the points of a $\lambda$-tuple lie
on a smooth conic $C$, then all seven points lie on $C$ and
$C$ is defined over $\Fq{q}$. We thus have that $\Delta_{c}$
is disjoint from $\Delta_{l}$.
We also see
that
\begin{equation*}
 |\Delta_{c}| = (q^5-q^2)(q^3-q)(q^2-q)(q^2-q-2),
\end{equation*}
so,
\begin{equation*}
 \left|\left( \Pts \right)^{\Frob \sigma} \right| = q^6-q^5-2q^4+3q^3+q^2-2q.
\end{equation*}

\subsection{The case \texorpdfstring{$\lambda=[1^2,2,3]$}{[1,1,2,3]}}
Throughout this section, $\lambda$ shall mean the partition
$[1^2,2^1,3^1]$. We shall use
the notation $P_1$, $P_2$, $P_3$ for the conjugate
triple, $Q_1$, $Q_2$ for the conjugate pair and use $R_6$ and $R_7$
to denote the two $\Fq{q}$-points.  
Let $U=\left( \Pn{2} \right)^7$. 

We decompose $\Delta_l$ as
\begin{equation*}
 \Delta_l = \Delta_{l,1} \cup \Delta_{l,2},
\end{equation*}
where $\Delta_{l,1}$ consists of $\lambda$-tuples such that
the three $\Fq{q^3}$-points lie on a $\Fq{q}$-line
and $\Delta_{l,2}$ consists of $\lambda$-tuples such that
the two $\Fq{q^2}$-points and one of the $\Fq{q}$-points lie on a $\Fq{q}$-line.

We have
\begin{equation*}
 |\Delta_{l,1}| = (q^2+q+1)^2(q^3-q)(q^4-q)(q^2+q).
\end{equation*}
We decompose $\Delta_{l,2}$ as
\begin{equation*}
 \Delta_{l,2} = \Delta_{l,2}^6 \cup \Delta_{l,2}^7,
\end{equation*}
where $\Delta_{l,2}^i$ consists of tuples such that the line through the two $\Fq{q^2}$-points passes
through $R_i$.
We have
\begin{equation*}
 |\Delta_{l,2}^6|=|\Delta_{l,2}^7|= (q^2+q+1)(q^2-q)(q+1)(q^6+q^3-q^2-q)(q^2+q).
\end{equation*}

We now turn to the double intersections.
We have
\begin{equation*}
 |\Delta_{l,1} \cap \Delta_{l,2}^6| = |\Delta_{l,1} \cap \Delta_{l,2}^7| = (q^2+q+1)^2(q^3-q)(q^2-q)(q+1)(q^2+q),
\end{equation*}
and
\begin{equation*}
 |\Delta_{l,2}^6 \cap \Delta_{l,2}^7| = (q^2+q+1)(q^2-q)(q+1)q(q^6+q^3-q^2-q).
\end{equation*}
Finally, we compute the cardinality of the intersection of all three sets
\begin{equation*}
 |\Delta_{l,1} \cap \Delta_{l,2}^6 \cap \Delta_{l,2}^7| = (q^2+q+1)^2(q^3-q)(q^2-q)(q+1)q.
\end{equation*}
This now allows us to compute $\Delta_{l}$.

If six points of a $\lambda$-tuple lie on a smooth conic $C$, then
the three $\Fq{q^3}$-points and the two $\Fq{q^2}$-points lie on
$C$ and $C$ is defined over $\Fq{q}$.
Thus, to compute $|\Delta_{c}|$ we begin by choosing a smooth
conic $C$ over $\Fq{q}$. We then choose a conjugate triple and
a conjugate pair of $\Fq{q^2}$-points on $C$. Then, we choose either
$R_6$ or $R_7$ and place the chosen point on $C$. Finally, we place
the remaining point anywhere we want. We thus obtain the number
\begin{equation*}
 2(q^5-q^2)(q^3-q)(q^2-q)(q+1)(q^2+q).
\end{equation*}
However, in the above we have counted the configurations where
all seven points lie on the conic twice. We thus have to take away
\begin{equation*}
 (q^5-q^2)(q^3-q)(q^2-q)(q+1)q,
\end{equation*}
in order to obtain $|\Delta_{c}|$.

It only remains to compute the cardinality of the intersection
$\Delta_{l} \cap \Delta_{c}$. We only have non\-empty
intersection between the set $\Delta_{c}$ and the set
$\Delta_{l,2}$. To compute the cardinality of this intersection, we only have to make sure to choose
the point $R_6$ (resp. $R_7$) on the line through the two $\Fq{q^2}$-points.
Hence, we have
\begin{equation*}
 |\Delta_{l,2}^6 \cap \Delta_{c}| = |\Delta_{l,2}^7 \cap \Delta_{c}| = (q^5-q^2)(q^3-q)(q^2-q)(q+1)^2,
\end{equation*}
and, therefore,
\begin{equation*}
 |\Delta_{l} \cap \Delta_{c}| = 2(q^5-q^2)(q^3-q)(q^2-q)(q+1)^2.
\end{equation*}
This gives us
\begin{equation*}
  \left|\left( \Pts \right)^{\Frob \sigma} \right| =  q^6-3q^5-q^4+5q^3-2q.
\end{equation*}

\subsection{The case \texorpdfstring{$\lambda=[1^4,3^1]$}{[1,1,1,1,3]}}
Throughout this section, $\lambda$ shall mean the partition
$[1^4,3^1]$. We shall denote the four $\Fq{q}$-points by $P_1$, $P_2$, $P_3$ and $P_4$
and denote the conjugate triple by $Q_1$, $Q_2$, $Q_3$.
Let $U \subset ( \Pn{2} )^7$ be the subset
consisting of septuples of points with the first four in general position.

A septuple in $\Delta_l$ will have the three $\Fq{q^3}$-points on a
$\Fq{q}$-line. Thus, to compute the size of $\Delta_{l}$, we only need to place
the four $\Fq{q}$-points in general position, choose $\Fq{q}$-line $L$ and place
the conjugate $\Fq{q^3}$-tuple on $L$. We thus have,
\begin{equation*}
 |\Delta_{l}|= (q^2+q+1)(q^2+q)q^2(q^2-2q+1)(q^2+q+1)(q^3-q).
\end{equation*}

A septuple in $\Delta_c$ will have the three $\Fq{q^3}$-points on a
smooth conic $C$ defined over $\Fq{q}$.
Thus, to compute $|\Delta_{c}|$, we first choose a smooth conic
$C$ defined over $\Fq{q}$ and then a conjugate triple on $C$. We then
choose one of the points $P_1$, $P_2$, $P_3$ and $P_4$ to possibly
not lie on $C$. Call this point $P$.
We then place the other three points on $C$. These three points
define three lines which, in total, contain $(q+1)+q+(q-1)=3q$ points. As
long as we choose $P$ away from these points, the four $\Fq{q}$-points will be
in general position. We thus obtain
\begin{equation*}
 4(q^5-q^2)(q^3-q)(q+1)q(q-1)(q^2-2q+1).
\end{equation*}
However, we have counted the septuples with all seven
points on a smooth conic four times. We thus need to take away
\begin{equation*}
 3(q^5-q^2)(q^3-q)(q+1)q(q-1)(q-2).
\end{equation*}
Since $\Delta_l$ and $\Delta_c$ are disjoint we are done and
conclude that
\begin{equation*}
 \left|\left( \Pts \right)^{\Frob \sigma} \right| = q^6-5q^5+10q^4-5q^3-11q^2+10q.
\end{equation*}

\subsection{The case \texorpdfstring{$\lambda=[1,2^3]$}{[1,2,2,2]}}
Throughout this section, $\lambda$ shall mean the partition
$[1^1,2^3]$. We shall denote the three conjugate pairs of $\Fq{q^2}$-points
by $P_1$, $P_2$, $Q_1$, $Q_2$ and $R_1$, $R_2$ and the $\Fq{q}$-point
by $O$.
Let $U \subset (\Pn{2})^7$ be the subset
consisting of septuples of points such that the first six points have no
three on a line.

We decompose $\Delta_l$ as
\begin{equation*}
 \Delta_{l} = \Delta_{l,1} \cup \Delta_{l,2},
\end{equation*}
where $\Delta_{l,1}$ consists of those septuples where
two conjugate $\Fq{q^2}$-points and the $\Fq{q}$-point lie on a $\Fq{q}$-line and
$\Delta_{l,2}$ consists of those septuples where
two conjugate $\Fq{q^2}$-lines, containing two $\Fq{q^2}$-points each, intersect in the point
defined over $\Fq{q}$. 

The set $\Delta_{l,1}$ naturally decomposes into three
equally large, but not disjoint, subsets:
\begin{itemize}
 \item the set $\Delta_{l,1}^a$ where $P_1$, $P_2$ and $O$ lie on a $\Fq{q}$-line,
 \item the set $\Delta_{l,1}^b$ where $Q_1$, $Q_2$ and $O$ lie on a $\Fq{q}$-line, and
 \item the set $\Delta_{l,1}^c$ where $R_1$, $R_2$ and $O$ lie on a $\Fq{q}$-line.
\end{itemize}
Similarly, the set $\Delta_{l,2}$ decomposes into
six disjoint and equally large subsets:
\begin{itemize}
 \item the two sets $\Delta_{l,2}^{P_1,Q_i}$ where the line through the points $P_1$ and $Q_i$ also passes through the point $O$,
 \item the two sets $\Delta_{l,2}^{P_1,R_i}$ where the line through the points $P_1$ and $R_i$ also passes through the point $O$, and
 \item the two sets $\Delta_{l,2}^{Q_1,R_i}$ where the line through the points $Q_1$ and $R_i$ also passes through the point $O$.
\end{itemize}
The cardinalities of these sets are easily computed to be
\begin{equation*}
 |\Delta_{l,1}^a|=|\Delta_{l,1}^b|=|\Delta_{l,1}^c|=(q^4-q)(q^4-q^2)(q^4-6q^2+q+8)(q+1),
\end{equation*}
and
\begin{equation*}
 |\Delta_{l,2}^{P_1,Q_i}| = |\Delta_{l,2}^{P_1,R_i}| = |\Delta_{l,2}^{Q_1,R_i}| = (q^4-q)(q^4-q^2)(q^4-6q^2+q+8).
\end{equation*}
To compute $|\Delta_{l,1}^a \cap \Delta_{l,1}^b|$ we note that if we place
the three pairs of $\Fq{q^2}$-points such that no three lie on a line, then
the line through $P_1$ and $P_2$ and the line through $Q_1$ and $Q_2$
will intersect in a $\Fq{q}$-point. By choosing this point as $O$ we obtain
an element of $\Delta_{l,1}^a \cap \Delta_{l,1}^b$. 
We now see that
\begin{equation*}
 |\Delta_{l,1}^a \cap \Delta_{l,1}^b| = |\Delta_{l,1}^a \cap \Delta_{l,1}^c| = 
 |\Delta_{l,1}^b \cap \Delta_{l,1}^c| = (q^4-q)(q^4-q^2)(q^4-6q^2+q+8).
\end{equation*}

To compute $|\Delta_{l,1}^a \cap \Delta_{l,2}^{Q_1,R_1}|$
we first choose a conjugate pair $Q_1$, $Q_2$ and then
a conjugate pair of $\Fq{q^2}$-points $R_1$, $R_2$ which do not
lie on the line through $Q_1$ and $Q_2$. We now only have one choice for
$O$. We choose a $\Fq{q}$-line $L$ through $O$. There are two possibilities:
either $L$ will pass through the intersection point $P$ of the line
through $Q_1$ and $R_2$ and the line through $Q_2$ and $R_1$ or it will not.
If $L$ passes through $P$, then we have $q^2-q$ possible choices
for $P_1$ and $P_2$ on $L$. Otherwise, we only have $q^2-q-2$ choices.
Hence
\begin{align*}
 & |\Delta_{l,1}^a \cap \Delta_{l,2}^{Q_1,R_i}| = 
 |\Delta_{l,1}^b \cap \Delta_{l,2}^{P_1,R_i}| = 
 |\Delta_{l,1}^c \cap \Delta_{l,2}^{P_1,Q_i}| =  \\
 & = (q^4-q)(q^4-q^2)(q^2-q)+ (q^4-q)(q^4-q^2)q(q^2-q-2).
\end{align*}

The only non\-empty triple intersection is $\Delta_{l,1}^a \cap \Delta_{l,1}^b \cap \Delta_{l,1}^c$.
A computation very similar to the one for $|\Delta_{l,1}^a \cap \Delta_{l,2}^{Q_1,R_1}|$ gives
\begin{align*}
 |\Delta_{l,1}^a \cap \Delta_{l,1}^b \cap \Delta_{l,1}^c| & = 2(q^4-q)(q^4-q^2)(q^2-q-2)+\\
  & + (q^4-q)(q^4-q^2)(q-3)(q^2-q-4).
\end{align*}
This finishes the investigation of $\Delta_{l}$.

We now turn to $\Delta_{c}$. If six points of a $\lambda$-tuple 
lie on a smooth conic $C$, then the six $\Fq{q^2}$-points lie on $C$ and
$C$ is defined over $\Fq{q}$.
Since no three points of a smooth conic can lie on a line, we shall
obtain an element of $\Delta_{c}$ simply by
choosing a smooth conic $C$, three conjugate pairs
on $C$ and, finally, a $\Fq{q}$-point anywhere. We thus have
\begin{equation*}
 |\Delta_{c}| = (q^5-q^2)(q^2-q)(q^2-q-2)(q^2-q-4)(q^2+q+1).
\end{equation*}

We shall now compute the cardinality of the intersection
between $\Delta_{l}$ and $\Delta_{c}$.
The intersections with the cases $\Delta_{l,1}^a$, $\Delta_{l,1}^b$, and $\Delta_{l,1}^c$ are
easily handled: we simply choose a smooth conic with three conjugate
pairs on it and then place $O$ on the line through the right conjugate
pair. We thus get
\begin{align*}
 |\Delta_{l,1}^a \cap \Delta_{c}| = |\Delta_{l,1}^b \cap \Delta_{c}| = |\Delta_{l,1}^c \cap \Delta_{c}| =  \\
 = (q^5-q^2)(q^2-q)(q^2-q-2)(q^2-q-4)(q+1).
\end{align*}
The intersections with the sets $\Delta_{l,2}^{a_1,Q_i}$, $\Delta_{l,2}^{P_1,R_i}$ and
$\Delta_{l,2}^{Q_1,R_i}$ are perhaps even simpler: once we have chosen
our conic $C$ and our conjugate pairs we have only one choice
for $O$. Hence,
\begin{align*}
 & |\Delta_{l,2}^{P_1,Q_i} \cap \Delta_{c}| = |\Delta_{l,2}^{P_1,R_i} \cap \Delta_{c}| = |\Delta_{l,2}^{Q_1,R_i} \cap \Delta_{c}| =  \\
 & = (q^5-q^2)(q^2-q)(q^2-q-2)(q^2-q-4).
\end{align*}
An analogous argument shows that
\begin{align*}
 & |\Delta_{l,1}^a \cap \Delta_{l,1}^b \cap \Delta_{c}| = |\Delta_{l,1}^a \cap \Delta_{l,1}^c \cap \Delta_{c}| 
 = |\Delta_{l,1}^b \cap \Delta_{l,1}^c \cap \Delta_{c}| = \\
 & = (q^5-q^2)(q^2-q)(q^2-q-2)(q^2-q-4).
\end{align*}

The remaining intersections are quite a bit harder than the previous
ones. We consider $\Delta_{l,1}^a \cap \Delta_{l,2}^{Q_1,R_1} \cap \Delta_{c}$,
but the other intersections of this type are completely analogous
and have the same size.

We first consider the case when $O$ is on the outside
of $C$. There are $q+1$ lines through $O$. Of these, precisely
two are tangents and $\frac{1}{2}(q-1)$ intersect $C$ in $\Fq{q}$-points.
Thus, the remaining $\frac{1}{2}(q-1)$ lines will intersect $C$
in two conjugate $\Fq{q^2}$-points. We thus pick one of these lines
and label one of the intersection points by $P_1$. 

Picking a $\Fq{q^2}$-point not defined over $\Fq{q}$ on $C$ will typically define a $\Fq{q^2}$-line
through $O$ which is not defined over $\Fq{q}$. However, some of these
choices will give $\Fq{q}$-lines and we saw above that the number of such $\Fq{q}$-lines
is $\frac{1}{2}(q-1)$. Thus, the number of $\Fq{q^2}$-lines, not defined over $\Fq{q}$, intersecting
$C$ in two $\Fq{q^2}$-points is
\begin{equation*}
 \frac{1}{2}(q^2-q)-\frac{1}{2}(q-1)=\frac{1}{2}(q^2-2q+1).
\end{equation*}
We pick one such line, label one of the intersection points $Q_1$
and the other intersection point $R_1$. This gives us a configuration
of the desired type. \linebreak Hence, the number of tuples in $\Delta_{l,1}^a \cap \Delta_{l,2}^{Q_1,R_1} \cap \Delta_{c}$
with $O$ on the outside of $C$ is
\begin{equation*}
 \frac{1}{2}(q^5-q^2)(q+1)q(q-1)(q^2-2q+1).
\end{equation*}

We now turn to the case when $O$ is on the inside of $C$.
Of the $q+1$ lines defined over $\Fq{q}$ which pass through $O$,
$\frac{1}{2}(q+1)$ will now intersect $C$ in $\Fq{q}$-points and
equally many in conjugate $\Fq{q^2}$-points. We thus pick a line that intersects
$C$ in two conjugate $\Fq{q^2}$-points and label one of them by $P_1$.

We now want to pick a $\Fq{q^2}$-line through $O$ which is not defined over $\Fq{q}$ 
and which intersects $C$ in two $\Fq{q^2}$-points that are not defined over $\Fq{q}$. 
To obtain such a line we pick a $\Fq{q^2}$-point which is not defined over $\Fq{q}$
on $C$. However, two such points define tangents to $C$ which pass
through $O$ and $\frac{1}{2}(q+1)$ of the lines obtained in this
way are actually defined over $\Fq{q}$. We thus have
\begin{equation*}
 \frac{1}{2}(q^2-q-2)-\frac{1}{2}(q+1)=\frac{1}{2}(q^2-2q-3)
\end{equation*}
choices. We pick such a line and label the intersection
points by $Q_1$ and $R_1$. Hence, the number of tuples in $\Delta_{l,1}^a \cap \Delta_{l,2}^{Q_1,R_1} \cap \Delta_{c}$
with $O$ on the inside of $C$ is
\begin{equation*}
 \frac{1}{2}(q^5-q^2)(q^2-q)(q+1)(q^2-2q-3).
\end{equation*}
This finishes the computation of $|\Delta_{l,1}^a \cap \Delta_{l,2}^{Q_1,R_i} \cap \Delta_{c}|$, 
$|\Delta_{l,1}^b \cap \Delta_{l,2}^{P_1,R_i} \cap \Delta_{c}|$
and $|\Delta_{l,1}^c \cap \Delta_{l,2}^{P_1,Q_i} \cap \Delta_{c}|$.

The only remaining intersection is $\Delta_{l,1}^a \cap \Delta_{l,1}^b \cap \Delta_{l,1}^c \cap \Delta_{c}$
which we shall handle in a way similar to that above. Fortunately,
much of the work has already been done. To start, if $O$ is on
the outside of $C$, then there are $\frac{1}{2}(q-1)$ lines though $O$
which are defined over $\Fq{q}$ and intersect $C$ in conjugate pairs of $\Fq{q^2}$-points.
Thus, there are
\begin{equation*}
 (q-1)(q-3)(q-5)
\end{equation*}
ways to pick three lines and label the intersection points
with $P_1$ and $P_2$, $Q_1$ and $Q_2$ and $R_1$ and $R_2$.
Hence, the number of $\lambda$-tuples in $\Delta_{l,1}^a \cap \Delta_{l,1}^b \cap \Delta_{l,1}^c \cap \Delta_{c}$
with $O$ on the outside of $C$ is
\begin{equation*}
 \frac{1}{2}(q^5-q^2)(q+1)q(q-1)(q-3)(q-5).
\end{equation*}
Similarly, if $O$ lies on the inside of $C$ we have seen that
there are $\frac{1}{2}(q+1)$ lines through $O$ which are defined
over $\Fq{q}$ and which intersect $C$ in a pair of conjugate $\Fq{q^2}$-points.
Thus, there are
\begin{equation*}
 (q+1)(q-1)(q-3)
\end{equation*}
ways to pick three lines and label the intersection points
with $P_1$ and $P_2$, $Q_1$ and $Q_2$ and $R_1$ and $R_2$.
Hence, the number of $\lambda$-tuples in $\Delta_{l,1}^a \cap \Delta_{l,1}^b \cap \Delta_{l,1}^c \cap \Delta_{c}$
with $O$ on the inside of $C$ is
\begin{equation*}
 \frac{1}{2}(q^5-q^2)(q+1)q(q+1)(q-1)(q-3).
\end{equation*}
We finally obtain
\begin{equation*}
 \left|\left( \Pts \right)^{\Frob \sigma} \right| = q^6-3q^5-6q^4+19q^3+6q^2-24q+7.
\end{equation*}

\subsection{The case \texorpdfstring{$\lambda=[1^3,2^2]$}{[1,1,1,2,2]}} 
Throughout this section, $\lambda$ shall mean the partition
$[1^3,2^2]$. We shall denote the $\Fq{q}$-points by $P_1$, $P_2$ and $P_3$
and the two conjugate pairs of $\Fq{q^2}$-points
by $Q_1$, $Q_2$ and $R_1$, $R_2$.
Let $U \subset (\Pn{2})^7$ be the subset
consisting of septuples of points such that the first five points lie
in general position. 

The set $\Delta_l$ can be decomposed as
\begin{equation*}
 \Delta_l = \Delta_{l,1} \cup \Delta_{l,2} \cup \Delta_{l,3},
\end{equation*}
where
 \begin{itemize}
  \item $\Delta_{l,1}$ consists of tuples such that the line through $R_1$ and $R_2$ also passes through $P_1$, $P_2$
  or $P_3$, 
  \item $\Delta_{l,2}$ consists of tuples such that the points $R_1$ and $R_2$ lie on the line through $Q_1$ and $Q_2$, and
  \item $\Delta_{l,3}$ consists of tuples such that a line through $Q_1$ and one of the points $P_1$, $P_2$ and $P_3$
  also contains $R_1$ or $R_2$.
 \end{itemize}
The set $\Delta_{l,1}$ 
decomposes as a union of the sets $\Delta_{l,1}^1$, $\Delta_{l,1}^2$ and $\Delta_{l,1}^3$ consisting
of tuples with the line through $R_1$ and $R_2$ passing through
$P_1$, $P_2$ and $P_3$, respectively. Similarly, the set $\Delta_{l,3}$ is
the union of the six sets $\Delta_{l,3}^{i,j}$, $i=1,2$, $j=1,2,3$, where
$\Delta_{l,3}^{i,j}$ contains all tuples such that $Q_1$, $R_i$ and $P_j$
lie on a line.

The cardinalities of the above sets are easily computed to be
\begin{align*}
 |\Delta_{l,1}^i| & = (q^2+q+1)(q^2+q)q^2(q^4-3q^3+3q^2-q)(q+1)(q^2-q), \\
 |\Delta_{l,2}| & = (q^2+q+1)(q^2+q)q^2(q^4-3q^3+3q^2-q)(q^2-q-2),
\end{align*}
and
\begin{equation*}
 |\Delta_{l,3}^{i,j}| = (q^2+q+1)(q^2+q)q^2(q^4-3q^3+3q^2-q)(q^2-1).
\end{equation*}
The cardinality of $\Delta_{l,1}^i \cap \Delta_{l,1}^j$, $i \neq j$, is also easily computed:
\begin{equation*}
 |\Delta_{l,1}^i \cap \Delta_{l,1}^j| = (q^2+q+1)(q^2+q)q^2(q^4-3q^3+3q^2-q)(q^2-q).
\end{equation*}
There is only non\-empty intersection between the set $\Delta_{l,1}^i$ and the
set $\Delta_{l,3}^{j,k}$ if $k \neq i$. We then place the first five points in general
position and choose a $\Fq{q}$-line through $P_k$ which does not pass through
$P_i$ in $q$ ways. This gives a tuple of the desired form. We thus see that
\begin{equation*}
 |\Delta_{l,1}^i \cap \Delta_{l,3}^{j,k}| = (q^2+q+1)(q^2+q)q^2(q^4-3q^3+3q^2-q)q.
\end{equation*}
We also have non\-empty intersection between the sets $\Delta_{l,3}^{1,i}$ and
the set $\Delta_{l,3}^{2,j}$ where $i \neq j$. Such a configuration is actually
given by specifying the first five points in general position since
we must then take $R_1$ as the intersection point of the line
between $Q_1$ and $P_i$ and the line between $Q_2$ and $P_j$
and similarly for $R_2$. Hence,
\begin{equation*}
 |\Delta_{l,3}^{1,i} \cap \Delta_{l,3}^{2,j}| = (q^2+q+1)(q^2+q)q^2(q^4-3q^3+3q^2-q).
\end{equation*}
Since the set $\Delta_{l,2}$ cannot intersect any of the other sets, because
this would require $Q_1$ and $Q_2$ to lie on a line through one of the 
$\Fq{q}$-points, it is now time to consider the triple intersections.

Since $P_1$, $P_2$ and $P_3$ do not lie on a line we have that
the intersection of $\Delta_{l,1}^1$, $\Delta_{l,1}^2$ and $\Delta_{l,1}^3$ is empty. We thus
only have two types of triple intersections, namely $\Delta_{l,1}^i \cap \Delta_{l,1}^j \cap \Delta_{l,3}^{r,s}$ 
and $\Delta_{l,1}^i \cap \Delta_{l,3}^{1,j} \cap \Delta_{l,3}^{2,k}$ where, of course, $i$, $j$ and $k$ are assumed
to be distinct.

An element of $\Delta_{l,1}^i \cap \Delta_{l,1}^j \cap \Delta_{l,3}^{r,s}$ is specified by choosing
the first five points in general position. The point $R_r$ must then
be chosen as the intersection point of the line between $P_i$ and $P_j$
and the line between $Q_1$ and $P_s$ and similarly for $\Frob R_r$. We thus
have
\begin{equation*}
 |\Delta_{l,1}^i \cap \Delta_{l,1}^j \cap \Delta_{l,3}^{r,s}|= (q^2+q+1)(q^2+q)q^2(q^4-3q^3+3q^2-q).
\end{equation*}

To compute the cardinality of the intersection $\Delta_{l,1}^i \cap \Delta_{l,3}^{1,j} \cap \Delta_{l,3}^{2,k}$
we first choose two $\Fq{q}$-points $P_j$ and $P_k$. We then
choose a conjugate pair of $\Fq{q^2}$-lines through each of these points.
The intersections of these lines give four $\Fq{q^2}$-points which
we only have one way to label with $Q_1$, $Q_2$, $R_1$ and $R_2$.
 We must now place the point $P_i$ somewhere
on the line $L$ through $R_1$ and $R_2$. The line through $P_i$ and
$P_k$ intersects $L$ in one $\Fq{q}$-point and the line through $Q_1$ and $Q_2$
intersects $L$ in another. Thus, we have $q-1$ choices for $P_i$. We thus
see that
\begin{equation*}
 |\Delta_{l,1}^i \cap \Delta_{l,3}^{1,j} \cap \Delta_{l,3}^{2,s}| = (q^2+q+1)(q^2+q)(q^2-q)^2(q-1).
\end{equation*}
This completes the investigation of $\Delta_l$.

If a smooth conic $C$ contains six of the points, then $C$ contains
 both the conjugate $\Fq{q^2}$-pairs and $C$ is defined over $\Fq{q}$.
Thus, to compute $|\Delta_c|$ we first choose a smooth conic $C$ over $\Fq{q}$
and then pick one of the points $P_1$, $P_2$ and $P_3$ to possibly lie outside $C$. 
We call the chosen point $P$. We then
place the other two points and the two $\Fq{q^2}$-pairs on $C$. Finally,
we must place $P$ somewhere to make $P_1$, $P_2$, $P_3$, $Q_1$ and
$Q_2$ lie in general position. Hence, we must choose $P$ away
from the line through the two other $\Fq{q}$-points and away from 
the line through $Q_1$ and $Q_2$. This gives us
\begin{equation*}
 (q^5-q^2)(q+1)q(q^2-q)(q^2-q-2)(q^2-q).
\end{equation*}
However, in the above we have counted the configurations where
all seven points lie on $C$ three times. We must therefore take
away
\begin{equation*}
 2 \cdot (q^5-q^2)(q+1)q(q-1)(q^2-q)(q^2-q-2)
\end{equation*}
in order to obtain $|\Delta_c|$.

The intersection $\Delta_{l,2} \cap \Delta_c$ is empty but the
intersections of $\Delta_c$ with the other sets
in the decomposition of $\Delta_l$ are not. 
To compute $|\Delta_{l,1}^i \cap \Delta_c|$ we shall first
assume that $P_i$ lies on the outside of $C$. Of the $q+1$ lines
through $P_i$ which are defined over $\Fq{q}$ we have that $2$ are tangent
to $C$ and $\frac{1}{2}(q-1)$ intersect $C$ in two $\Fq{q}$-points. Thus,
there are $\frac{1}{2}(q-1)$ lines left which must intersect $C$ in
a pair of conjugate $\Fq{q^2}$-points. We pick such a line and label the
intersection points by $R_1$ and $R_2$ in one of two ways. 
We shall now place the other two $\Fq{q}$-points on $C$. There
are $\frac{1}{2}(q+1)q$ ways to choose two $\Fq{q}$-points on $C$ of which
$\frac{1}{2}(q-1)$ pairs lie on a $\Fq{q}$-line through $P_i$. There
are thus $\frac{1}{2}(q^2+1)$ pairs which do not lie on a line
through $P_i$ and, since there are two ways to label each pair,
we thus have $q^2+1$ choices for the two $\Fq{q}$-points. Finally, we shall
place $Q_1$ and $Q_2$ somewhere on $C$ but we have to make sure that the points
$P_1$, $P_2$, $P_3$, $Q_1$ and $Q_2$ are in general position.
Since the lines between $P_i$ and the other two $\Fq{q}$-points
intersect $C$ only in $\Fq{q}$-points, the only thing that might go
wrong when choosing $Q_1$ and $Q_2$ is that the line through
$Q_1$ and $Q_2$ might also go through $P_i$. As seen above,
there are exactly $q-1$ choices for $Q_1$ and $Q_2$ for which
this happens, so the remaining $q^2-q-(q-1)=q^2-2q+1$ choices
will give a configuration of the desired type. We thus have that
the number of elements in $\Delta_{l,1}^i \cap \Delta_c$ such that
$P_i$ lies on the outside of $C$ is
\begin{equation*}
 \frac{1}{2}(q^5-q^2)(q+1)q(q-1)(q^2+1)(q^2-2q+1).
\end{equation*}

We now assume that $P_i$ lies on the inside of $C$. We proceed
similarly to the above. First we observe that the number of $\Fq{q}$-lines
through $P_i$ is $q+1$ of which half intersect $C$ in two $\Fq{q}$-points
and half intersect $C$ in conjugate pairs of $\Fq{q^2}$-points. We choose
a line which intersects $C$ in two conjugate $\Fq{q^2}$-points and label
the intersection points by $R_1$ and $R_2$. We now choose
a $\Fq{q}$-point $P_j$ on $C$ in one of $q+1$ ways. The line through $P_i$
and $P_j$ intersects $C$ in another $\Fq{q}$-point and we choose
the final $\Fq{q}$-point away from this intersection point and $P_j$.
Finally, we shall place the points $Q_1$ and $Q_2$ on $C$ in a 
way so that the points $P_1$, $P_2$, $P_3$, $Q_1$ and $Q_2$ are
in general position. As above, the only thing that might go wrong is that
the line through $Q_1$ and $Q_2$ might go through $P_i$ and there
are precisely $q+1$ choices for $Q_1$ and $Q_2$ for which this happens.
Thus, there are $q^2-q-(q+1)=q^2-2q-1$ valid choices for $Q_1$ and
$Q_2$. Hence, there are 
\begin{equation*}
 \frac{1}{2}(q^5-q^2)(q^2-q)(q+1)(q+1)(q-1)(q^2-2q-1)
\end{equation*}
elements in $\Delta_{l,1}^i \cap \Delta_c$ such that
$P_i$ lies on the inside of $C$.

To compute the intersection $\Delta_{l,3}^{i,j} \cap \Delta_c$
we note that if we place $P_j$ outside of $C$ and then
choose two $\Fq{q}$-points on $C$ and two conjugate $\Fq{q^2}$-points
$Q_1$ and $Q_2$ on $C$ such that $P_1$, $P_2$, $P_3$, $Q_1$ and
$Q_2$ are in general position, then we must choose
$R_i$ as the other intersection point of $C$ with the line through 
$Q_1$ and $P_j$. We may thus use constructions analogous to those
above to see that there are
\begin{equation*}
 \frac{1}{2}(q^5-q^2)(q+1)q(q^2+1)(q^2-2q+1)
\end{equation*}
elements in $\Delta_{l,1}^{i,j} \cap \Delta_c$ with $P_j$ on
the outside of $C$ and
\begin{equation*}
 \frac{1}{2}(q^5-q^2)(q^2-q)(q+1)(q-1)(q^2-2q-3)
\end{equation*}
elements with $P_j$ on the inside of $C$. 

We may now put all the pieces together to obtain
\begin{equation*}
 \left|\left( \Pts \right)^{\Frob \sigma} \right| = q^6-7q^5+10q^4+15q^3-26q^2-8q+15.
\end{equation*}

 \subsection{The case \texorpdfstring{$\lambda=[1^5,2]$}{[1,1,1,1,1,2]}} 
 Throughout this section, $\lambda$ shall mean the partition 
 $[1^5,2]$. We shall denote the $\Fq{q}$-points by $P_1$, $P_2$, $P_3$, $P_4$ and $P_5$
 and the points of the conjugate pair of $\Fq{q^2}$-points by $Q_1$ and $Q_2$.
Let $U \subset (\Pn{2})^7$ be the subset
consisting of septuples of points such that the first five points lie
in general position.

 If three points of a conjugate $\lambda$-tuple in $U(\lambda)$ lie
 on a line, then $Q_1$ and $Q_2$ lie on a line passing
 through one of the $\Fq{q}$-points. 
There are
\begin{equation*}
 (q+1)+q+(q-1)+(q-2)+(q-3)=5q-5,
\end{equation*}
$\Fq{q}$-lines passing through $P_1$, $P_2$, $P_3$, $P_4$ or $P_5$
(or possibly two of them). Each of these lines contains $q^2-q$
conjugate pairs and no conjugate pair lies on two such lines.
We thus have
\begin{equation*}
 |\Delta_{l}| = (q^2+q+1)(q^2+q)q^2(q^2-2q+1)(q^2-5q+6)(5q-5)(q^2-q).
\end{equation*}

 If six of the points of a
 conjugate $\lambda$-tuple lie on a smooth conic $C$, then
 $C$ is defined over $\Fq{q}$ and contains $Q_1$ and $Q_2$. 
Therefore, to compute the cardinality of $\Delta_{c}$ we first choose
a smooth conic $C$ defined over $\Fq{q}$ and one of the points
$P_1$, $P_2$, $P_3$, $P_4$ or $P_5$ to possibly lie outside $C$. We call the
chosen point $P$.
Then, we choose four $\Fq{q}$-points and a conjugate pair on $C$.
Finally, we choose $P$ away from the six lines through pairs of the other
four $\Fq{q}$-points. We thus get
\begin{equation*}
 5(q^5-q^2)(q+1)q(q-1)(q-2)(q^2-q)(q^2-5q+6).
\end{equation*}
In the above we have counted the $\lambda$-tuples with all seven
points on a conic five times. We therefore must take away
\begin{equation*}
 4(q^5-q^2)(q+1)q(q-1)(q-2)(q-3)(q^2-q),
\end{equation*}
in order to obtain $|\Delta_{c}|$.

To compute the size of the intersection $\Delta_{l} \cap \Delta_{c}$
we shall decompose this set into a disjoint disjoint union of five subsets $A_i$, $i=1,\ldots, 5$,
where $A_i$ consists of those tuples where $P_i$ does not lie on the
conic $C$ through the other six points. Each of the sets $A_i$ is then
decomposed further into a  union of the sets $A_i^{\mathrm{out}}$ and $A_i^{\mathrm{in}}$
where $A_{i}^{\mathrm{out}}$ consists of those tuples with $P_i$ on the outside of $C$
and $A_i^{\mathrm{in}}$ consists of those with $P_i$ on the inside of $C$.
Finally, we shall decompose $A_i^{\mathrm{out}}$ into a union of the three disjoint subsets:
\begin{itemize}
 \item the set $A_{i,0}^{\mathrm{out}}$ consisting of $\lambda$-tuples such that the tangent lines to
 $C$ passing through $P_i$ do not pass through any of the other points of the $\lambda$-tuple,
 \item the set $A_{i,1}^{\mathrm{out}}$ consisting of $\lambda$-tuples such that exactly one of the tangent lines to
 $C$ passing through $P_i$ pass through one of the other points of the $\lambda$-tuple,
 \item the set $A_{i,2}^{\mathrm{out}}$ consisting of $\lambda$-tuples such that both the tangent lines to
 $C$ passing through $P_i$ passes through  another point of the $\lambda$-tuple.
\end{itemize}

To compute $|A_i^{\mathrm{out}}|$,
we first choose a smooth conic $C$ defined over $\Fq{q}$ in $q^5-q^2$
ways and then a point $P_i$ outside $C$ in $\frac{1}{2}(q+1)q$ ways.
As seen many times before, there are exactly $\frac{1}{2}(q-1)$
lines through $P_i$ which are defined over $\Fq{q}$ and which intersect
$C$ in a conjugate pair of points. We pick such a line and label
the points $Q_1$ and $Q_2$ in one of two ways.
From this point on, the computations are a little bit different
for the three subsets of $A_i^{\mathrm{out}}$.

\vspace{7pt}
 \noindent \textbf{The subset $A_{i,0}^{\mathrm{out}}$.}
We shall now pick the other four $\Fq{q}$-points of the $\lambda$-tuple.
Since we should not pick points whose tangents pass through $P_i$, we have $q-1$ choices
for the first point. For the second point, we should stay away from
the tangent points, the first point and the other intersection point
of $C$ and the line through $P_i$ and the first point. Hence, we have $q-3$
choices. In a similar way, we see that we have $q-5$ choices for the
third point and $q-7$ for the fourth. Hence,
\begin{equation*}
 |A_{i,0}^{\mathrm{out}}| = \frac{1}{2}(q^5-q^2)(q+1)q(q-1)(q-1)(q-3)(q-5)(q-7).
\end{equation*}

\vspace{7pt}
 \noindent \textbf{The subset $A_{i,1}^{\mathrm{out}}$.}
We begin by choosing one of the four $\Fq{q}$-points to lie on a tangent
to $C$ passing through $P_i$ and then we pick the tangent it should lie
on. For the first of the remaining three points we now have $q-1$
choices and, similarly to the above case, we have $q-3$ choices
for the second and $q-5$ for the third. Thus,
\begin{equation*}
 |A_{i,1}^{\mathrm{out}}| = 4 \cdot 2 \cdot \frac{1}{2}(q^5-q^2)(q+1)q(q-1)(q-1)(q-3)(q-5).
\end{equation*}

\vspace{7pt}
 \noindent \textbf{The subset $A_{i,2}^{\mathrm{out}}$.}
We begin by choosing two of the four $\Fq{q}$-points to lie on tangents
to $C$ passing through $P_i$ and then we pick which point should lie on
which tangent. For the first of the remaining two points we now have $q-1$
choices and we then have $q-3$ choices for the second. Thus,
\begin{equation*}
 |A_{i,2}^{\mathrm{out}}| = \binom{4}{2} \cdot 2 \cdot \frac{1}{2}(q^5-q^2)(q+1)q(q-1)(q-1)(q-3).
\end{equation*}

It remains to compute $|A_i^{\mathrm{in}}|$. We first choose
a smooth conic $C$ defined over $\Fq{q}$ in $q^5-q^2$ ways and then
a point $P_i$ on the inside of $C$ in $\frac{1}{2}(q^2-q)$ ways.
We have already seen that there now are $\frac{1}{2}(q+1)$ lines
passing through $P_i$ which are defined over $\Fq{q}$ and which intersect
$C$ in a conjugate pair of points. We thus pick such a line
and label the intersection points by $Q_1$ and $Q_2$. Since any
$\Fq{q}$-line through $P_i$ will intersect $C$ in precisely two points,
we have $(q+1)(q-1)(q-3)(q-5)$ choices for the remaining four $\Fq{q}$-points
of the $\lambda$-tuple. We thus see that
\begin{equation*}
 |A_i^{\mathrm{in}}| = \frac{1}{2}(q^5-q^2)(q^2-q)(q+1)(q-1)(q-3)(q-5).
\end{equation*}
We now conclude that
\begin{equation*}
 \left|\left( \Pts \right)^{\Frob \sigma} \right| = q^6-15q^5+90q^4-265q^3+374q^2-200q+15.
\end{equation*}

  \subsection{The case \texorpdfstring{$\lambda=[1^7]$}{[1,1,1,1,1,1,1]}} 
  Throughout this section, $\lambda$ shall mean the partition
  $[1^7]$.
  Since we shall almost exclusively be interested in objects
  defined over $\Fq{q}$, we shall often omit the decoration ``$\Fq{q}$''.
  For instance, we shall simply write ``point'' to mean ``$\Fq{q}$-point''.
  Let 
  $U \subset (\Pn{2})^7$ be the subset consisting of septuples of points such that 
  the first four points lie in general position.
  We thus have
  \begin{equation*}
   |U(\lambda)| = (q^2+q+1)(q^2+q)q^2(q^2-2q+1)(q^2+q-3)(q^2+q-4)(q^2+q-5).
  \end{equation*}

  The following notation will be quite convenient.

  \pagebreak[2]
  \begin{definition}
   If $P$ and $Q$ are two points in $\mathbb{P}^2$, then the line through $P$
   and $Q$ shall be denoted $\gline{P}{Q}$.
  \end{definition}
  \pagebreak[2]  
  
  Since we shall often want to stay away from lines through two of the first four
  points we define
  \begin{equation*}
   \mathscr{S} = \bigcup_{1\leq i < j \leq 4} \gline{P_i}{P_j}.
  \end{equation*}
  We note that $\mathscr{S}$ contains
  \begin{equation*}
   6(q-2)+4+3=6q-5
  \end{equation*}
  points.

  \subsubsection{The set \texorpdfstring{$\Delta_l$}{with three points on a line}} 

  The set $\Delta_l$ decomposes into a disjoint union of three sets
  \begin{equation*}
   \Delta_l = \Delta_{l,1} \cup \Delta_{l,2} \cup \Delta_{l,3},
  \end{equation*}
  where
 \begin{itemize}
  \item the points of $\Delta_{l,1}$ are such that at least one of the points $P_5$, $P_6$ or $P_7$ lies in $\mathscr{S}$,
  \item the points of $\Delta_{l,2}$ are such that one of the lines $\gline{P_i}{P_j}$, $5 \leq i < j \leq 7$, contains
  one of the points $P_1$, $P_2$, $P_3$ and $P_4$, but $\{P_5,P_6,P_7\} \cap \mathscr{S} = \emptyset$, and
  \item the points of $\Delta_{l,3}$ are such that the three points $P_5$, $P_6$ and $P_7$ lie on a line which does not pass
  through $P_1$, $P_2$, $P_3$ or $P_4$.
 \end{itemize}
 We shall consider the three subsets separately.
 
 \vspace{7pt}
 \noindent \textbf{The set $\Delta_{l,1}$.}
 For each subset $I \subset \{5,6,7\}$, let $\Delta_{l,1}(I)$
 denote the set of points in $\Delta_{l,1}$ such that
 $P_i \in \mathscr{S}$ for all $i \in I$. We can then
 decompose $\Delta_{l,1}$ further as
 \begin{equation*}
  \Delta_{l,1} = \Delta_{l,1}(\{5\}) \cup \Delta_{l,1}(\{6\}) \cup \Delta_{l,1}(\{7\}).
 \end{equation*}
 Clearly, $\Delta_{l,1}(\{i\}) \cap \Delta_{l,1}(\{j\}) = \Delta_{l,1}(\{i,j\})$.

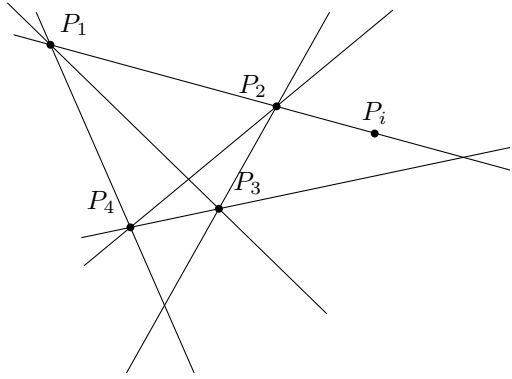
\begin{figure}[!htbp]
\begin{center}
\begin{tikzpicture}[scale=0.3]
 \def\ptsize{5pt}
 \draw (5,0) coordinate (a_1)--(14,16) coordinate (a_2);
 \draw (8,0) coordinate (b_1)--(1,16) coordinate (b_2);
 \draw (0,15) coordinate (c_1)--(22,9) coordinate (c_2);
 \draw (3,6) coordinate (d_1)--(22,10) coordinate (d_2);
 \coordinate (p_1) at (intersection of b_1--b_2 and c_1--c_2);
 \coordinate (p_2) at (intersection of a_1--a_2 and c_1--c_2);
 \coordinate (p_3) at (intersection of a_1--a_2 and d_1--d_2);
 \coordinate (p_4) at (intersection of b_1--b_2 and d_1--d_2);
 \fill[black] (p_1) circle (\ptsize);
 \fill[black] (p_2) circle (\ptsize);
 \fill[black] (p_3) circle (\ptsize);
 \fill[black] (p_4) circle (\ptsize);
 \draw [shorten <=-0.8cm, shorten >=-2cm] (p_1) -- (p_3);
 \draw [shorten <=-0.8cm, shorten >=-2cm] (p_4) -- (p_2);
 \coordinate (v) at (intersection of p_1--p_2 and p_3--p_4);
 \node [above right] at (p_1) {$P_1$};
 \node [above left] at (p_2) {$P_2$};
 \node [above right=0.1cm] at (p_3) {$P_3$};
 \node [above left=0.1cm] at (p_4) {$P_4$};
 \coordinate (p_i) at (16,10.62);
 \fill[black] (p_i) circle (\ptsize);
 \node [above] at (p_i) {$P_i$};
\end{tikzpicture}
\end{center}
\caption{A typical element of $\Delta_{l,1}(\{i\})$.}
\label{caseA}
\end{figure}

A typical element of $\Delta_{l,1}(\{i\})$ is illustrated in Figure \ref{caseA} above.
To compute $|\Delta_{l,1}(\{i\})|$ we first place the first four points in general position,
then choose $P_i$ as any point in $\mathscr{S}$ and finally place the
remaining two points anywhere. Hence
\begin{equation*}
 |\Delta_{l,1}(\{i\})| = \underbrace{(q^2+q+1)(q^2+q)q^2(q^2-2q+1)}_{|\pglt|}(6q-9)(q^2+q-4)(q^2+q-5).
\end{equation*}
Similarly, we have
\begin{equation*}
 |\Delta_{l,1}(\{i,j\})| = |\pglt| \cdot (6q-9)(6q-10)(q^2+q-5),
\end{equation*}
and
\begin{equation*}
 |\Delta_{l,1}(\{5,6,7\})| = |\pglt| \cdot (6q-9)(6q-10)(6q-11).
\end{equation*}
This allows us to compute $|\Delta_{l,1}|$ as
\begin{equation*}
 |\Delta_{l,1}|=|\pglt| \cdot (18q^5-99q^4+252q^3-414q^2+417q-180).
\end{equation*}

\begin{figure}[!htbp]
\begin{center}
\begin{tikzpicture}[scale=0.3]
 \def\ptsize{5pt}
 \def\sptsize{3pt}
 \draw [gray, dashed] (5,0) coordinate (a_1)--(14,16) coordinate (a_2);
 \draw [gray, dashed] (8,0) coordinate (b_1)--(1,16) coordinate (b_2);
 \draw [gray, dashed] (0,15) coordinate (c_1)--(22,9) coordinate (c_2);
 \draw [gray, dashed] (3,6) coordinate (d_1)--(22,10) coordinate (d_2);
 \coordinate (p_1) at (intersection of b_1--b_2 and c_1--c_2);
 \coordinate (p_2) at (intersection of a_1--a_2 and c_1--c_2);
 \coordinate (p_3) at (intersection of a_1--a_2 and d_1--d_2);
 \coordinate (p_4) at (intersection of b_1--b_2 and d_1--d_2);
 \coordinate (h) at (16,7);
 \coordinate (h11) at (7.5,0);
 \coordinate (h12) at (7.5,15);
 \coordinate (h21) at (12,0);
 \coordinate (h22) at (12,15);
 \fill[black] (p_1) circle (\ptsize);
 \fill[black] (p_2) circle (\sptsize);
 \fill[black] (p_3) circle (\sptsize);
 \fill[black] (p_4) circle (\sptsize);
 \draw [gray, dashed, shorten <=-0.8cm, shorten >=-2cm] (p_1) -- (p_3);
 \draw [gray, dashed, shorten <=-0.8cm, shorten >=-2cm] (p_4) -- (p_2);
 \coordinate (v) at (intersection of p_1--p_2 and p_3--p_4);
 \draw [shorten <=-0.8cm, shorten >=-2cm] (p_1) -- (h);
 \coordinate (p_i) at (intersection of h11--h12 and p_1--h);
 \coordinate (p_j) at (intersection of h21--h22 and p_1--h);
 \fill[black] (p_i) circle (\ptsize);
 \fill[black] (p_j) circle (\ptsize);
 \node [above] at (p_i) {$P_i$};
 \node [above] at (p_j) {$P_j$};
 \node [below] at (p_1) {$P_r$};
\end{tikzpicture}
\end{center}
\caption{A typical element of $\Delta_{l,2}^r(\{i,j\})$.}
\label{caseB}
\end{figure}
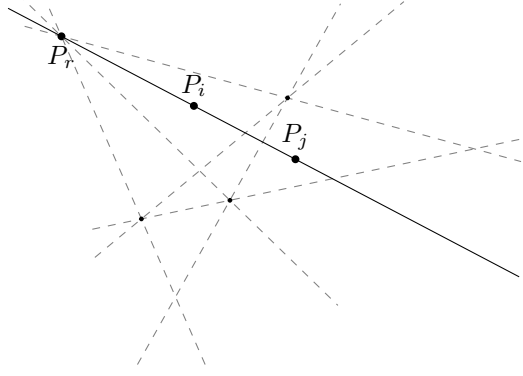

\vspace{7pt}
\noindent \textbf{The set $\Delta_{l,2}$.}
Let $\{i,j \} \in \{5,6,7\}$, $r \in \{1,2,3,4\}$ and let $\Delta_{l,2}^r(\{i,j\})$ be the
subset of points in $\Delta_{l,2}$ such that $\gline{P_i}{P_j} \cap \{P_1,P_2,P_3,P_4\}= \{P_r\}$. 
We also define
\begin{equation*}
 \Delta_{l,2}(\{i,j\}) = \bigcup_{r=1}^4 \Delta_{l,2}^r(\{i,j\}).
\end{equation*}
A typical element of $\Delta_{l,2}^r(\{i,j\})$ is illustrated in Figure
\ref{caseB}. To obtain an element of $\Delta_{l,2}^r(\{i,j\})$ we first
place $P_1$, $P_2$, $P_3$ or $P_4$ in general position.  There are
$q+1$ lines through $P_r$ of which $3$ are contained in $\mathscr{S}$.
We choose $\gline{P_i}{P_j}$ as one of the remaining $q-2$ lines.
Note that $\gline{P_i}{P_j}$ will not pass through any of the points
\begin{equation}
 Q_1 = \gline{P_1}{P_4} \cap \gline{P_2}{P_3}, \quad Q_2 = \gline{P_2}{P_4} \cap \gline{P_1}{P_3}, \quad Q_3 = \gline{P_3}{P_4} \cap \gline{P_1}{P_2}.
 \label{q1q2q3}
\end{equation}
Hence, $\gline{P_i}{P_j}$ will intersect $\mathscr{S}$ in $P_r$
and three further points. There are thus $q-3$ ways to choose
$P_i$ and then $q-4$ ways to choose $P_j$. Finally, there
are 
\begin{equation*}
|\mathbb{P}^2 \setminus \mathscr{S}|-2=q^2+q+1-(6q-5)-2=q^2-5q+4,
\end{equation*}
choices for the seventh point. We thus have
\begin{equation*}
 |\Delta_{l,2}^r(\{i,j\})| = |\pglt| \cdot (q-2)(q-3)(q-4)(q^2-5q+4).
\end{equation*}

We have counted some tuples several times. To begin
with, the points of
\begin{equation*}
 \Delta_{l,2}^r(\{5,6\}) \cap \Delta_{l,2}^r(\{5,7\}) \cap \Delta_{l,2}^r(\{6,7\}),
\end{equation*}
have been counted three times. There are
\begin{equation*}
 |\pglt| \cdot (q-2)(q-3)(q-4)(q-5),
\end{equation*}
of these. 

\begin{figure}[!htbp]
\begin{center}
\begin{tikzpicture}[scale=0.3]
 \def\ptsize{5pt}
 \def\sptsize{3pt}
 \draw [gray, dashed] (5,0) coordinate (a_1)--(14,16) coordinate (a_2);
 \draw [gray, dashed] (8,0) coordinate (b_1)--(1,16) coordinate (b_2);
 \draw [gray, dashed] (0,15) coordinate (c_1)--(22,9) coordinate (c_2);
 \draw [gray, dashed] (3,6) coordinate (d_1)--(22,10) coordinate (d_2);
 \coordinate (p_1) at (intersection of b_1--b_2 and c_1--c_2);
 \coordinate (p_2) at (intersection of a_1--a_2 and c_1--c_2);
 \coordinate (p_3) at (intersection of a_1--a_2 and d_1--d_2);
 \coordinate (p_4) at (intersection of b_1--b_2 and d_1--d_2);
 \fill[black] (p_1) circle (\sptsize);
 \fill[black] (p_2) circle (\ptsize);
 \fill[black] (p_3) circle (\ptsize);
 \fill[black] (p_4) circle (\sptsize);
 \draw [gray, dashed, shorten <=-0.8cm, shorten >=-2cm] (p_1) -- (p_3);
 \draw [gray, dashed, shorten <=-0.8cm, shorten >=-2cm] (p_4) -- (p_2);
 \coordinate (v) at (intersection of p_1--p_2 and p_3--p_4);
 \coordinate (h_2) at (2,9);
 \draw [shorten <=-0.8cm, shorten >=-1.4cm] (p_3) -- (h_2);
 \coordinate (h_3) at (2,10);
 \draw [shorten <=-0.8cm, shorten >=-1.2cm] (p_2) -- (h_3);
 \coordinate (p_i) at (intersection of p_3--h_2 and p_2--h_3);
 \coordinate (p_j1) at (2.7,8.8);
 \coordinate (p_j2) at (2.3,10.05);
 \fill[black] (p_i) circle (\ptsize);
 \fill[black] (p_j1) circle (\ptsize);
 \fill[black] (p_j2) circle (\ptsize);
 \node [above] at (p_i) {$P_i$};
 \node [below] at (p_j1) {$P_{j_1}$};
 \node [above] at (p_j2) {$P_{j_2}$};
 \node [above] at (p_2) {$P_s$};
 \node [below] at (p_3) {$P_r$};
\end{tikzpicture}
\end{center}
\caption{A typical element of $\Delta_{l,2}^r(\{i,j\}) \cap \Delta_{l,2}^s(\{i,k\})$.}
\label{caseBb}
\end{figure}
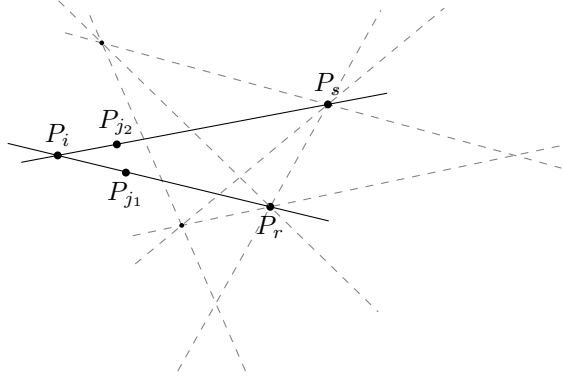

Further, the sets $\Delta_{l,2}^r(\{i,j\})$ and
$\Delta_{l,2}^s(\{i,k\})$ will intersect if $r \neq s$ and $j \neq k$.  A typical
element is illustrated in Figure \ref{caseBb}.

To compute $|\Delta_{l,2}^r(\{i,j\}) \cap \Delta_{l,2}^s(\{i,k\})|$
we begin by choosing $P_1$, $P_2$, $P_3$ and $P_4$ in general position and continue by
choosing $P_i$ outside $\mathscr{S}$ in $q^2-5q+6$ ways.
This gives us two lines $\gline{P_i}{P_r}$ and $\gline{P_i}{P_s}$ which
intersect $\mathscr{S}$ in four points each. We choose $P_{j}$ on
$\gline{P_i}{P_r}$ away from $P_i$ and $\mathscr{S}$ in $q-4$ ways
and similarly for $P_{k}$. This gives
\begin{equation*}
 |\Delta_{l,2}^r(\{i,j\}) \cap \Delta_{l,2}^s(\{i,k\})| = |\pglt| \cdot (q^2-5q+6) (q-4)^2.
\end{equation*}

\begin{figure}[!htbp]
\begin{center}
\begin{tikzpicture}[scale=0.3]
 \def\ptsize{5pt}
 \def\sptsize{3pt}
 \draw [gray, dashed] (5,0) coordinate (a_1)--(14,16) coordinate (a_2);
 \draw [gray, dashed] (8,0) coordinate (b_1)--(1,16) coordinate (b_2);
 \draw [gray, dashed] (0,15) coordinate (c_1)--(22,9) coordinate (c_2);
 \draw [gray, dashed] (3,6) coordinate (d_1)--(22,10) coordinate (d_2);
 \coordinate (p_1) at (intersection of b_1--b_2 and c_1--c_2);
 \coordinate (p_2) at (intersection of a_1--a_2 and c_1--c_2);
 \coordinate (p_3) at (intersection of a_1--a_2 and d_1--d_2);
 \coordinate (p_4) at (intersection of b_1--b_2 and d_1--d_2);
 \fill[black] (p_1) circle (\sptsize);
 \fill[black] (p_2) circle (\ptsize);
 \fill[black] (p_3) circle (\ptsize);
 \fill[black] (p_4) circle (\ptsize);
 \draw [gray, dashed, shorten <=-0.8cm, shorten >=-2cm] (p_1) -- (p_3);
 \draw [gray, dashed, shorten <=-0.8cm, shorten >=-2cm] (p_4) -- (p_2);
 \coordinate (v) at (intersection of p_1--p_2 and p_3--p_4);
 \coordinate (h_1) at (2,11);
 \draw [shorten <=-0.8cm, shorten >=-2cm] (p_4) -- (h_1);
 \coordinate (h_2) at (2,9);
 \draw [shorten <=-0.8cm, shorten >=-1.4cm] (p_3) -- (h_2);
 \coordinate (h_3) at (2,10);
 \draw [shorten <=-0.8cm, shorten >=-1.2cm] (p_2) -- (h_3);
 \coordinate (p_5) at (intersection of p_4--h_1 and p_3--h_2);
 \coordinate (p_6) at (intersection of p_4--h_1 and p_2--h_3);
 \coordinate (p_7) at (intersection of p_3--h_2 and p_2--h_3);
 \fill[black] (p_5) circle (\ptsize);
 \fill[black] (p_6) circle (\ptsize);
 \fill[black] (p_7) circle (\ptsize);
 \node [below] at (p_5) {$P_7$};
 \node [above] at (p_6) {$P_6$};
 \node [above] at (p_7) {$P_5$};
 \node [above] at (p_2) {$P_r$};
 \node [below] at (p_3) {$P_s$};
 \node [below] at (p_4) {$P_t$};
\end{tikzpicture}
\end{center}
\caption{A typical element of $\Delta_{l,2}^r(\{5,6\}) \cap \Delta_{l,2}^s(\{5,7\}) \cap \Delta_{l,2}^t(\{6,7\})$.}
\label{caseBc}
\end{figure}
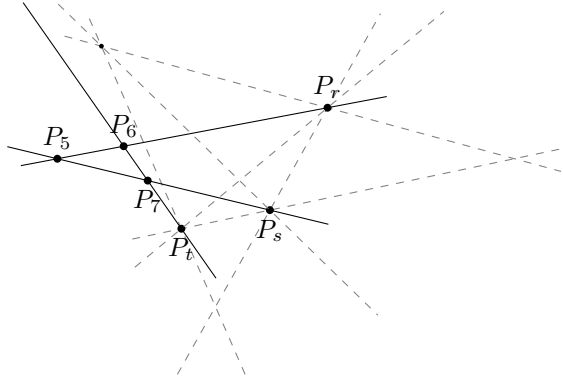

Finally, we must compute the cardinality of the triple intersection
\begin{equation*}
 \Delta_{l,2}^r(\{5,6\}) \cap \Delta_{l,2}^s(\{5,7\}) \cap \Delta_{l,2}^t(\{6,7\}),
\end{equation*}
where $r$, $s$ and $t$ are distinct.
A typical element of the intersection is illustrated in Figure
\ref{caseBc}.

This is where we have to pay for the awkward requirement that 
$P_5$, $P_6$ and $P_7$ should not be in $\mathscr{S}$. We shall
view $\Delta_{l,2}^r(\{5,6\}) \cap \Delta_{l,2}^s(\{5,7\}) \cap \Delta_{l,2}^t(\{6,7\})$
as an open subset of the set $T^{r,s,t}$ consisting of
tuples such that 
\begin{itemize}
 \item the line $\gline{P_5}{P_6}$ passes through $P_r$,
$\gline{P_5}{P_7}$ passes through $P_s$ and $\gline{P_6}{P_7}$ passes
through $P_t$ but,
\item we allow $P_5$, $P_6$ and $P_7$ to lie
in $\mathscr{S}$, but,
\item we do not allow the lines $\gline{P_i}{P_j}$, $5 \leq i < j \leq 7$ to be contained
in $\mathscr{S}$. 
\end{itemize}
The complement of $\Delta_{l,2}^r(\{5,6\}) \cap \Delta_{l,2}^s(\{5,7\}) \cap \Delta_{l,2}^t(\{6,7\})$ in
$T^{r,s,t}$ can be decomposed
into a union of three subsets $T^{r,s,t}_i$, $i=5,6,7$,
consisting of those tuples with $P_i$ in $\mathscr{S}$. 

We begin with the computation of $|T^{r,s,t}|$. To obtain such a tuple,
we begin by choosing a line $L_r$ through $P_r$ in $q-2$ ways.
We shall then choose a line $L_s$ through $P_s$. There are however
two cases that may occur. Typically, the intersection point
$P_5=L_r \cap L_s$ will lie outside $\mathscr{S}$ but for one
choice of $L_s$ it will lie in $\mathscr{S}$. The situation is
illustrated in Figure \ref{speccase}.

\begin{figure}[!htbp]
\begin{center}
\begin{tikzpicture}[scale=0.3]
 \def\ptsize{5pt}
 \def\sptsize{3pt}
 \draw [gray, dashed] (5,0) coordinate (a_1)--(14,16) coordinate (a_2);
 \draw [gray, dashed] (8,0) coordinate (b_1)--(1,16) coordinate (b_2);
 \draw [gray, dashed] (0,15) coordinate (c_1)--(22,9) coordinate (c_2);
 \draw [gray, dashed] (3,6) coordinate (d_1)--(22,10) coordinate (d_2);
 \coordinate (p_1) at (intersection of b_1--b_2 and c_1--c_2);
 \coordinate (p_2) at (intersection of a_1--a_2 and c_1--c_2);
 \coordinate (p_3) at (intersection of a_1--a_2 and d_1--d_2);
 \coordinate (p_4) at (intersection of b_1--b_2 and d_1--d_2);
 \fill[black] (p_1) circle (\sptsize);
 \fill[black] (p_2) circle (\ptsize);
 \fill[black] (p_3) circle (\ptsize);
 \fill[black] (p_4) circle (\sptsize);
 \draw [gray, dashed, shorten <=-0.8cm, shorten >=-2cm] (p_1) -- (p_3);
 \draw [gray, dashed, shorten <=-0.8cm, shorten >=-2cm] (p_4) -- (p_2);
 \coordinate (v) at (intersection of p_1--p_2 and p_3--p_4);
 \coordinate (h_2) at (2,10);
 \draw [shorten <=-0.8cm, shorten >=-1.4cm] (p_3) -- (h_2);
 \coordinate (h_3) at (2,8.63);
 \draw [dotted,shorten <=-0.8cm, shorten >=-1.2cm] (p_2) -- (h_3);
 \coordinate (k_2) at (2,11);
 \draw [shorten <=-0.8cm, shorten >=-1.2cm] (p_2) -- (k_2);
 \coordinate (p_i) at (intersection of p_3--h_2 and p_2--h_3);
 \coordinate (p_j) at (intersection of p_3--h_2 and p_2--k_2);
 \fill[black] (p_i) circle (\sptsize);
 \fill[black] (p_j) circle (\ptsize);
 \node [above] at (p_2) {$P_s$};
 \node [below] at (p_3) {$P_r$};
\end{tikzpicture}
\end{center}
\caption{}
\label{speccase}
\end{figure}

There are $q-3$ ways to
choose $L_s$ so that $L_r \cap L_s$ lies outside $\mathscr{S}$.
When we choose the line $L_t$ through $P_t$ we must make sure
that $L_t$ is not contained in $\mathscr{S}$ and that $L_t$
does not pass through $L_r \cap L_s$, since we want to end up with
three distinct intersection points. We thus have $q-3$ choices.
On the other hand, if we choose $L_s$ as the one line making
the intersection point $L_r \cap L_s$ lie in $\mathscr{S}$ we
only need to make sure that $L_t$ is not contained in $\mathscr{S}$
and we thus have $q-2$ choices. Hence, we see that
\begin{equation*}
 |T^{r,s,t}| = |\pglt| \cdot \left( (q-2)(q-3)^2 + (q-2)^2 \right).
\end{equation*}

We now turn to the computation of $|T^{r,s,t}_i|$, $i=5,6,7$. 
We then begin by choosing a line $L_r$ through $P_r$
in $q-2$ ways. The line $L_s$ through $P_s$ is then completely
determined since we must have $P_i \in \mathscr{S}$. This
gives us $q-2$ choices for the final line $L_t$ through $P_t$. Hence,
\begin{equation*}
 |T^{r,s,t}_i|=|\pglt| \cdot (q-2)^2.
\end{equation*}

We now turn to the computation of $|T^{r,s,t}_{i} \cap T^{r,s,t}_j|$, $5 \leq i < j \leq 7$.
As above, we begin by choosing
a line $L_r$ through $P_r$ in $q-2$ ways. Since $P_i$ must lie
in $\mathscr{S}$ we have only one choice for $L_s$. Since $P_j = L_s \cap L_t$
we see that we now have precisely one choice for $L_t$ also. Hence,
\begin{equation*}
 |T^{r,s,t}_{i} \cap T^{r,s,t}_j| = |\pglt| \cdot (q-2).
\end{equation*}

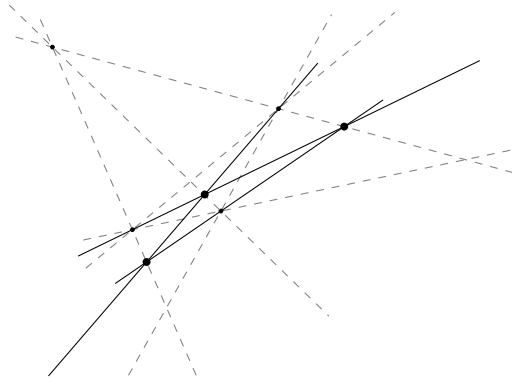
\begin{figure}[!htbp]
\begin{center}
\begin{tikzpicture}[scale=0.3]
 \def\ptsize{5pt}
 \def\sptsize{3pt}
 \draw [gray, dashed] (5,0) coordinate (a_1)--(14,16) coordinate (a_2);
 \draw [gray, dashed] (8,0) coordinate (b_1)--(1,16) coordinate (b_2);
 \draw [gray, dashed] (0,15) coordinate (c_1)--(22,9) coordinate (c_2);
 \draw [gray, dashed] (3,6) coordinate (d_1)--(22,10) coordinate (d_2);
 \coordinate (p_1) at (intersection of b_1--b_2 and c_1--c_2);
 \coordinate (p_2) at (intersection of a_1--a_2 and c_1--c_2);
 \coordinate (p_3) at (intersection of a_1--a_2 and d_1--d_2);
 \coordinate (p_4) at (intersection of b_1--b_2 and d_1--d_2);
 \fill[black] (p_1) circle (\sptsize);
 \fill[black] (p_2) circle (\sptsize);
 \fill[black] (p_3) circle (\sptsize);
 \fill[black] (p_4) circle (\sptsize);
 \draw [gray, dashed, shorten <=-0.8cm, shorten >=-2cm] (p_1) -- (p_3);
 \draw [gray, dashed, shorten <=-0.8cm, shorten >=-2cm] (p_4) -- (p_2);
 \coordinate (v) at (intersection of p_1--p_2 and p_3--p_4);
 \coordinate (h_2) at (11.6,9);
 \draw [shorten <=-1.7cm, shorten >=-1.7cm] (p_3) -- (h_2);
 \coordinate (h_1) at (intersection of b_1--b_2 and p_3--h_2);
 \draw [shorten <=-0.8cm, shorten >=-2cm] (p_2) -- (h_1);
 \coordinate (h_3) at (intersection of c_1--c_2 and p_3--h_2);
 \draw [shorten <=-0.8cm, shorten >=-2cm] (p_4) -- (h_3);
 \coordinate (p_5) at (intersection of p_4--h_1 and p_3--h_2);
 \coordinate (p_6) at (intersection of p_4--h_3 and p_2--h_1);
 \coordinate (p_7) at (intersection of p_3--h_2 and p_2--h_3);
 \fill[black] (p_5) circle (\ptsize);
 \fill[black] (p_6) circle (\ptsize);
 \fill[black] (p_7) circle (\ptsize);
\end{tikzpicture}
\end{center}
\caption{The only element in $T^{r,s,t}_5 \cap T^{r,s,t}_6 \cap T^{r,s,t}_7$.}
\label{thesingularity}
\end{figure}

We now consider $T^{r,s,t}_5 \cap T^{r,s,t}_6 \cap T^{r,s,t}_7$.
It turns out that once the four points $P_1$, $P_2$, $P_3$ and $P_4$
have been placed in general position, there is precisely one such tuple.
The situation is illustrated in Figure \ref{thesingularity}.

This finally allows us to compute
\begin{equation*}
 \Delta_{l,2} = |\pglt| \cdot (12q^5-212q^4+1504q^3-5320q^2+9296q-6360).
\end{equation*}

\vspace{7pt}
\noindent \textbf{The set $\Delta_{l,3}$.}
Recall the definition of the three points $Q_1$, $Q_2$ and $Q_3$
from Equation \ref{q1q2q3}. Using these three points we may decompose
$\Delta_{l,3}$ into a disjoint union of the following
subsets:
\begin{itemize}
 \item $\Delta_{l,3}(\{Q_r,Q_s\})$ consisting of those tuples of $\Delta_{l,3}$
 where $P_5$, $P_6$ and $P_7$ lie on the line $\gline{Q_r}{Q_s}$, $1 \leq r < s \leq 3$, and,
 \item $\Delta_{l,3}(\{Q_r\})$ consisting of those tuples of $\Delta_{l,3}$ with
 $P_5$, $P_6$ and $P_7$ on a line through $Q_r$, $1 \leq r \leq 3$, which does not pass
 through any of the other $Q_i$, and
 \item $\Delta_{l,3}(\emptyset)$ consisting of those tuples of $\Delta_{l,3}$
 with $P_5$, $P_6$ and $P_7$ on a line which does not pass through
 $Q_1$, $Q_2$ or $Q_3$.
\end{itemize}

We begin by considering $\Delta_{l,3}(\{Q_r,Q_s\})$. The line $\gline{Q_r}{Q_s}$
contains $q+1$ points of which four lie in $\mathscr{S}$. There are
thus $q-3$ choices for $P_5$, $q-4$ choices for $P_6$ and
$q-5$ choices for $P_7$. Hence,
\begin{equation*}
 |\Delta_{l,3}(\{Q_r,Q_s\})| = |\pglt| \cdot (q-3)(q-4)(q-5).
\end{equation*}
We continue with $|\Delta_{l,3}(\{Q_r\})|$. There are $q+1$ lines
through $Q_r$ of which two are contained in $\mathscr{S}$ and
two are the lines through the other two $Q_i$. Hence,
there are $q-3$ choices for a line $L$ though $Q_r$. The
line $L$ intersects $\mathscr{S}$ in five points so we have
$q-4$ choices for $P_5$, $q-5$ choices for $P_6$ and
$q-6$ choices for $P_7$. We conclude that
\begin{equation*}
 |\Delta_{l,3}(\{Q_r\})| = |\pglt| \cdot (q-3)(q-4)(q-5)(q-6).
\end{equation*}
To compute $|\Delta_{l,3}(\emptyset)|$ we begin by
choosing a line $L$ which does not pass through any of the points
$P_1$, $P_2$, $P_3$, $P_4$, $Q_1$, $Q_2$ and $Q_3$.
There are $q^2+q+1$ lines in $\mathbb{P}^2$, of which
$q+1$ passes through $P_i$, $i=1,2,3,4$. There is
exactly one line through each pair of these points so there
are
\begin{equation*}
 q^2+q+1-4(q+1)+6=q^2-3q+3
\end{equation*}
lines which do not pass through $P_1$, $P_2$, $P_3$ and $P_4$.
Of the $q+1$ lines through $Q_i$, $i=1,2,3$, precisely two have
been removed above and the line $\gline{Q_i}{Q_j}$ passes
through both $Q_i$ and $Q_j$. Hence, we have
\begin{equation*}
 q^2-3q+3-3(q-1)+3 = q^2-6q+9,
\end{equation*}
choices for $L$. 

The line $L$ intersects $\mathscr{S}$ in six points. We therefore have
$q-5$ choices for $P_5$, $q-6$ choices for $P_6$ and $q-7$ choices
for $P_7$. Hence,
\begin{equation*}
 |\Delta_{l,3}(\emptyset)|=|\pglt| \cdot (q^2-6q+9)(q-5)(q-6)(q-7).
\end{equation*}
We now add everything together to obtain
\begin{equation*}
 |\Delta_{l,3}| = |\pglt| \cdot (q^5-21q^4+173q^3-693q^2+1338q-990)
\end{equation*}
and, finally,
\begin{equation*}
 |\Delta_{l}|=|\pglt| \cdot (31q^5-332q^4+1929q^3-6427q^2+11051q-7530).
\end{equation*}

\subsubsection{The set \texorpdfstring{$\Delta_c$}{with six points on a conic}} 
We decompose $\Delta_c$ as
\begin{equation*}
 \Delta_c = \Delta_{c,1} \cup \Delta_{c,2},
\end{equation*}
where $\Delta_{c,1}$ consists of tuples where
six points lie on a smooth conic $C$ with one of the points
 $P_1$, $P_2$, $P_3$ or $P_4$ possibly outside $C$ and
 $\Delta_{c,2}$ consists of tuples where six points lie on a smooth conic $C$ with one of the points
 $P_5$, $P_6$ or $P_7$ possibly outside $C$.

To obtain an element of $\Delta_{c,1}$ we first choose one of
the points $P_1$, $P_2$, $P_3$ and $P_4$ and call it $P$. Then 
we choose a smooth conic $C$ in $q^5-q^2$ ways and place
all of the seven points except $P$ on $C$ in
\begin{equation*}
 (q+1)q(q-1)(q-2)(q-3)(q-4),
\end{equation*}
ways. There are three lines through pairs of points in $\{P_1,P_2,P_3,P_4\} \setminus \{P\}$
which together contain $3q$ points. These lines do not contain
$P_5$, $P_6$ and $P_7$ so we have
\begin{equation*}
 q^2+q+1-3q-3=q^2-2q-2,
\end{equation*}
choices for $P$. Multiplying everything together we obtain
\begin{equation*}
 N_1:=4(q^5-q^2)(q+1)q(q-1)(q-2)(q-3)(q-4)(q^2-2q-2),
\end{equation*}
which is almost $|\Delta_{c,1}|$ except that we have counted
the tuples where all seven points lie on $C$ four times.

To obtain an element of $\Delta_{c,2}$ we first
choose $P_5$, $P_6$ and $P_7$ and call the chosen point $P$.
We then choose a smooth conic $C$ and place all but the chosen points
on $C$. Finally, we place $P$ anywhere in $\mathbb{P}^2$
except at the six chosen points. In this way we obtain the number
\begin{equation*}
 N_2 := 3(q^5-q^2)(q+1)q(q-1)(q-2)(q-3)(q-4)(q^2+q-5),
\end{equation*}
which is almost equal to $|\Delta_{c,2}|$ except that we have counted
the tuples with all seven points on $C$ three times.

We now want to compute the number of tuples with all seven points
on a smooth conic $C$. We thus choose a smooth conic $C$ and
place all seven points on it in
\begin{equation*}
 N_7:=(q^5-q^2)(q+1)q(q-1)(q-2)(q-3)(q-4)(q-5),
\end{equation*}
ways. We thus have

\begin{equation*}
 |\Delta_{c}| = |\pglt| \cdot (7q^5-74q^4+288q^3-517q^2+446q-168)
\end{equation*}

\subsubsection{The set \texorpdfstring{$\Delta_l \cap \Delta_c$}{with three points on a line and six points on a conic}} 
We introduce the filtration $\mathscr{F}_3 \subset \mathscr{F}_2 \subset \mathscr{F}_1=\Delta_l \cap \Delta_c$ where
\begin{itemize}
 \item the set $\mathscr{F}_1$ consists  of tuples such that at least one line contains three points of the tuple,
 \item the set $\mathscr{F}_2$ consists  of tuples such that at least two lines contain three points of the tuple,
 \item the set $\mathscr{F}_3$ consists  of tuples such that at least three lines contain three points of the tuple.
\end{itemize}
The strategy will be to compute the numbers:
\begin{align*}
 N_1 & = |\mathscr{F}_1|+|\mathscr{F}_2|+|\mathscr{F}_3|, \\
 N_2 & = |\mathscr{F}_2|+2|\mathscr{F}_3|, \\
 N_3 & = |\mathscr{F}_3|,
\end{align*}
and thereby obtain the desired cardinality.

Since the points $P_1$, $P_2$, $P_3$ and $P_4$ are assumed
to constitute a frame, we must do things a little bit differently
depending on whether the point not on the conic is one of these
four or not. We therefore make further subdivisions.

\vspace{7pt}
\noindent \textbf{The subsets with $P_5$, $P_6$ or $P_7$ not on the conic.}
We shall denote the subsets in question by $\mathscr{F}_i^{5,6,7}$ and, similarly
\begin{align*}
 N_1^{5,6,7} & = |\mathscr{F}_1^{5,6,7}|+|\mathscr{F}_2^{5,6,7}|+|\mathscr{F}_3^{5,6,7}|, \\
 N_2^{5,6,7} & = |\mathscr{F}_2^{5,6,7}|+2|\mathscr{F}_3^{5,6,7}|, \\
 N_3^{5,6,7} & = |\mathscr{F}_3^{5,6,7}|.
\end{align*}
To compute $N_1^{5,6,7}$, we first choose one of the points $P_5$, $P_6$ or $P_7$ to be the point
$P$ not on the smooth conic $C$ and call the remaining two points $P_i$ and $P_j$. 
We then choose $C$ in $q^5-q^2$ ways and choose two points 
among $\{P_1,P_2,P_3,P_4,P_i,P_j\}$ and call them $R_1$ and $R_2$.
There are $(q+1)q$ ways to place $R_1$ and $R_2$ on $C$ and there
are then $q-1$ ways to place $P$ on the line $\gline{R_1}{R_2}$. Finally,
we place the remaining four points on $C$ in $(q-1)(q-2)(q-3)(q-4)$
ways. Multiplying everything together we obtain
\begin{equation*}
 N_1^{5,6,7} := 3 \cdot \binom{6}{2} \cdot (q^5-q^2)(q+1)q(q-1)^2(q-2)(q-3)(q-4).
\end{equation*}

In order to compute $N_2^{5,6,7}$, we first choose one of the points $P_5$, $P_6$ or $P_7$ to be the point
$P$ not on the smooth conic $C$ and call the remaining two points $P_i$ and $P_j$.
We then choose $C$ in $q^5-q^2$ ways and choose two unordered pairs of unordered points 
among $\{P_1,P_2,P_3,P_4,P_i,P_j\}$. This can be done in $\frac{1}{2} \cdot\binom{6}{4} \cdot \binom{4}{2}$ ways.
We call the points of the first pair  $R_1$ and $R_2$ and those of the second $O_1$ and $O_2$.
There are $(q+1)q(q-1)(q-2)(q-3)(q-4)$ ways to place $\{P_1,P_2,P_3,P_4,P_i,P_j\}$ on $C$
and the point $P$ is then completely determined as $P=\gline{R_1}{R_2} \cap \gline{O_1}{O_2}$.
Thus
\begin{equation*}
 N_2^{5,6,7}= 3 \cdot \frac{1}{2} \cdot \binom{6}{4} \cdot \binom{4}{2} \cdot (q^5-q^2)(q+1)q(q-1)(q-2)(q-3)(q-4).
\end{equation*}

The computation of $N_3^{5,6,7}$ is slightly more complicated since we need
to subdivide into two subcases depending on if $P$ is on the outside or on the inside
of $C$. We call the two corresponding numbers $N_{3,\mathrm{out}}^{5,6,7}$ and $N_{3,\mathrm{in}}^{5,6,7}$.

To compute $N_{3,\mathrm{out}}^{5,6,7}$ we first choose one of the points $P_5$, $P_6$ or $P_7$ to be the point
$P$ not on the smooth conic $C$. We
proceed by choosing the smooth conic $C$ in $q^5-q^2$ ways
and then the point $P$ on the outside of $C$ in $\frac{1}{2}(q+1)q$ ways.
We now place $P_1$ at one of the $q-1$ points of $C$ whose tangent does not
pass through $P$ and choose one of the remaining $5$ points as the other intersection
point in $C \cap \gline{P_1}{P}$. There are now four remaining points $P_i$, $P_j$, $P_k$ and $P_l$
to place on $C$. We place $P_i$ at one of the $q-3$ remaining points of $C$
whose tangent does not pass through $P$ and choose one of the remaining three points
as the other intersection point in $C \cap \gline{P_i}{P}$. There are now two points
$P_r$ and $P_s$ to place on $C$. We place $P_r$ at one of the $q-5$ possible points
and the point $P_s$ is then determined.
We thus have
\begin{equation*}
 N_{3,\mathrm{out}}^{5,6,7} = 3 \cdot (q^5-q^2) \cdot \frac{1}{2}(q+1)q \cdot (q-1) \cdot 5 \cdot (q-3) \cdot 3 \cdot (q-5).
\end{equation*}

We proceed by computing $N_{3,\mathrm{in}}^{5,6,7}$. We first choose one of the points $P_5$, $P_6$ or $P_7$ to be the point
$P$ not on the smooth conic $C$. We
proceed by choosing the smooth conic $C$ in $q^5-q^2$ ways
and then the point $P$ on the outside of $C$ in $\frac{1}{2}(q-1)q$ ways.

We now place $P_1$ at one of the $q+1$ points of $C$ whose tangent does not
pass through $P$ and choose one of the remaining $5$ points as the other intersection
point in $C \cap \gline{P_1}{P}$. There are now four remaining points $P_i$, $P_j$, $P_k$ and $P_l$
to place on $C$. We place $P_i$ at one of the $q-1$ remaining points of $C$
and choose one of the remaining three points
as the other intersection point in $C \cap \gline{P_i}{P}$. There are now two points
$P_r$ and $P_s$ to place on $C$. We place $P_r$ at one of the $q-3$ possible points
and the point $P_s$ is then determined.
We now
see that
\begin{equation*}
 N_{3,\mathrm{in}}^{5,6,7} = 3 \cdot (q^5-q^2) \cdot \frac{1}{2}(q-1)q \cdot (q+1) \cdot 5 \cdot (q-1) \cdot 3 \cdot (q-3).
\end{equation*}

\vspace{7pt}
\noindent \textbf{The subsets with $P_1$, $P_2$, $P_3$ or $P_4$ not on the conic.}
We shall denote the subsets in question by $\mathscr{F}_i^{1,2,3,4}$ and, similarly
\begin{align*}
 N_1^{1,2,3,4} & = |\mathscr{F}_1^{1,2,3,4}|+|\mathscr{F}_2^{1,2,3,4}|+|\mathscr{F}_3^{1,2,3,4}|, \\
 N_2^{1,2,3,4} & = |\mathscr{F}_2^{1,2,3,4}|+2|\mathscr{F}_3^{1,2,3,4}|, \\
 N_3^{1,2,3,4} & = |\mathscr{F}_3^{1,2,3,4}|.
\end{align*}
In order to compute $N_1^{1,2,3,4}$, we first choose one of the 
points $P_1$, $P_2$, $P_3$ or $P_4$ to be the point
$P$ not on the smooth conic $C$ and call the remaining three points $P_r$, $P_s$ and $P_t$.
We continue by choosing a smooth conic $C$ in $q^5-q^2$ ways.

We first assume that $P$ lies on a line $\gline{R_1}{R_2}$ where
$\{R_1, R_2 \} \subset \{P_5,P_6,P_7\}$. We therefore choose
the two points in $3$ ways and call the remaining point $P_i$.
We then place $R_1$ and $R_2$ on $C$ in $(q+1)q$ ways.
We continue by choosing the three points $P_r$, $P_s$ and
$P_t$ on $C$ in $(q-1)(q-2)(q-3)$ ways. The lines
$\gline{P_r}{P_s}$, $\gline{P_r}{P_t}$ and $\gline{P_s}{P_t}$
intersect the line $\gline{R_1}{R_2}$ in three distinct points 
so there are $q-4$ ways to choose the point $P$ on $\gline{R_1}{R_2}$
but away from these three points and $R_1$ and $R_2$. Finally, we place
$P_i$ at one of the $q-4$ remaining points of $C$. Multiplying everything
together we get
\begin{equation*}
 4 \cdot 3 \cdot (q^5-q^2)(q+1)q(q-1)(q-2)(q-3)(q-4)^2.
\end{equation*}

We now assume that $P$ lies on a line $\gline{a}{b}$ with
$a \in \{P_1,P_2,P_3,P_4\}$ and $b \in \{P_5,P_6,P_7\}$.
We thus first choose $a$ as one of the
points in $\{P_r,P_s,P_t\}$ and the point $b$ as
one of the points $\{P_5,P_6,P_7\}$ and place $a$ and
$b$ on $C$ in one of $(q+1)q$ ways. We then place the
remaining two points, $c$ and $d$, of $\{P_1,P_2,P_3,P_4\}$
on $C$ in $(q-1)(q-2)$ ways. The line $\gline{c}{d}$ intersects
$\gline{a}{b}$ in a point outside of $C$ so there are
$q-2$ ways to choose $P$ on $\gline{a}{b}$ but away
from this intersection point and $a$ and $b$. Finally,
we place the remaining two points of $\{P_5,P_6,P_7\}$ on $C$
in one of $(q-3)(q-4)$ ways. Multiplying everything together
we obtain
\begin{equation*}
 4 \cdot 3  \cdot 3 \cdot (q^5-q^2)(q+1)q(q-1)(q-2)^2(q-3)(q-4).
\end{equation*}
We now add the two answers above together to get
\begin{equation*}
 N_{1}^{1,2,3,4}= 24 q^3(q-2)(q-3)(q-4)(2q-5)(q+1)(q^2+q+1)(q-1)^2.
\end{equation*}

To compute $N_2^{1,2,3,4}$, we first  choose one of the 
points $P_1$, $P_2$, $P_3$ or $P_4$ to be the point
$P$ not on the smooth conic $C$ and call the remaining three points $P_r$, $P_s$ and $P_t$.
We continue by choosing a smooth conic $C$ in $q^5-q^2$ ways.

We first assume that $P$ lies on two lines  $\gline{R_1}{R_2}$ and $\gline{O_1}{O_2}$ where 
$\{R_1, O_1 \} \subset \{P_5,P_6,P_7\}$ and $\{R_2,O_2\} \subset \{P_r,P_s,P_t\}$.
We now choose two points among $\{P_5,P_6,P_7\}$ in three ways and choose two points 
among $\{P_r,P_s,P_t\}$ in three ways and rename
the remaining two points to $P_u$ and $P_v$.
There are now two possible ways to label the four chosen points $R_1,R_2,O_1$ and $O_2$
in such a way that $\{R_1,O_1\} \subset \{P_5,P_6,P_7\}$ and $\{R_2,O_2\} \subset \{P_r,P_s,P_t\}$
and we choose one of them.
We then place the four points $R_1$, $R_2$, $O_1$ and $O_2$ on the conic $C$
in $(q+1)q(q-1)(q-2)$ ways. The point $P$ is now given as $P=\gline{R_1}{R_2} \cap \gline{O_1}{O_2}$
and no matter how we place $P_u$ and $P_v$, the three lines $\gline{P_r}{P_s}$, 
$\gline{P_r}{P_t}$ and $\gline{P_s}{P_t}$ will not go through $P$.
We can now multiply everything together to obtain
\begin{equation*}
 4 \cdot 3 \cdot 3 \cdot 2 \cdot (q^5-q^2)(q+1)q(q-1)(q-2)(q-3)(q-4).
\end{equation*}

The other possibility is that $P$ lies on two lines $\gline{R_1}{R_2}$ and $\gline{R_3}{b}$ where
$\{R_1,R_2,R_3\}$ is the set $\{P_5,P_6,P_7\}$ and $b \in \{P_r,P_s,P_t\}$.
We thus choose $b$ in three ways and rename the remaining two points
in $\{P_r,P_s,P_t\}$ to $P_u$ and $P_v$.
From now on, we must differentiate between when $P$ is on the outside
and on the inside of $C$.

First, we choose $P$ on the outside of $C$ in $\frac{1}{2}(q+1)q$ ways.
We then choose $b$ as a point on $C$ whose tangent does not
pass through $P$ in $q-1$ ways. We then choose one of the points
$P_5$, $P_6$ and $P_7$ to become the second intersection point in $C \cap \gline{b}{P}$.
Then, we place the remaining two points among $\{P_5,P_6,P_7\}$ on $C$ such that
the line through them passes through $P$ in $q-3$ ways.
There are now $q-5$ ways to choose $P_u$ and $P_v$
such that the line $\gline{P_u}{P_v}$ will pass through $P$.
Thus, the remaining $(q-3)(q-4)-(q-5)=q^2-8q+17$ choices must
give $P_u$ and $P_v$ such that none of the lines
$\gline{P_r}{P_s}$, $\gline{P_r}{P_t}$ and $\gline{P_s}{P_t}$ will
contain $P$. We may now multiply everything together to obtain
\begin{equation*}
 4 \cdot (q^5-q^2) \cdot 3 \cdot \frac{1}{2}(q+1)q \cdot (q-1) \cdot 3 \cdot (q-3) \cdot (q^2-8q+17).
\end{equation*}

Now we choose $P$ on the inside of $C$ in one of $\frac{1}{2}(q-1)q$ ways.
We then choose $b$ as a point on $C$ whose tangent does not
pass through $P$ in $q+1$ ways. We then choose one of the points
$P_5$, $P_6$ and $P_7$ to become the second intersection point in $C \cap \gline{b}{P}$.
Then, we place the remaining two points among $\{P_5,P_6,P_7\}$ on $C$ such that
the line through them passes through $P$ in $q-1$ ways.
There are now $q-3$ ways to choose $P_u$ and $P_v$
such that the line $\gline{P_u}{P_v}$ will pass through $P$.
Thus, the remaining $(q-3)(q-4)-(q-3)=(q-3)(q-5)$ choices must
give $P_u$ and $P_v$ such that none of the lines
$\gline{P_r}{P_s}$, $\gline{P_r}{P_t}$ and $\gline{P_s}{P_t}$ will
contain $P$. We may now multiply everything together to obtain
\begin{equation*}
 4 \cdot (q^5-q^2) \cdot 3 \cdot \frac{1}{2}(q-1)q \cdot (q+1) \cdot 3 \cdot (q-1) \cdot (q-3)(q-5).
\end{equation*}
We may now add everything together to get
\begin{equation*}
 N_2^{1,2,3,4} = 36q^3(q+1)(q^2+q+1)(5q^3-37q^2+82q-60)(q-1)^2.
\end{equation*}

Finally, we need to compute $N_3^{1,2,3,4}$. We begin by choosing one of the 
points $P_1$, $P_2$, $P_3$ or $P_4$ to be the point
$P$ not on the smooth conic $C$ and call the remaining three points $P_r$, $P_s$ and $P_t$.
We continue by choosing a smooth conic $C$ in $q^5-q^2$ ways.

Here, we only have the possibility that $P$ lies
on three lines $\gline{R_1}{O_1}$, $\gline{R_2}{O_2}$
and $\gline{R_3}{O_3}$ where $\{R_1,R_2,R_3\} = \{P_5,P_6,P_7\}$
and $\{O_1,O_2,O_3\}=\{P_r,P_s,P_t\}$. However, we must take
care of the case that $P$ is on the outside of $C$ and
the case that $P$ is on the inside of $C$ separately.
We call the corresponding numbers $N_{3,\mathrm{out}}^{1,2,3,4}$
and $N_{3,\mathrm{in}}^{1,2,3,4}$.

We begin by computing $N_{3,\mathrm{out}}^{1,2,3,4}$. We thus choose
the point $P$ as a point on the outside of $C$ in $\frac{1}{2}(q+1)q$
ways. 
We begin by placing $P_5$ at one of the points of $C$ whose tangent
does not pass through $P$ in $q-1$ ways. We label the second intersection
point of $C \cap \gline{P_5}{P}$ with $P_r$, $P_s$ or $P_t$ and call the
remaining two points $P_u$ and $P_v$.
We then place $P_6$ at one of the $q-3$ remaining points of $C$ whose tangent
does not pass through $P$ and then choose one of the points $P_u$ and $P_v$
to become the other intersection point of $C \cap \gline{P_6}{P}$.
Finally, we place $P_7$ at one of the remaining $q-5$ points
and label the other point of $C \cap \gline{P_7}{P}$ in the only possible
way. We thus have
\begin{equation*}
 N_{3,\mathrm{out}}^{1,2,3,4} = 4 \cdot (q^5-q^2) \cdot \frac{1}{2}(q+1)q \cdot (q-1) \cdot 3 \cdot (q-3) \cdot 2 \cdot (q-5).
\end{equation*}

We now turn to computing $N_{3,in}^{1,2,3,4}$. We thus choose
the point $P$ as a point on the inside of $C$ in $\frac{1}{2}(q-1)q$
ways.
We begin by placing $P_5$ at one of the points of $C$ whose tangent
does not pass through $P$ in $q+1$ ways. We label the second intersection
point of $C \cap \gline{P_5}{P}$ with $P_r$, $P_s$ or $P_t$ and call the
remaining two points $P_u$ and $P_v$.
We then place $P_6$ at one of the $q-1$ remaining points of $C$  
and then choose one of the points $P_u$ and $P_v$
to become the other intersection point of $C \cap \gline{P_6}{P}$.
Finally, we place $P_7$ at one of the remaining $q-3$ points
and label the other point of $C \cap \gline{P_7}{P}$ in the only possible
way.
We now see that
\begin{equation*}
 N_{3,\mathrm{in}}^{1,2,3,4} = 4 \cdot (q^5-q^2) \cdot \frac{1}{2}(q-1)q \cdot  (q+1) \cdot 3 \cdot (q-1) \cdot 2 \cdot (q-3),
\end{equation*}
and we get
\begin{equation*}
 N_3^{1,2,3,4} = 192q^3(q+1)(q^2+q+1)(q^2-3q+3)(q-1)^2.
\end{equation*}
We now obtain
\begin{equation*}
 |\Delta_l \cap \Delta_c| = |\pglt| \cdot (93q^4-1245q^3+6195q^2-13470q+10737),
\end{equation*}
and, finally,
 \begin{equation*}
  \left|\left( \Pts \right)^{\Frob \sigma} \right| = q^6-35q^5+490q^4-3485q^3+13174q^2-24920q+18375.
 \end{equation*}
 
This concludes the equivariant point count of $\Q[2]$. In Section~\ref{summary}
we provide a summary of the results of the computations.

\section{The hyperelliptic locus}

Up to this point we have almost exclusively discussed plane
quartics. We shall now briefly turn our attention to the other type
of genus $3$ curves - the hyperelliptic curves.
There are many possible ways to approach the computation of the cohomology
of $\Hg{3}[2]$. Our choice is by means of equivariant point counts as in the
previous section.

Recall that a hyperelliptic curve $C$ of genus $g$ is determined, up to isomorphism,
by $2g+2$ distinct points on $\mathbb{P}^1$, up to projective equivalence
and that any such collection $S$ of $2g+2$ points determines a double cover 
$\pi: C \to \Pn{1}$ branched precisely over $S$ (and $C$ is thus a hyperelliptic curve).
Moreover, if we pick $2g+2$ ordered points $P_1, \ldots, P_{2g+2}$ on $\Pn{1}$, the curve $C$ 
also attains a level $2$-structure. In the genus $3$ case, we get $8$ points $Q_i=\pi^{-1}(P_i)$ which determine $\binom{8}{2}=28$ odd theta characteristics $Q_i + Q_j$, $i <j$
 and $\{Q_1+Q_8, \ldots,Q_7+Q_8\}$ is an ordered Aronhold basis, see
 \cite{grossharris}
 and \cite{arbarellocgh}, Appendix B.32-33, 
 and an ordered Aronhold basis determines a level $2$-structure.
 
 However, not all level $2$-structures on the  hyperelliptic curve $C$ arise from different orderings of the points.
 Nevertheless, there is an intimate
 relationship between the moduli space $\Hg{g}[2]$ of
 hyperelliptic curves with level $2$-structure and the moduli space 
 $\mathcal{M}_{0,2g+2}$ of $2g+2$ ordered points on $\Pn{1}$
 given by the following theorem
 which can be found in \cite{dolgachevortland}, Theorem VIII.1.

\pagebreak[2]
\begin{thm}
\label{m08irrthm}
 Each irreducible component of $\Hg{g}[2]$ is isomorphic to the
 moduli space $\mathcal{M}_{0,2g+2}$ of $2g+2$ ordered points on the projective line.
\end{thm}
\pagebreak[2]

Dolgachev and Ortland \cite{dolgachevortland} pose the question
whether the irreducible components of $\Hg{g}[2]$ also are the connected
components or, in other words, if $\Hg{g}[2]$ is smooth.
In the complex case, the question was answered positively by Tsuyumine in \cite{tsuyumine}
and later, by a shorter argument, by Runge in \cite{runge}.
Using the results of \cite{andreatta}, the argument of Runge carries over word for
word to an algebraically closed field of positive characteristic different from $2$.

\pagebreak[2]
\begin{thm}
 \label{m08connthm}
 If $g \geq 2$, then each irreducible component of $\Hg{g}[2]$ is
 also a connected component.
\end{thm}
\pagebreak[2]

We have a natural action of $S_{2g+2}$ on the space $\mathcal{M}_{0,2g+2}$. 
Since different orderings of the points correspond
to different symplectic level $2$ structures, $S_{2g+2}$ sits naturally inside
$\symp{2g,\Z/2\Z}$ and, in fact, for $g=3$ and for even $g$
it is a maximal subgroup, see \cite{dye2}. With Theorems \ref{m08irrthm} and \ref{m08connthm} at hand, the following
slight generalization of a corollary in \cite{dolgachevortland} (p.145) is clear.

\pagebreak[2]
\begin{cor}
\label{m08conncor}
 Let $g \geq 2$ and let $X_{[\tau]} = \mathcal{M}_{0,2g+2}$ for each left coset $[\tau] \in \mathscr{T}:=\symp{2g,\Z/2\Z}/S_{2g+2}$.
 Then
 \begin{equation*}
  \Hg{g}[2] \cong \coprod_{[\tau] \in \mathscr{T}} X_{[\tau]},
 \end{equation*}
 and the group $\symp{2g,\Z/2\Z}$ acts transitively on the set of connected
 components $X_{[\tau]}$ of $\Hg{g}[2]$. In particular, there are
 \begin{equation*}
  \frac{\left|\symp{2g,\Z/2\Z}\right|}{\left|S_{2g+2}\right|} = 
  \frac{2^{g^2}\left(2^{2g}-1\right)\left(2^{2g-2}-1\right) \cdots \left(2^{2}-1\right)}{(2g+2)!},
 \end{equation*}
 connected components of $\Hg{g}[2]$.
\end{cor}
\pagebreak[2]

\begin{rem}
As pointed out in \cite{runge}, the argument to prove the corollary stated in \cite{dolgachevortland}
is not quite correct in full generality as it is given there. However, it is enough
to prove the result for $g=3$ and for even $g$, 
and in \cite{runge} it is explained how to obtain the full result.
\end{rem}

Let us now, once and for all, choose a set $T$ of representatives
of \linebreak $\symp{2g,\Z/2\Z}/S_{2g+2}$. 
If we denote the elements of $X_{[\id]}$ by $x$, then any element
in $X_{[\tau]}$ can be written as $\tau x$ for some $x \in X_{[\id]}$.
Let $\alpha$ be any element of $\symp{2g,\Z/2\Z}$. Then
\begin{equation*}
 \alpha \tau  = \tau' \sigma,
\end{equation*}
for some $\sigma \in S_{2g+2}$ and some $\tau' \in T$.  Since the Frobenius commutes with the action
of $\symp{2g,\Z/2\Z}$ we have that
\begin{equation*}
 F \alpha (\tau x) = \tau x,
\end{equation*}
if and only if 
\begin{equation*}
 F (\tau' \sigma x) = \tau' (F \sigma x) = \tau x.
\end{equation*}
But the Frobenius acts on each of the components of $\Hg{g}[2]$ so we
see that $F \alpha (\tau x)=\tau x$ if and only if $\tau' = \tau$ and
$F \sigma x = x$. 

We now translate the above observation into more standard representation theoretic
vocabulary. Define a class function $\psi$ on $S_{2g+2}$ by
\begin{equation*}
 \psi(\sigma) = |X_{\id}^{F \sigma}|,
\end{equation*}
and define a class function $\hat{\psi}$ on $\symp{2g,\Z/2\Z}$
by setting
\begin{equation*}
 \hat{\psi}(\alpha) = |\Hg{g}[2]^{F \alpha}|,
\end{equation*}
for any $\alpha \in \symp{6}{\Z/2\Z}$.
By the above observation we have that
\begin{equation*}
 \hat{\psi}(\alpha) = \sum_{\tau \in T} \widetilde{\psi}(\tau^{-1} \alpha \tau),
\end{equation*}
where 
\begin{equation*}
 \widetilde{\psi}(\beta) = \left\{ \begin{array}{ll}
                                     \psi(\beta) & \text{if } \beta \in S_{2g+2}, \\
                                     0, & \text{otherwise}.
                                    \end{array} \right.
\end{equation*}
In other words, $\hat{\psi}$ is the class function $\psi$ induced from $S_{2g+2}$
up to $\symp{2g,\Z/2\Z}$. Thus, to make an $S_{2g+2}$-equivariant
point count of $\Hg{g}[2]$ we can make an $S_{2g+2}$-equivariant
point count of $\mathcal{M}_{0,2g+2}$ and then use the representation
theory of $S_{2g+2}$ and $\symp{2g,\Z/2\Z}$ in order to
first induce the class function up to $\symp{2g,\Z/2\Z}$
and then restrict it down again to $S_{2g+2}$. Once this is done,
we can obtain the $S_{2g+1}$-equivariant point count by restricting
from $S_{2g+2}$ to $S_{2g+1}$. 

Using Lemma~\ref{conjlem}, the $S_8$-equivariant point
count of $\Hg{3}[2]$ is very easy.
We first compute the number of $\lambda$-tuples of $\mathbb{P}^1$ for 
each partition of $\lambda$ of
$8$ and then divide by $|\mathrm{PGL}(2)|$
in order to obtain $|\mathcal{M}_{0,8}^{\Frob \sigma}|$, where
$\sigma$ is a permutation in $S_8$ of cycle type $\lambda$.
The result is given in Table~\ref{m08equivarianttab}.
Once this is done, we induce up to $\symp{6}{\Z/2\Z}$ in order
to obtain the $\symp{6}{\Z/2\Z}$-equivariant cohomology of $\Hg{3}[2]$.
The results are given in Table~\ref{hyps8eqtab}.
Finally, we
restrict to $S_7$ to get the results of Table~\ref{hypS7tab} and~\ref{Hypcohtable}.
The computations present no difficulties whatsoever.
We also mention that the equivariant Poincar\'e polynomials of
$\mathcal{M}_{0,n}$ and $\overline{\mathcal{M}}_{0,n}$ have been
computed for all $n \geq 3$ in \cite{getzler}.

It is not very hard to see that $M_{0,2g+2}$ is isomorphic to the complement of a
hyperplane arrangement. One way to see this is to start by placing the first three
points at $0$, $1$ and $\infty$. Then $\mathcal{M}_{0,2g+2}$ is isomorphic
to $\left( \mathbb{A}^1\setminus \{0,1\} \right)^{2g-1} \setminus \Delta$,
where $\Delta \subset \left( \mathbb{A}^1\setminus \{0,1\} \right)^{2g-1}$ is the
subset of points where at least two coordinates are equal. Thus,
by the results of Section~\ref{finiteminimalpurity} we can deduce the cohomology
of $\Hg{3}[2]$ from the equivariant point counts.
In Section~\ref{summary}
we provide a summary of the results of the computations.

\section{The total moduli space}
We now know the cohomology groups of both $\Q[2]$ and $\Hg{3}[2]$ as representations of $S_7$.
Unfortunately, we have not been able to obtain the cohomology of $\Mg{3}[2]$.
In order to say something, we shall use the comparison theorem in étale cohomology
and switch to work over the complex numbers and with re Rham cohomology. 
There, we have the following result.

\begin{lem}[Looijenga, \cite{looijenga}]
 \label{gysinlemma}
 Let $X$ be a variety of pure dimension and let $Y \subset X$ be a hypersurface.
 Then there is a Gysin exact sequence of mixed Hodge structures
 \begin{equation*}
  \cdots \to \derham{k-2}{Y}(\text{-}1) \to \derham{k}{X} \to \derham{k}{X\setminus Y} \to \derham{k-1}{Y}(\text{-}1) \to \cdots  
 \end{equation*}
 \end{lem}
 
Looijenga, \cite{looijenga}, also showed that $\derham{k}{\Qbtg}$ is of pure Tate
type $(k,k)$. The space $\Ht{3}$ is
isomorphic to a disjoint union of complements of hyperplane arrangements
and the $k$'th cohomology group of such spaces are known to have
Tate type $(k,k)$, also by a result of Looijenga \cite{looijenga}.
Thus, if we apply Lemma \ref{gysinlemma} to $X=\Mg{3}[2]$, $Y=\Hg{3}[2]$ and 
$X\setminus Y=\Q[2]$
we have that the long exact sequence splits into four term sequences
\begin{equation*}
 0 \to W_k\derham{k}{X} \to \derham{k}{X \setminus Y} \to \derham{k-1}{Y}(\text{-}1) \to W_k\derham{k+1}{X} \to 0,
\end{equation*}
where $W_k\derham{k}{X}$ denotes the weight $(k,k)$ part
of $\derham{k}{X}$. 
Moreover, let $m^k_X(\lambda)$ denote the multiplicity of $s_{\lambda}$
in $\derham{k}{X}$ and let 
$n^k(\lambda)= m^k_{\Q[2]}(\lambda) - m^{k-2}_{\Hg{3}[2]}(\lambda)$.
If $n^k(\lambda)\geq 0$, then $s_{lambda}$ occurs with multiplicity at least 
$n^k(\lambda)$ in $W_k\derham{k}{\Mg{3}[2]}$ and if $n^k(\lambda) \leq 0$, then
$s_{\lambda}$ occurs with multiplicity at least $-n^k(\lambda)$ in
$W_k\derham{k+1}{\Mg{3}[2]}$.
Thus, Tables \ref{Qcohtable} and \ref{Hypcohtable} provide explicit bounds for
the cohomology groups of $\Mg{3}[2]$. 

\section{Summary of computations}
\label{summary}
 We summarize the computations related to $\Q[2]$ in Table \ref{q2equivarianttab} and
 in Proposition~\ref{qtpoincareserre} we give the Poincaré
 polynomial of $\Q[2]$.
  In Table \ref{Qcohtable} we
 give the cohomology of $\Q[2]$ as a representation of $S_7$. The rows correspond
 to the cohomology groups and the columns correspond to the irreducible representations of $S_7$.
 The symbol $s_{\lambda}$ denotes the irreducible representation of $S_7$ corresponding to the partition
 $\lambda$ and a number $n$ in row $H^k$ and column $s_{\lambda}$ means that $s_{\lambda}$ occurs
 in $H^k$ with multiplicity $n$.
 
 \pagebreak[2]
 \begin{prop}
 \label{qtpoincareserre}
 The Poincaré polynomial of $\Q[2]$ is
 \begin{equation*}
 \poincareserre{\Q[2]}(t) =
  1+35t+490t^2+3485t^3+13174t^4+24920t^5+18375t^6.
 \end{equation*}
 \end{prop}
 \pagebreak[2] 
 
 Tables~\ref{m08equivarianttab},~\ref{hyps8eqtab} and~\ref{hypS7tab} give equivariant
 point counts for various spaces and groups related to the equivariant point count of
 $\Hg{3}[2]$ and in Table \ref{Hypcohtable}
 we give the cohomology groups of $\Hg{3}[2]$ as representations of $S_7$.  
 For convenience of we also give the Poincaré polynomial of $\Hg{3}[2]$.
 
 \pagebreak[2]
 \begin{prop}
 \label{hgtpoincareserre} 
 The Poincaré polynomial of $\Hg{3}[2]$ is
 \begin{equation*}
 \poincareserre{\Hg{3}[2]}(t) =
 36+720t+5580t^2+20880t^3+37584t^4+25920t^5.
 \end{equation*}
 \end{prop}
 \pagebreak[2] 
 
 \begin{table}[htbp]
 \begin{align*}
  & \lambda & \, & |\Q[2]^{F \cdot \sigma_{\lambda}}| \\
  \hline\\[-10pt]
  & \left[7\right] & \, & q^6+q^3 \\
  & \left[6,1\right] & \, & q^6-2q^3+1 \\
  & \left[5,2\right] & \, & q^6-q^2 \\
  & \left[5,1^2\right] & \, & q^6-q^2 \\
  & \left[4,3\right] & \, & q^6-q^5-2q^4+q^3+q^2 \\
  & \left[4,2,1\right] & \, & q^6-q^5-2q^4+q^3-2q^2+3 \\
  & \left[4,1^3\right] & \, & q^6-q^5-2q^4+q^3-2q^2+3 \\
  & \left[3^2,1\right] & \, & q^6-2q^5-2q^4-8q^3+16q^2+10q+21 \\
  & \left[3,2^2\right] & \, & q^6-q^5-2q^4+3q^3+q^2-2q \\
  & \left[3,2,1^2\right] & \, & q^6-3q^5+5q^3-q^2-2q \\
  & \left[3,1^4\right] & \, & q^6-5q^5+10q^4-5q^3-11q^2+10q \\
  & \left[2^3,1\right] & \, & q^6-3q^5-6q^4+19q^3+6q^2-24q+7 \\
  & \left[2^2,1^3\right] & \, & q^6-7q^5+10q^4+15q^3-26q^2-8q+15 \\
  & \left[2,1^5\right] & \, & q^6-15q^5+90q^4-265q^3+374q^2-200q+15 \\
  & \left[1^7\right] & \, & q^6-35q^5+490q^4-3485q^3+13174q^2-24920q+18375
\end{align*}
\caption{The $S_7$-equivariant point count of $\Q[2]$. 
We use $\sigma_{\lambda}$ to denote any permutation in $S_7$ of cycle type $\lambda$.}
\label{q2equivarianttab}
\end{table}

\begin{table}[htbp]
\begin{equation*}
\begin{array}{r|rrrrrrrrrr} 
\, & s_{7} & s_{6,1} & s_{5,2} & s_{5,1^2} & s_{4,3} & s_{4,2,1} & s_{4,1^3} & s_{3^2,1} & s_{3,2^2} & s_{3,2,1^2} \\
\hline
H^0 & 1 & 0 & 0 & 0 & 0 & 0 & 0 & 0 & 0 & 0 \\ 
H^1 & 1 & 1 & 1 & 0 & 1 & 0 & 0 & 0 & 0 & 0 \\ 
H^2 & 0 & 3 & 4 & 4 & 3 & 5 & 1 & 3 & 1 & 1 \\
H^3 & 1 & 8 & 14 & 18 & 14 & 30 & 16 & 16 & 12 & 18 \\ 
H^4 & 4 & 20 & 44 & 47 & 44 & 99 & 56 & 56 & 54 & 83 \\ 
H^5 & 6 & 33 & 76 & 76 & 72 & 178 & 97 & 104 & 105 & 169 \\ 
H^6 & 6 & 23 & 51 & 54 & 54 & 127 & 74 & 76 & 77 & 126 \\
\hline
\, & s_{3,1^4} & s_{2^3,1} & s_{2^2,1^3} & s_{2,1^5} & s_{1^7} & \, & \, & \, & \, & \, \\
\hline
H^0 & 0 & 0 & 0 & 0 & 0 & \,&\,&\,&\,&\,\\ 
H^1 & 0 & 0 & 0 & 0 & 0 & \,&\,&\,&\,&\,\\ 
H^2 & 0 & 0 & 0 & 0 & 0 & \,&\,&\,&\,&\,\\
H^3 & 4 & 6 & 3 & 0 & 0 & \,&\,&\,&\,&\,\\ 
H^4 & 32 & 31 & 25 & 6 & 1 & \,&\,&\,&\,&\,\\ 
H^5 & 71 & 65 & 64 & 26 & 3 & \,&\,&\,&\,&\,\\ 
H^6 & 54 & 54 & 50 & 22 & 5 & \,&\,&\,&\,&\,
\end{array}
\end{equation*}
\caption{The cohomology of $\Q[2]$ as a representation of $S_7$.}
\label{Qcohtable}
\end{table}

  \begin{table}[htbp]
 \begin{align*} 
  & \lambda & \, & |\mathcal{M}_{0,8}^{F \cdot \sigma_{\lambda}}|\\
  \hline \\[-10pt]
  & \left[8\right] & \, & \left( {q}^{2}+1 \right) {q}^{3} \\
  & \left[7,1\right] & \, & \left( q+1 \right)  \left( {q}^{2}+q+1 \right)  \left( {q}^{2}-q+1 \right) \\
  & \left[6,2\right] & \, & q \left( q-1 \right) \left( {q}^{3}+q-1 \right) \\
  & \left[6,1^2\right] & \, & q \left( q+1 \right)  \left( {q}^{3}+q-1 \right) \\
  & \left[5,3\right] & \, &  q \left( q-1 \right) \left( q+1 \right)  \left( {q}^{2}+1 \right) \\
  & \left[5,2,1\right] & \, & q \left( q-1 \right) \left( q+1 \right)  \left( {q}^{2}+1 \right) \\
  & \left[5,1^3\right] & \, & q \left( q-1 \right) \left( q+1 \right)  \left( {q}^{2}+1 \right) \\
  & \left[4^2\right] & \, & q \left( {q}^{4}-{q}^{2}-4 \right) \\
  & \left[4,3,1\right] & \, & \left( q-1 \right) {q}^{2} \left( q+1 \right) ^{2} \\
  & \left[4,2^2\right] & \, & \left( q-1 \right)  \left( q-2 \right) \left( q+1 \right) {q}^{2}\\
  & \left[4,2,1^2\right] & \, & \left( q-1 \right)  \left( q+1 \right) {q}^{3} \\
  & \left[4,1^4\right] & \, & \left( q-1 \right)  \left( q-2 \right) \left( q+1 \right) {q}^{2} \\
  & \left[3^2,2\right] & \, & q \left( q-1 \right) \left( {q}^{3}-q-3 \right) \\
  & \left[3^2,1^2\right] & \, & q \left( q+1 \right) \left( {q}^{3}-q-3 \right) \\
  & \left[3,2^2,1\right] & \, & q \left( q-1 \right) \left( q-2 \right) \left( q+1 \right)^{2} \\
  & \left[3,2,1^3 \right] & \, & \left( q+1 \right)  q^2 \left( q-1 \right)^{2} \\
  & \left[3,1^5\right] & \, & q \left( q-1 \right)  \left( q-2 \right)  \left( q-3 \right) \left( q+1 \right) \\
  & \left[2^4\right] & \, & \left( q-2 \right) \left( q-3 \right) \left( q+2 \right) \left( {q}^{2}-q-4 \right) \\
  & \left[2^3,1^2\right] & \, & q \left( q-2 \right) \left( q+1 \right) \left( {q}^{2}-q-4 \right) \\
  & \left[2^2,1^4\right] & \, & q \left( q-1 \right)  \left( q+1 \right)  \left( q-2 \right) ^{2} \\
  & \left[2,1^6\right] & \, & q \left( q-1 \right)  \left( q-2 \right)  \left( q-3 \right)  \left( q-4 \right) \\
  & \left[1^8\right] & \, & \left( q-2 \right)  \left( q-3 \right)  \left( q-4 \right)  \left( q-5 \right)  \left( q-6 \right)
\end{align*}
\caption{The $S_8$-equivariant point count of $\mathcal{M}_{0,8}$. 
We use $\sigma_{\lambda}$ to denote any permutation in $S_8$ of cycle type $\lambda$.}
\label{m08equivarianttab}
\end{table}

 \begin{table}[htbp]
 \begin{align*} 
  & \lambda & \, & |\mathcal{H}_3[2]^{F \cdot \sigma_{\lambda}}|\\
  \hline \\[-10pt]
  & \left[8\right] & \, &2\,{q}^{5}+2\,{q}^{3} \\
  & \left[7,1\right] & \, & {q}^{5}+{q}^{4}+{q}^{3}+{q}^{2}+q+1 \\
  & \left[6,2\right] & \, & 3\,{q}^{5}+3\,{q}^{3}-6\,{q}^{2}-3\,{q}^{4}+3\,q \\
  & \left[6,1^2\right] & \, & {q}^{5}+{q}^{4}+{q}^{3}-q \\
  & \left[5,3\right] & \, & {q}^{5}-q \\
  & \left[5,2,1\right] & \, & {q}^{5}-q \\
  & \left[5,1^3\right] & \, & {q}^{5}-q \\
  & \left[4^2\right] & \, & 4\,{q}^{5}-16\,q-4\,{q}^{3} \\
  & \left[4,3,1\right] & \, & 2\,{q}^{5}+2\,{q}^{4}-2\,{q}^{3}-2\,{q}^{2} \\
  & \left[4,2^2\right] & \, & 6\,{q}^{5}+12\,{q}^{2}-12\,{q}^{4}-6\,{q}^{3} \\
  & \left[4,2,1^2\right] & \, & 2\,{q}^{5}-2\,{q}^{3} \\
  & \left[4,1^4\right] & \, & 2\,{q}^{5}-4\,{q}^{4}-2\,{q}^{3}+4\,{q}^{2} \\
  & \left[3^2,2\right] & \, & {q}^{5}-{q}^{4}-{q}^{3}-2\,{q}^{2}+3\,q \\
  & \left[3^2,1^2\right] & \, & 3\,{q}^{5}+3\,{q}^{4}-3\,{q}^{3}-12\,{q}^{2}-9\,q \\
  & \left[3,2^2,1\right] & \, & 2\,{q}^{5}-2\,{q}^{4}-6\,{q}^{3}+2\,{q}^{2}+4\,q \\
  & \left[3,2,1^3 \right] & \, & 4\,{q}^{5}-4\,{q}^{4}-4\,{q}^{3}+4\,{q}^{2} \\
  & \left[3,1^5\right] & \, & 6\,{q}^{5}-30\,{q}^{4}+30\,{q}^{3}+30\,{q}^{2}-36\,q \\
  & \left[2^4\right] & \, & 12\,{q}^{5}+48\,q-60\,{q}^{3}+336\,{q}^{2}-48\,{q}^{4}-576 \\
  & \left[2^3,1^2\right] & \, & 4\,{q}^{5}-8\,{q}^{4}-20\,{q}^{3}+24\,{q}^{2}+32\,q \\
  & \left[2^2,1^4\right] & \, & 8\,{q}^{5}-32\,{q}^{4}+24\,{q}^{3}+32\,{q}^{2}-32\,q \\
  & \left[2,1^6\right] & \, & 16\,{q}^{5}-160\,{q}^{4}+560\,{q}^{3}-800\,{q}^{2}+384\,q \\
  & \left[1^8\right] & \, & 36\,{q}^{5}-720\,{q}^{4}+5580\,{q}^{3}-20880\,{q}^{2}+37584\,q-25920
\end{align*}
\caption{The $S_8$-equivariant point count of $\mathcal{H}_3[2]$. 
We use $\sigma_{\lambda}$ to denote any permutation in $S_8$ of cycle type $\lambda$.}
\label{hyps8eqtab}
\end{table}

 \begin{table}[htbp]
 \begin{align*} 
  & \lambda & \, & |\mathcal{H}_3[2]^{F \cdot \sigma_{\lambda}}| \\
  \hline \\[-10pt]
  & \left[7\right] & \, & {q}^{5}+{q}^{4}+{q}^{3}+{q}^{2}+q+1 \\
  & \left[6,1\right] & \, & {q}^{5}+{q}^{4}+{q}^{3}-q \\
  & \left[5,2\right] & \, & {q}^{5}-q \\
  & \left[5,1^2\right] & \, & {q}^{5}-q \\
  & \left[4,3\right] & \, & 2\,{q}^{5}+2\,{q}^{4}-2\,{q}^{3}-2\,{q}^{2} \\
  & \left[4,2,1\right] & \, & 2\,{q}^{5}-2\,{q}^{3} \\
  & \left[4,1^3\right] & \, & 2\,{q}^{5}-4\,{q}^{4}-2\,{q}^{3}+4\,{q}^{2} \\
  & \left[3^2,1\right] & \, & 3\,{q}^{5}+3\,{q}^{4}-3\,{q}^{3}-12\,{q}^{2}-9\,q \\
  & \left[3,2^2\right] & \, & 2\,{q}^{5}-2\,{q}^{4}-6\,{q}^{3}+2\,{q}^{2}+4\,q \\
  & \left[3,2,1^2\right] & \, & 4\,{q}^{5}-4\,{q}^{4}-4\,{q}^{3}+4\,{q}^{2} \\
  & \left[3,1^4\right] & \, & 6\,{q}^{5}-30\,{q}^{4}+30\,{q}^{3}+30\,{q}^{2}-36\,q \\
  & \left[2^3,1\right] & \, & 4\,{q}^{5}-8\,{q}^{4}-20\,{q}^{3}+24\,{q}^{2}+32\,q \\
  & \left[2^2,1^3\right] & \, & 8\,{q}^{5}-32\,{q}^{4}+24\,{q}^{3}+32\,{q}^{2}-32\,q \\
  & \left[2,1^5\right] & \, & 16\,{q}^{5}-160\,{q}^{4}+560\,{q}^{3}-800\,{q}^{2}+384\,q \\
  & \left[1^7\right] & \, & 36\,{q}^{5}-720\,{q}^{4}+5580\,{q}^{3}-20880\,{q}^{2}+37584\,q-25920
\end{align*}
\caption{The $S_7$-equivariant point count of $\mathcal{H}_3[2]$. 
We use $\sigma_{\lambda}$ to denote any permutation in $S_7$ of cycle type $\lambda$.}
\label{hypS7tab}
\end{table}

\begin{table}[htbp]
\begin{equation*}
\begin{array}{r|rrrrrrrrrr} 
\, & s_{7} & s_{6,1} & s_{5,2} & s_{5,1^2} & s_{4,3} & s_{4,2,1} & s_{4,1^3} & s_{3^2,1} & s_{3,2^2} & s_{3,2,1^2} \\
\hline
H^0 & 2&1&1&0&1&0&0&0&0&0\\ 
H^1 & 2&7&9&5&5&7&1&3&2&1\\ 
H^2 & 3&18&30&31&25&50&20&26&19&26\\ 
H^3 & 6&35&74&80&72&162&86&92&83&129\\
H^4 & 8&48&114&117&109&271&150&157&158&254\\
H^5 & 5&31&72&77&72&180&103&108&108&180\\
\hline
\, & s_{3,1^4} & s_{2^3,1} & s_{2^2,1^3} & s_{2,1^5} & s_{1^7} & \, & \, & \, & \, & \, \\
\hline
H^0 & 0&0&0&0&0& \,&\,&\,&\,&\,\\ 
H^1 &0&0&0&0&0& \,&\,&\,&\,&\,\\ 
H^2 &5&7&4&0&0& \,&\,&\,&\,&\,\\ 
H^3 &43&45&36&10&1& \,&\,&\,&\,&\,\\ 
H^4 &105&96&92&35&4& \,&\,&\,&\,&\,\\ 
H^5 &77&72&72&31&5& \,&\,&\,&\,&\,
\end{array}
\end{equation*}
\caption{The cohomology of $\Hg{3}[2]$ as a representation of $S_7$.}
\label{Hypcohtable}
\end{table}


\bibliographystyle{acm}

\renewcommand{\bibname}{References} 

\clearpage
\bibliography{references}

\end{document}